\documentclass[11pt,twoside]{article}

\usepackage{epsfig,amsfonts,color}
\usepackage{amsmath}
\usepackage{mathtools}

\usepackage{amssymb, palatino, geometry,url}
\usepackage{amsthm}
\usepackage{algorithmic}
\usepackage{booktabs}
\usepackage{siunitx}
\usepackage{xcolor}
\usepackage{url}
\usepackage[most]{tcolorbox}
\usepackage[T1]{fontenc}
\newtcolorbox{highlighted}{colback=cyan,coltext=black,breakable}
\newcommand\hl[1]{\textcolor{black}{#1}}
\newcommand\hlp[1]{\textcolor{black}{#1}}
\usepackage[noresetcount,lined,boxed]{algorithm2e} 

\usepackage[colorlinks=true,linkcolor=blue,citecolor=blue,urlcolor=blue]{hyperref}
\usepackage{subcaption}
\usepackage{cleveref}
\usepackage{fancyhdr}
\newtheorem{theorem}{Theorem}
\theoremstyle{definition}
\newtheorem{definition}{Definition}[section]

\newcommand{\R}{\mbox{I}\!\mbox{R}}

\newcommand{\amin}{\mathop{\mbox{argmin}}}

\newcommand{\be}{\begin{equation}}
\newcommand{\ee}{\end{equation}}
\newcommand{\bea}{\begin{eqnarray}}
\newcommand{\eea}{\end{eqnarray}}

\newcommand{\bvec}{\left(\begin{array}{c}}
\newcommand{\evec}{\end{array}\right)}
\newcommand{\bsub}{\begin{subequations}}
\newcommand{\esub}{\end{subequations}}

\usepackage[thinc]{esdiff}

\usepackage{lineno}

\usepackage[outdir=./Output]{epstopdf}

\begin{document}

\title{Topological Data Analysis: \\ {\Large Concepts, Computation, and Applications in Chemical Engineering}}

\author{Alexander D. Smith${}^{\P}$, Pawe\l{} D\l{}otko${}^\dag$, and Victor M. Zavala${}^{\P}$\thanks{Corresponding Author: victor.zavala@wisc.edu}\\
{\small ${}^{\P}$Department of Chemical and Biological Engineering}\\
{\small \;University of Wisconsin-Madison, 1415 Engineering Dr, Madison, WI 53706, USA}\\
{\small ${}^{\dag}$Diouscuri Center for Topological Data Analysis}\\
{\small \;Institute of Mathematics, Polish Academy of Sciences, Śniadeckich 8, 00-656 Warsaw, Poland}}

\date{}
\maketitle

\begin{abstract}
A primary hypothesis that drives scientific and engineering studies is that data has structure. The dominant paradigms for describing such structure are statistics (e.g., moments, correlation functions) and signal processing (e.g., convolutional neural nets, Fourier series). Topological Data Analysis (TDA) is a field of mathematics that analyzes data from a fundamentally different per-spective. TDA represents datasets as geometric objects and provides dimensionality reduction techniques that project such objects onto low-dimensional \hl{descriptors in 2D. The key properties of these descriptors (also known as topological features) are that they provide multiscale information and that they are stable under perturbations (e.g., noise, translation, and rotation)}. In this work, we review the key mathematical concepts and methods of TDA and present different applications in chemical engineering.

\end{abstract}

{\bf Keywords:} data, topology, space, time, geometry.

\section{Introduction}
\begin{figure}[!htp]
\begin{subfigure}{.5\textwidth}
  \centering
  \includegraphics[width=.7\linewidth]{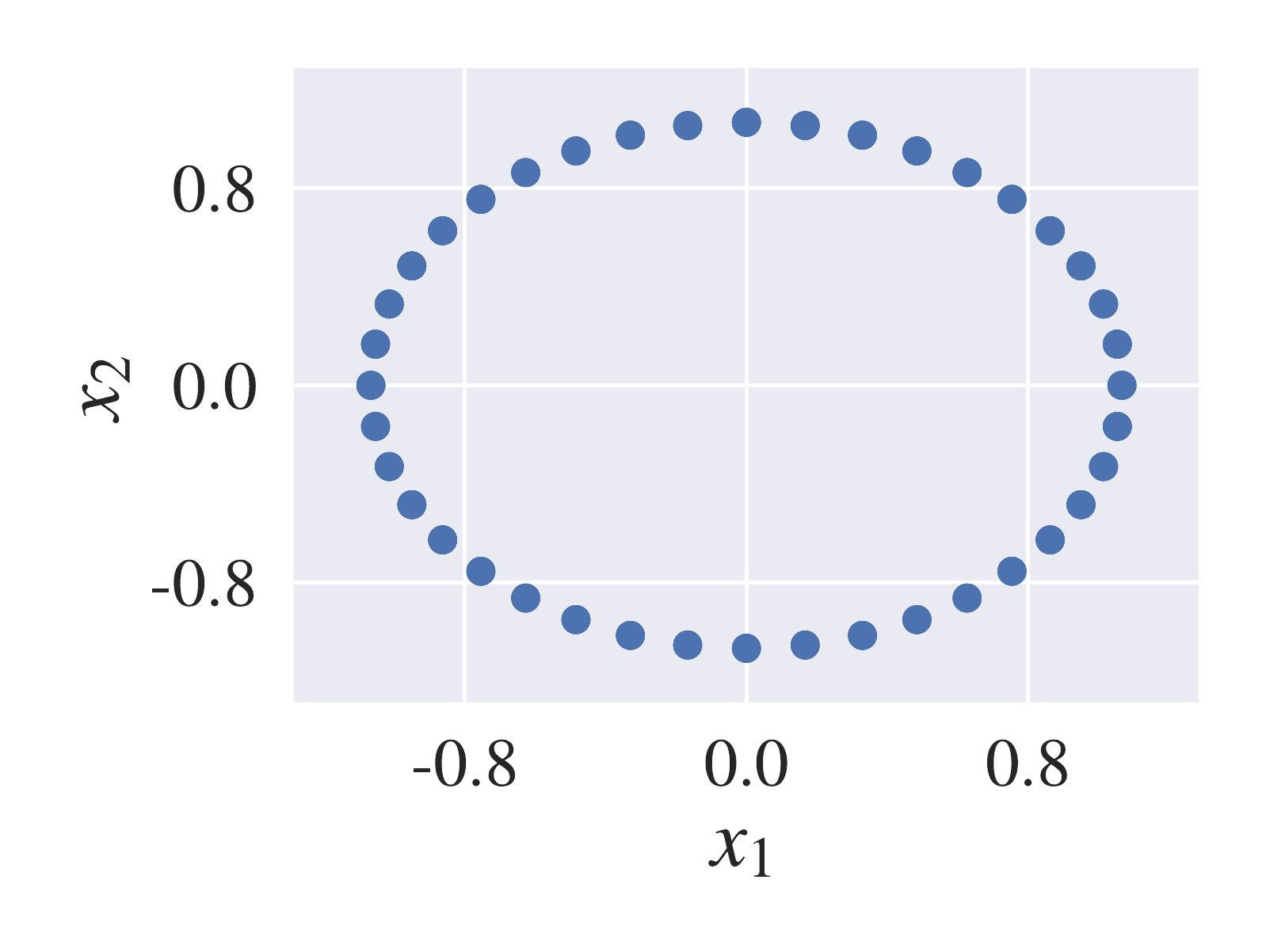}  
  \label{fig:sub-first}
\end{subfigure}
\begin{subfigure}{.5\textwidth}
  \centering
  \includegraphics[width=.7\linewidth]{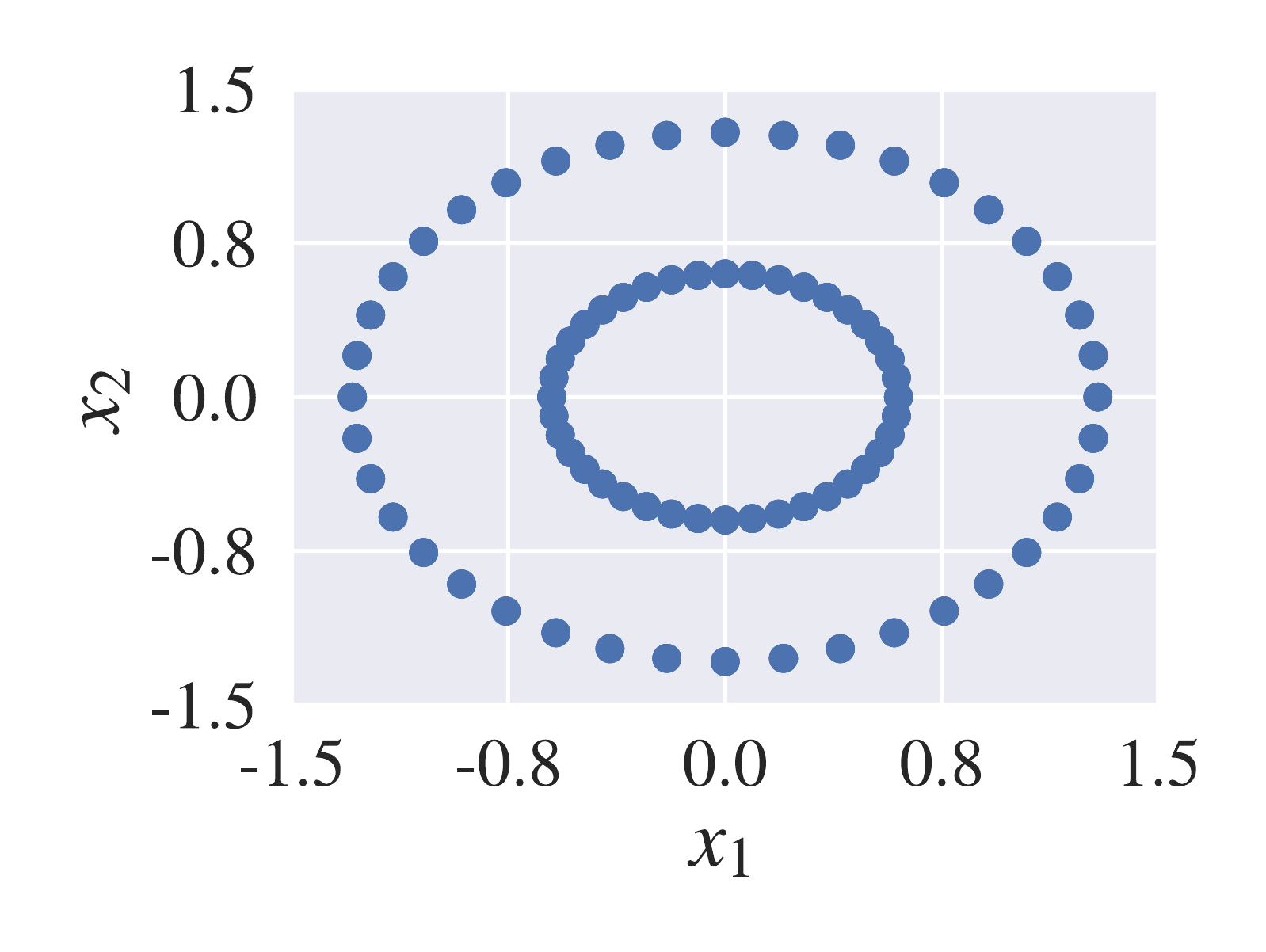}  
  \label{fig:sub-second}
\end{subfigure}
\caption{Datasets with the same mean along the $x_1,x_2$ dimensions (zero), the same standard deviation along both dimensions ($1/2$) and the same correlation between dimensions (zero). While the statistical descriptor values  (i.e., first and second moments of a 2D Gaussian ellipse) are identical, the geometric objects that they define are different.}
\label{fig:samestat}
\end{figure}

Statistical and signal processing techniques are the dominant paradigms used to analyze data. Unfortunately, these techniques provide limited capabilities to analyze certain types of datasets. \hl{Some} interesting examples that highlight this limitation are the Anscombe quartet and the \textit{Datasaurus dozen} datasets \cite{anscombe1973graphs,matejka2017same}. These datasets are {\em visually} distinct but they have the same descriptive statistics (e.g., \textit{mean, standard deviation} and \textit{correlation}). An illustration of this issue is provided in Figure \ref{fig:samestat}; here, the two datasets have the same mean and standard deviation along both dimensions and have the same correlation between dimensions. However, it is clear that these datasets define objects with different geometric features (shape). 

The recent application of algebraic and computational topology to data science has led to the development of a new field known as Topological Data Analysis (TDA) \cite{carlsson2009topology}. TDA techniques are based on the observation that data (e.g., a set of points in a Euclidean space) can be interpreted as elements of a geometric object. As the name suggests, TDA utilizes techniques from computational topology to quantify the \hl{shape} of data \cite{edelsbrunner2010computational}. Fundamentally, topology studies geometric and spatial relations that persist (are stable) in the face of continuous deformations of an object (e.g., stretching, twisting, and bending).This perspective brings a number of advantages over other data analysis techniques \cite{carlsson2009topology,zomorodian2012topological}:

\begin{itemize}
    \item Topology studies data in a manner that is independent of the chosen coordinates. 
    \item Topology studies data in a way that minimizes sensitivity to the choice of metric. 
     \item \hlp{Topology generalizes known graph theory techniques to high-dimensional spaces}.
     \item Topology is robust to large quantities of noise.
\end{itemize}

The main focus of this paper is a technique in the field of TDA known as  \textit{persistence homology} \cite{ghrist2008barcodes,carlsson2005persistence}. The goal of persistent homology is to {\em identify and quantify} topologically dominant features within the data in the form of basic (low-dimensional) topological features such as connected components, holes, voids, and their generalizations. This information can then be used by statistical and machine learning techniques to perform regression, classification, hypothesis testing, and clustering tasks   \cite{bubenik2012statistical,bubenik2010statistical,bubenik2017persistence,blumberg2014robust,adams2017persistence,reininghaus2015stable}.   The TDA methodology is summarized in Figure \ref{fig:summary}. It is important to emphasize that TDA is a dimensionality reduction technique that maps data from its original high-dimensional space to a low-dimensional space that it is easier to understand and visualize. This is similar in spirit to principal component analysis (PCA), which is a statistical technique that projects data into a \hlp{low dimensional space} by extracting latent variables (principal components) that contain maximum information in terms of variance. 

TDA has been used in different scientific and engineering domains. In the medical imaging field, persistent homology has been used in studying brain dendrograms \cite{lee2012persistent} and in identifying brain networks in children with ADHD and autism \cite{lee2011discriminative}. In material science, these techniques have been used to characterize complex craze formations \cite{ichinomiya2017persistent}, to analyze  hierarchical structures in glasses \cite{nakamura2015persistent}, and in materials informatics \cite{buchet2018persistent}. \hlp{They have also been used in high throughput screening of nanoporous materials, such as zeolites \cite{lee2018high}}. TDA has been used in the analysis of dynamical systems and time series  \cite{seversky2016time,perea2015sliding} to study gene expression \cite{perea2015sw1pers} and the dynamics of Kolmogorov flows and Raleigh-Bernard convection \cite{kramar2016analysis}.  TDA has also been used in the analysis of time-varying functional networks \cite{stolz2017persistent}. \hlp{In chemistry and biochemistry, persistent homology has been used to characterize protein structure, flexibility, and folding \cite{xia2014persistent}; as well as a metric for understanding membrane fusion \cite{kasson2007persistent}, and used to predict fullerene stability \cite{xia2015persistent}.}

In this work, we provide a concise summary of relevant concepts and computational methods of TDA from the perspective of chemical engineering applications. We show how to apply persistent homology techniques to analyze datasets described by point clouds and functions in high dimensions and we discuss fundamental stability results of topological features in the face of perturbations. We present multiple case studies with complex synthetic and real experimental datasets to demonstrate the benefits of TDA. Specifically, we show that TDA extracts informative features from complex datasets that correlate strongly with emerging features of practical interest. For instance, we show that the topological features of a 3D solvent environment explains reactivity in such an environment and that topological features of liquid crystals explain composition of its environment. These two results are presented for the first time in this paper. Moreover, we show how to characterize topological features of scatter fields from flow cytometry experiments.  Our work seeks to open new research directions and applications of TDA in chemical engineering. 

\section{TDA Basics}

\begin{figure}[!ht]
    \centering
    \includegraphics[width = 1\textwidth]{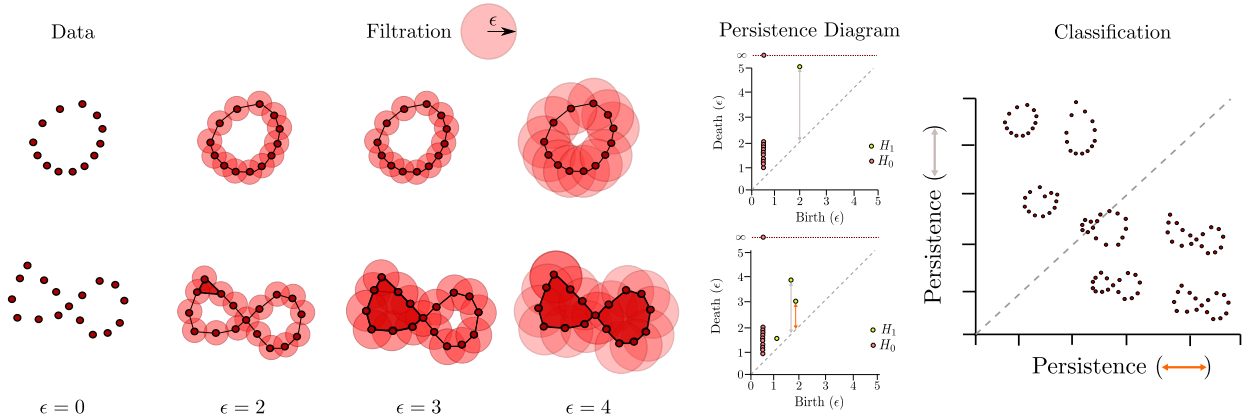}
    \caption{\hl{Persistence homology methodology for point clouds: each point cloud is converted into a geometric object via a filtration where the topology is measured at each point in the filtration. At certain points in the filtration topological features, such as the holes above, appear and are eventually filled. The $\epsilon$ value at the appearance and disappearance of these features are recorded as birth and death in the filtration. The birth and death of the topological features are represented as points in a persistence diagram, with $x$=birth and $y$=death and persistence defined as ($y - x$). The persistence diagram encodes the topological evolution of the data during the filtration and can be used directly to separate point clouds of different shape and cluster those of similar shape. In this illustration we create a representative classification plot that demonstrates the separation of example point clouds based upon the persistence of the largest and second largest (which is zero in some cases) hole(s) that appear and disappear during the filtration.}}
    \label{fig:summary}
\end{figure}

\hl{
This document is partitioned into two main parts, the first develops the mathematical basis for TDA (see Sections \ref{seq:mathbegin}-\ref{seq:mathend}) and the second explores the application of TDA to different chemical engineering focused problems (see Section \ref{section:examples}).  We develop this section as a simplified overview of TDA in order to equip the reader with the knowledge necessary to immediately review TDA applications of interest. The focus of TDA is to capture and record the evolution of the topology of a dataset at different scales which is measured through a filtration (see Section \ref{seq:mathend}). Figure \ref{fig:summary} is an example of applying TDA to two simple point clouds. On point clouds, a filtration is done by expanding balls of radius ($\epsilon$) around each data point, and connecting those points for which the given balls overlap. This changes the topology of the data representation as the \hlp{expansion} proceeds, resulting in the appearance and disappearance of topological features such as connected components (represented by $H_0$) and holes (represented by $H_1$). A similar filtration can also be constructed on continous functions (see Section \ref{sec:contfunc}) or on discrete representations of continuous functions such as images (see Section \ref{sec:cubicalcomp}). Regardless of how the filtration is performed, the time of appearance (birth) and the time of disappearance (death) of the topological features during the filtration constructs the persistence diagram. The persistence diagram is a scatter plot for which the $x$ axis is birth and the $y$ axis is death. All topological features that appear and disappear during the filtration are recorded as points with coordinates $x$ = birth and $y$ = death and persistence ($y-x$) within their representative group (e.g., $H_0$, $H_1$). This scatter plot encodes (reduces) the high dimensional structure of the dataset during the filtration and can be vectorized (see Section \ref{sec:pds}) for use in statistical and machine learning models and methods as we demonstrate in Section \ref{section:examples}.}

\section{Fundamental Concepts of TDA}
\label{seq:mathbegin}

We discuss fundamental concepts and computational methods of TDA. We first introduce the notion of simplicial and cubical complexes, which are the basic geometric constructs used to represent data objects. The representation of data objects as a simplicial or cubical complex enables the use of methods from \textit{simplicial and cubical homology} to quantify the shape of the data \hlp{in terms of connectedness and topologically important features} \cite{hatcher2005algebraic,kan1955abstract}.  In the following discussion, we use $\mathbb{R}$ to denote the set of real numbers and $\mathbb{Z}$ to denote the set of integer numbers..

\begin{figure}[!htp]
    \centering
    \includegraphics[width = .8\textwidth]{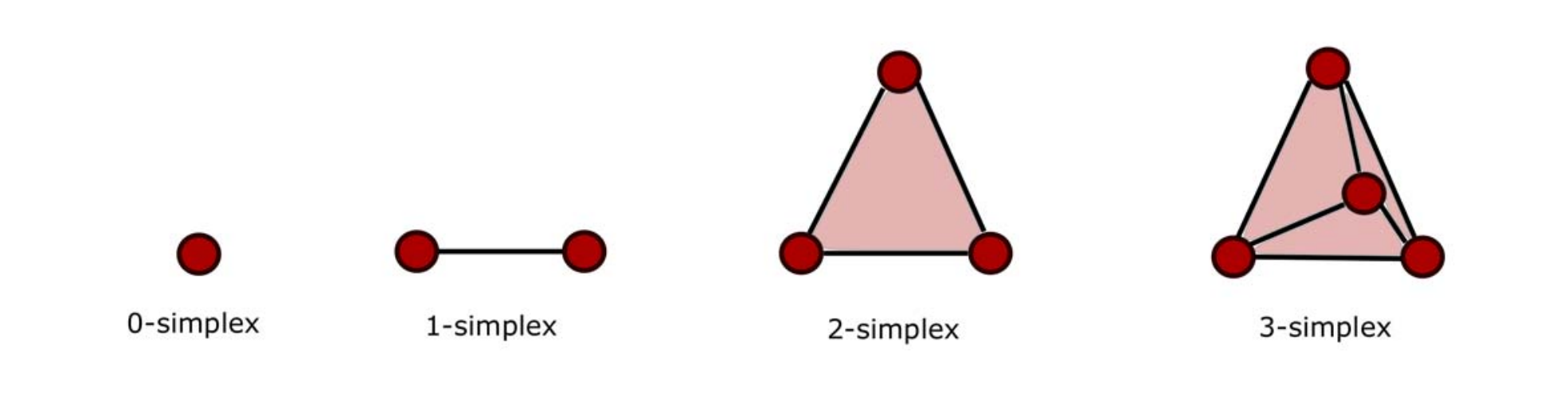}
    \caption{\hl{Examples of $k$-dimensional simplices for $k=0,1,2,3$. A simplex is a generalization of a triangle in high dimensions. $0$ simplices are vertices (points), $1$-simplices are edges, $2$-simplices are triangles, and $3$-simplices are tetrahedra.}}
    \label{fig:simplex}
\end{figure}

\subsection{Simplicies and Simplicial Complexes}

A simplex is a generalization of a triangle from 2D to other dimensions (e.g., a tetrahedra is a 3-dimensional simplex). \hl{Simplices} spanning dimension $k=0$ to dimension $k=3$ are shown in Figure \ref{fig:simplex}.  The formal definition of a simplex is as follows.

\theoremstyle{definition}
\begin{definition}{\textbf{$\mathbf{k}$-simplex}: A $k$-\textit{simplex} is a convex hull spanned by $k$+1 affinely independent points $v\in \mathbb{R}^m$ and is denoted as:}

\begin{equation}
    \sigma = [v_0,v_2,...,v_{k}]
\end{equation}

\end{definition}

Some interesting properties of \hl{simplices} are:
\begin{enumerate}
    \item \hl{An {$m$-face} of simplex $\sigma$ is the convex hull of any of its nonempty subsets.}
    \item The $m$-face is a simplex.
    \item A {$0$-face} is a \textbf{vertex}.
    \item A {$1$-face} is an \textbf{edge}.
\end{enumerate}

\begin{figure}[!htp]
\begin{subfigure}{.5\textwidth}
  \centering
  \includegraphics[width=.5\linewidth]{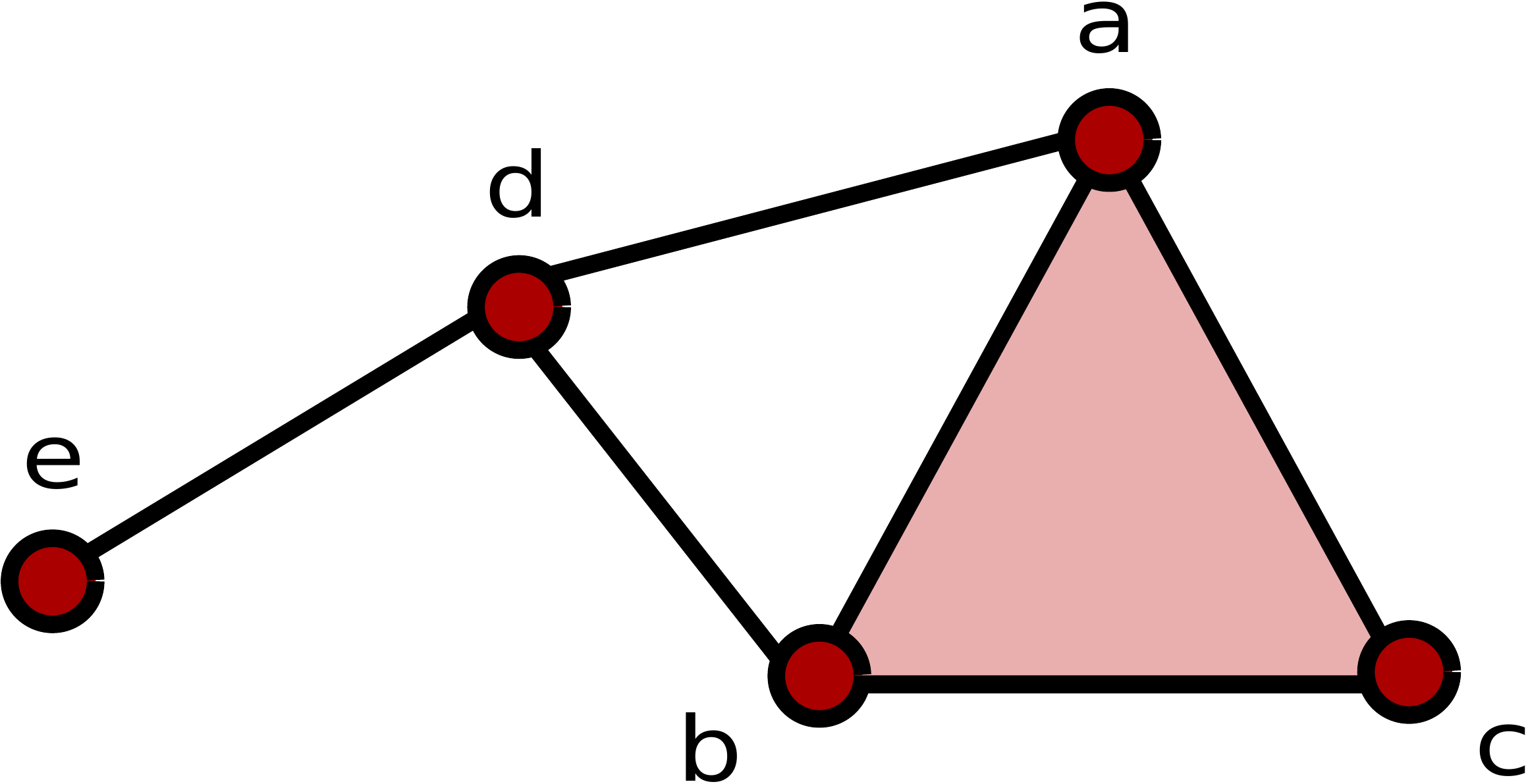}  
  \caption{}
  \label{fig:simplicial_complex}
\end{subfigure}
\begin{subfigure}{.5\textwidth}
  \centering
  \includegraphics[width=.6\linewidth]{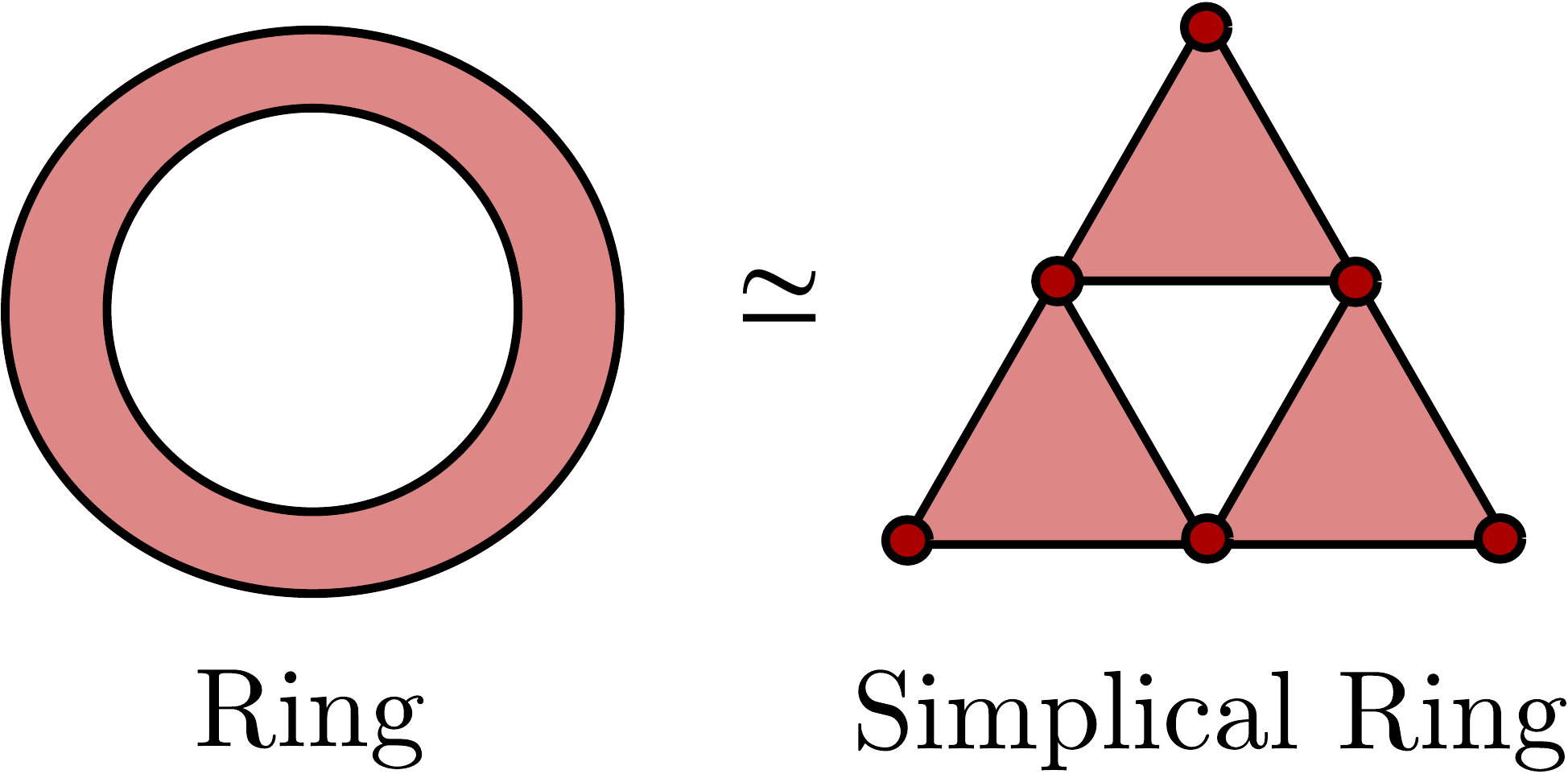}  
  \caption{}
  \label{fig:Annulus}
\end{subfigure}
\caption{(a) A simplicial $2$-complex created by connecting a $2$-simplex (a triangle) and multiple $1$-simplices (edges). This simplicial complex \hl{contains one hole}. (b) A geometric object (a ring) and its simplicial complex representation using $2$-\hl{simplices}. Both shapes contain an empty hole and are homotopy equivalent.}
\label{fig:loopholeex}
\end{figure}

A \textit{simplicial complex} (denoted as $\mathcal{K}$) is obtained by connecting (glueing) \hl{simplices}, as shown in Figure \ref{fig:simplicial_complex}. We denote the dimension of a given simplex or simplicial complex as dim$(\cdot)$.

\theoremstyle{definition}
\begin{definition}{\textbf{Simplicial Complex}:} A simplicial complex $\mathcal{K} \subset R^m$ is a collection of simplices that satisfies the following properties:
\begin{enumerate}
    \item Every {face} of an elemnt of $\mathcal{K}$ is also in $\mathcal{K}$
    \item A nonempty intersection of simplices $\sigma_1,\sigma_2 \in \mathcal{K}$ is a {face} of $\sigma_1$ and $\sigma_2$
    \item The {dimension} of $\mathcal{K}$ is the highest dimension of its simplices: $\dim(\mathcal{K})=\text{max}(\dim(\sigma): \sigma \in \mathcal{K})$
\end{enumerate}
\end{definition}

\hl{A simplicial complex is used to represent the topological characteristics of objects}. \hl{One example of this is found in finite element analysis where  simplicial complexes (also known as triangulations) are used to represent domains over which partial differential equations are solved \cite{cook2007concepts}}. In \hl{Figure  \ref{fig:Annulus}} we represent a geometric object (a ring) as a simplicial complex. We see that the central hole of the ring, which is its main topological feature, is preserved in its representation as the simplicial complex. These objects are thus said to be homotopy equivalent (represented by the notation $\simeq$). Homotopy equivalence identifies spaces which can be deformed continuously into one another without cutting or tearing (e.g., via stretching) \cite{munkres2018elements}. The flexibility of simplicial complexes allows us to create a homotopically equivalent representation of any shape encountered in practice. As we will see, algebraic calculations can be applied to simplicial complexes to quantify the features of the original object \cite{hatcher2005algebraic}. 

\begin{figure}[!htp]
\begin{subfigure}{.5\textwidth}
  \centering
  \includegraphics[width=.4\linewidth]{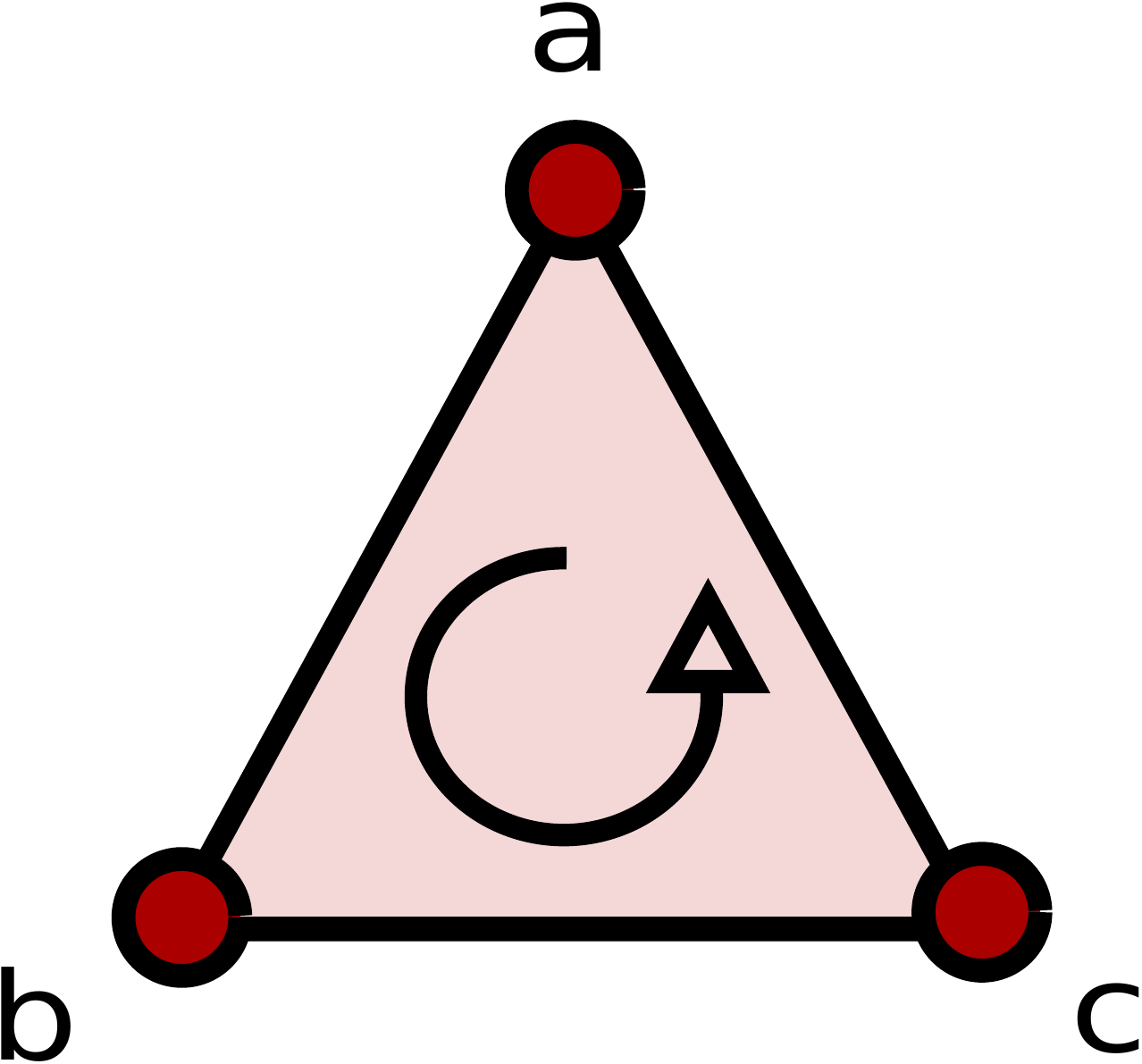}  
  \caption{Simplex represented by ordered set $\{a,b,c\}$}
  \label{fig:sub-first}
\end{subfigure}
\begin{subfigure}{.5\textwidth}
  \centering
  \includegraphics[width=.4\linewidth]{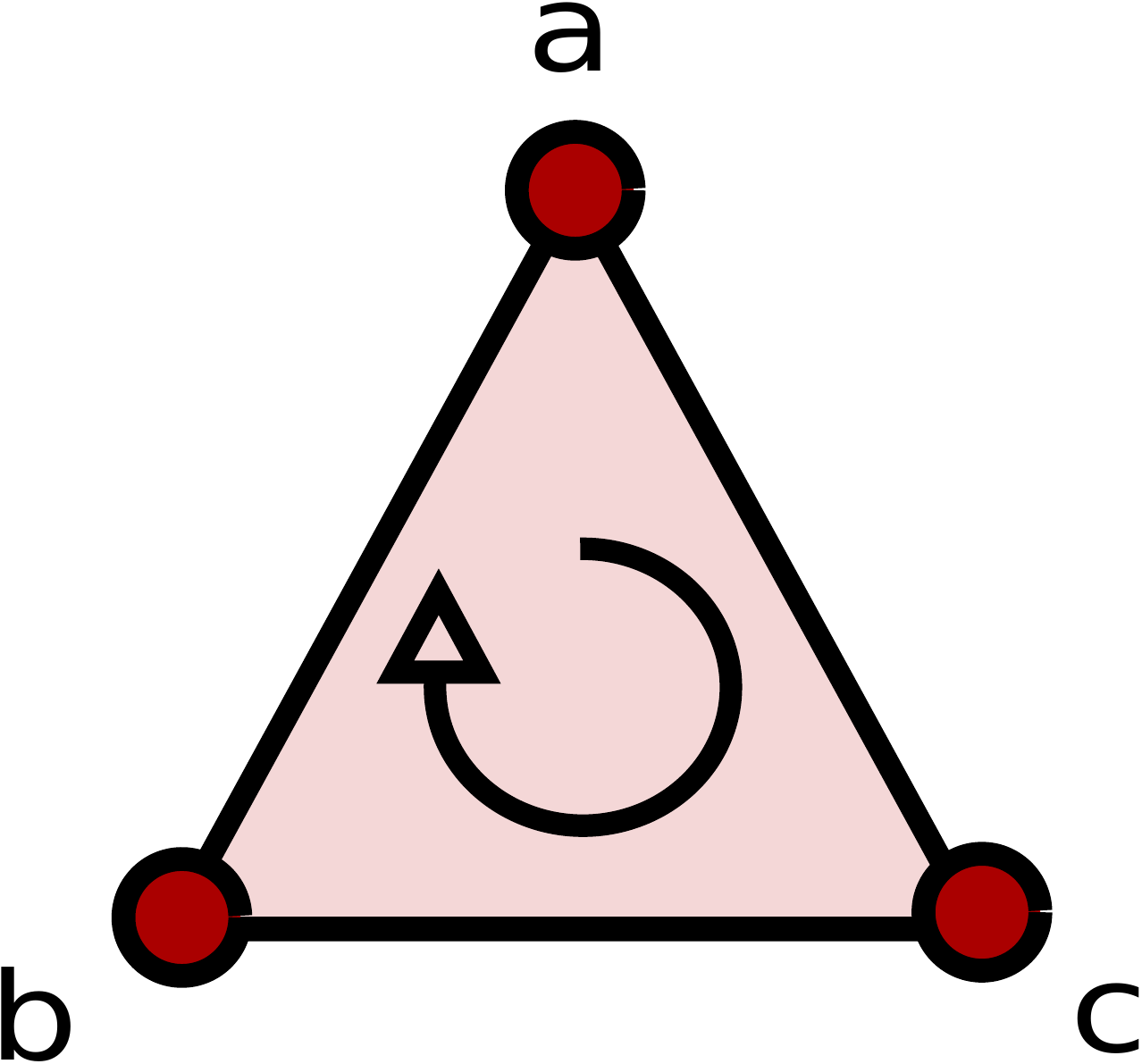}  
  \caption{Simplex represented by ordered set $\{b,a,c\}$}
  \label{fig:sub-second}
\end{subfigure}
\caption{Possible orientations of a $2$-simplex.}
\label{fig:orientation}
\end{figure}

\subsection{Simplicial Homology}

Simplicial homology provides computational techniques to study topological spaces that are represented  as simplicial complexes. We make some basic definitions that are necessary to explain the working principles of these techniques.

\theoremstyle{definition}
\begin{definition}{\textbf{Simplex Orientation}:} The orientation of \hl{a} $k$-simplex is given by the ordering of the vertices in the simplex $[v_1,v_2,...,v_{k+1}]$. Two orderings can define the same orientation if and only if they differ by an even permutation; thus, there are only two allowable orientations of a simplex.
\end{definition}
An example of the possible orientations of a $2$-simplex is shown in Figure \ref{fig:orientation}. We can see that only two orientations are possible for this simplex. Here, each oriented simplex is equal to the negative of the simplex with opposite orientation; mathematically, this is stated as $[a,b,c] = -[b,a,c]$. 

\begin{figure}[!htp]
\begin{subfigure}{.5\textwidth}
  \centering
  \includegraphics[width=.5\linewidth]{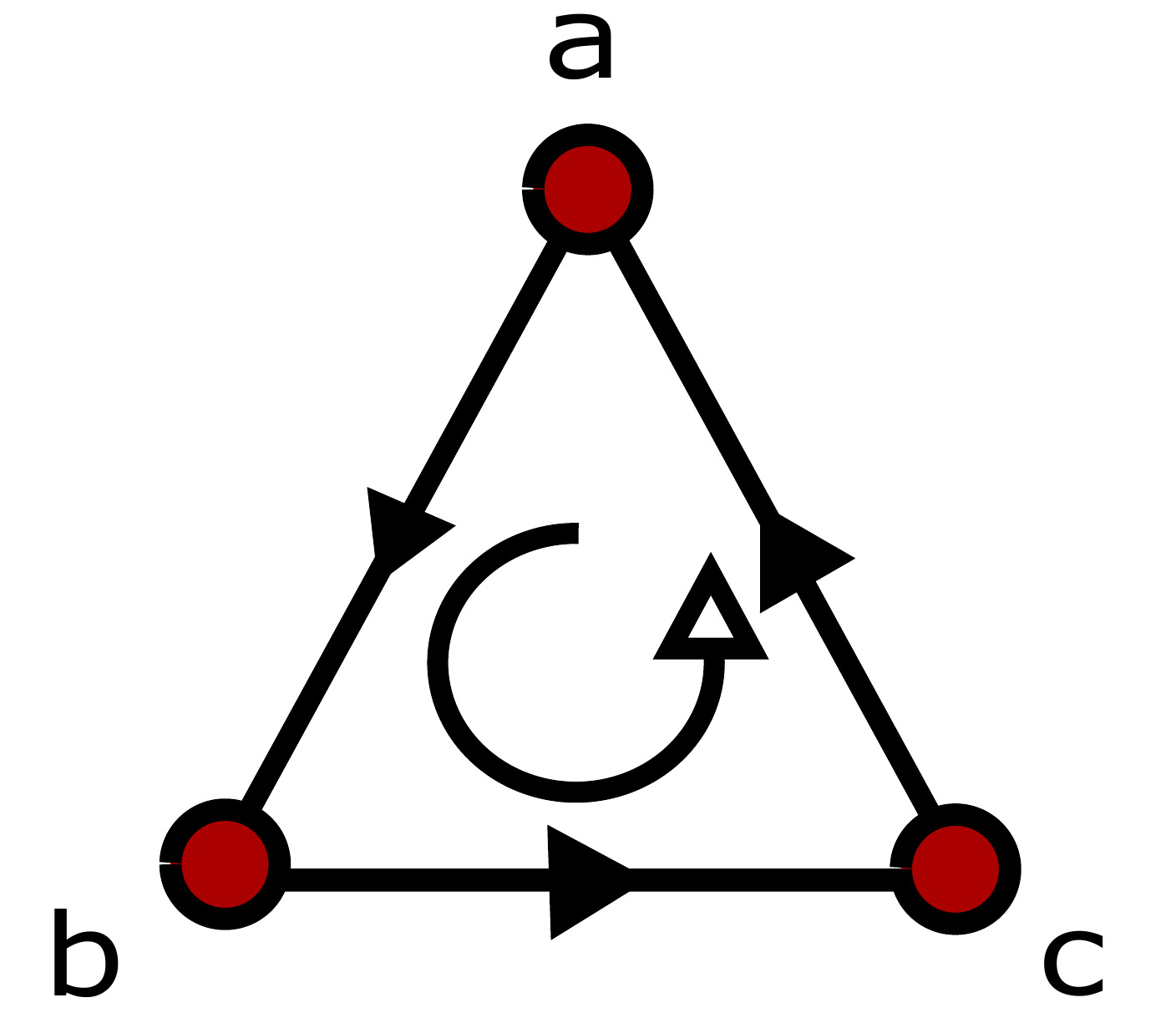}  
  \caption{$1$-Complex ($\mathcal{K}_1$)}
  \label{fig:sub-first}
\end{subfigure}
\begin{subfigure}{.5\textwidth}
  \centering
  \includegraphics[width=.4\linewidth]{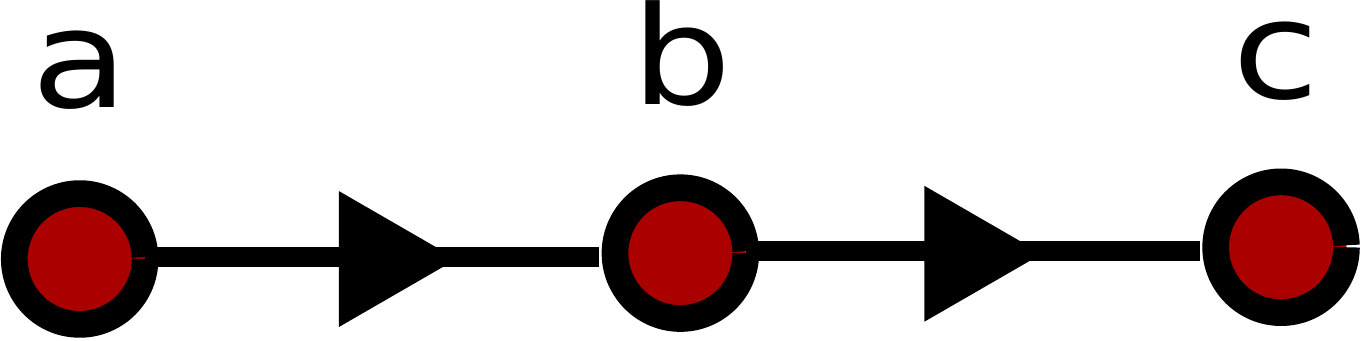}  
  \caption{$1$-Complex ($\mathcal{K}_2$)}
  \label{fig:sub-second}
\end{subfigure}
\caption{Examples of 1-simplicial complexes with the same number of vertices; $\mathcal{K}_1$ contains a hole while $\mathcal{K}_2$ does not.}
\label{fig:loopholeex}
\end{figure}

\begin{figure}[!htp]
    \centering
    \includegraphics[width = .4\textwidth]{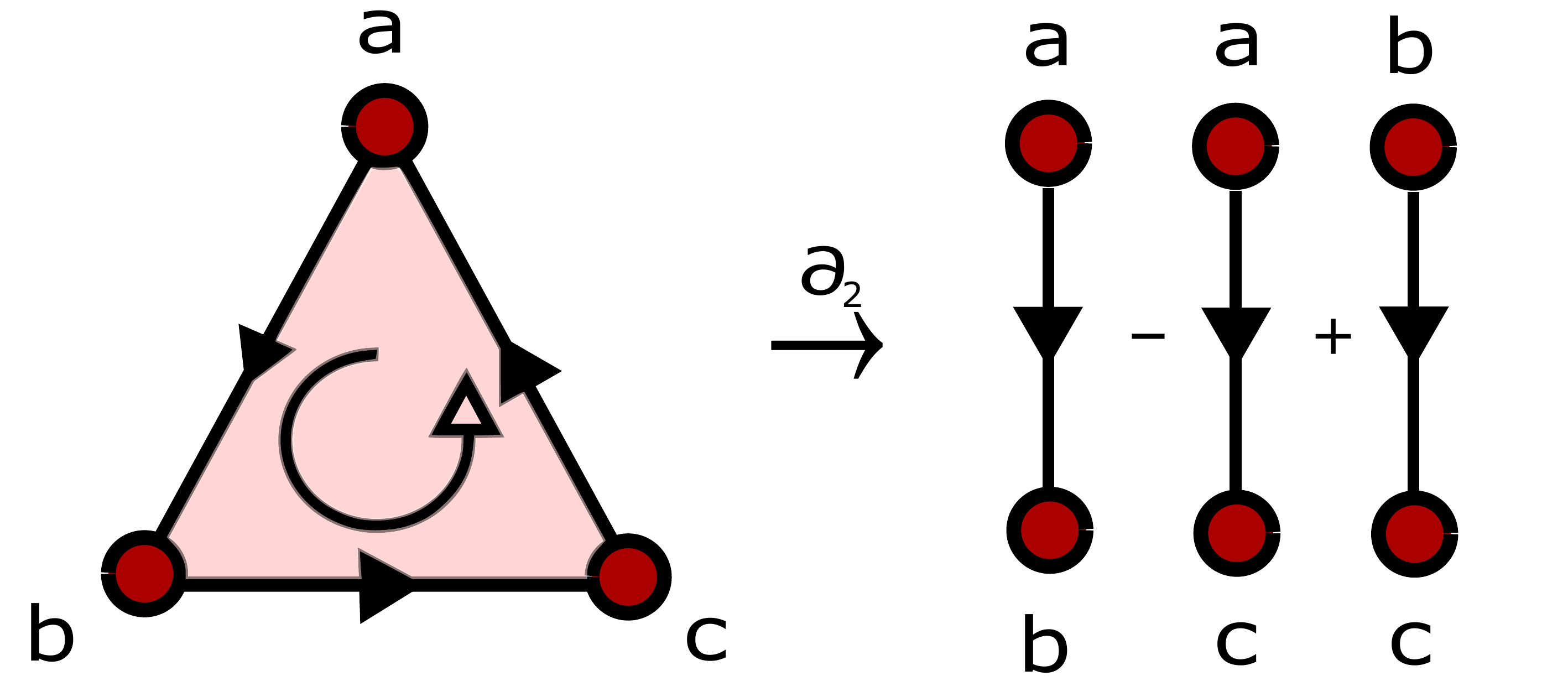}
    \caption{Visualization of the boundary operation ($\partial_2$) applied to a $2$-simplex. The boundary operator maps the $2$-simplex onto its bounding $1$-simplices. The $[a,c]$ simplex is inverted to retain lexicographical ordering.}
    \label{fig:boundary}
\end{figure}

\subsection{Cycles, Holes, and Homology Groups}

In simplicial homology, we want to identify cycles in a given object and we want to know whether a given cycle is bounding a collection of higher dimensional simplices (if it does not, then the cycles is a hole). An example of this concept is presented in Figure \ref{fig:loopholeex}; we see that complex $\mathcal{K}_1$ is a cycle that bounds an empty space which constitutes a hole, while the complex $\mathcal{K}_2$ represents a line with no holes. We now proceed to explore concepts that will allow us to systematically identify the presence of holes and cycles. 

A simplicial $k$-chain on a complex $\mathcal{K}$ is used in the identification of holes in a simplicial complex, and is defined as follows. 
\theoremstyle{definition}
\begin{definition}{\textbf{Simplicial $k$-chains}:} A $k$-chain is a finite weighted sum defined on all $k$-simplicies within a complex $\mathcal{K}$:
\begin{equation}
    \sum_{i=1}^{N}c_i\sigma_i
\end{equation}
where $c_i \in \mathbb{Z}$ and $N$ is the number of $k$-dimensional simplices. The set of $k$-chains on $\mathcal{K}$ is written as $C_k$. Typically, the coefficients are given by $c_i \in \{-1,0,1\}$ and we recall that a value of $-1$ inverts the simplex. 
\end{definition}

The \textit{boundary operator} is a linear operator that maps the $k$-\hl{chains} of a complex to its boundaries. The boundaries are the associated $(k-1)$-\hl{chains} that make up the higher dimensional $k$-\hl{chain}. A visualization of the boundary operator is presented in Figure \ref{fig:boundary}.

\theoremstyle{definition}
\begin{definition}{\textbf{Boundary Operator}:} For a set of $k$-chains (denoted as $C_k$) we define the boundary operator $\partial$ as the mapping:
\begin{equation}
    \partial: C_k \rightarrow C_{k-1}
\end{equation}
The boundary operation on a general simplex $\sigma$ with vertices $[v_0,v_2,...,v_{k}]$ is shown in  \eqref{eq:boundcomp}, where the vertex $\hat{v}_i$ is removed from the set of vertices in the summation. The boundary operation on a $k$-simplex maps the simplex to a summation of its $k-1$ faces.
\begin{equation}
    \partial([v_0,v_1,...,v_k) = \sum_{i=0}^{k+1}(-1)^i[v_0,...,\hat{v}_i,...,v_{k}]
    \label{eq:boundcomp}
\end{equation}
\end{definition}
We use the short-hand notation  $ \partial(C_k) = \partial_k$ to represent a boundary operation. We show an example for a complex $\mathcal{K}$ with $dim(\mathcal{K}) = k$ in \eqref{eq:boundary}; here, \hl{we see that the boundary operator is a mapping from the chains of a higher dimension to chains of lower dimension within the simplicial complex}. We also note that the $k$-chain for dimensions greater than $k$ and less than $0$ are zero and that the boundary maps $\partial_{k+1}$ and $\partial_{0}$ are zero maps. 

\begin{equation}
    0 \xrightarrow{\partial_{k+1}} C_k(\mathcal{K}) \xrightarrow{\partial_{k}} C_{k-1}(\mathcal{K}) \xrightarrow{\partial_{k-1}} ... \xrightarrow{\partial_1} C_0(\mathcal{K})\xrightarrow{\partial_0} 0
    \label{eq:boundary}
\end{equation}

As an example, we apply the boundary operator to the \hl{simplicial chains in the complexes} in Figure \ref{fig:loopholeex}. We note that these  complexes are built from sets of $1$-simplices; thus, when we apply the boundary operator $\partial_1$, we obtain a set of vertices ($0$-simplices). \hl{The result of the operation on complex $\mathcal{K}_1$ is found in \ref{eq:bound1}, where we invert the orientation of the $[a,c]$ boundary in order to retain lexicographical ordering:} 

\begin{equation}
    \partial_1(\mathcal{K}_1) = \partial([a,b]) + \partial([b,c]) - \partial([a,c]) = (a - b) + (b - c) - (a - c) = 0
    \label{eq:bound1}
\end{equation}
and the operation on the 1-chain in the complex $\mathcal{K}_2$ is:
\begin{equation}
    \partial_1(\mathcal{K}_2) = \partial([a,b]) + \partial([b,c])  = (a - b) + (b - c)  = a-c
    \label{eq:bound2}
\end{equation}

This example illustrates how simplicial homology uses basic boundary operations to identify cycles in a complex. We can see that the cycle formed by $\mathcal{K}_1$ is mapped to zero. In simplicial homology, a cycle is defined as a chain that is mapped to the \textit{null space} or \textit{kernel} of the boundary operator (denoted as $ker(\partial)$) and which is zero. With these basic definitions, we can now formally define \textit{cycles} and \textit{boundaries}.

\theoremstyle{definition}
\begin{definition}{\textbf{Cycles}}: \hl{The $k$-dimensional \emph{cycles} are given by}:
\begin{equation}
    Z_k = \text{ker}(\partial_k)
\end{equation}
where $\text{ker}(\partial_k)$ is the kernel of the operator $\partial_k$.
\end{definition}

\begin{definition}{\textbf{Boundaries}: \hl{The \textit{boundaries} are given by:}}
\begin{equation}
    B_k = \text{im}(\partial_{k+1})
\end{equation}
where $\text{im}(\partial_{k+1})$ is the image of the operator $\partial_{k+1}$.
\end{definition}

In summary, the information contained in $Z_k$ gives us the cycles of dimension $k$ within a given complex and the information in $B_{k}$ tells us whether or not a cycle is the boundary of a collection of higher dimensional simplices. It is also important to note that $B_{k}$ is a subgroup of $Z_k$. Also, if a cycle is not a boundary, then it is known as a \textit{hole}. We can summarize this information for any complex by defining the $k$-homology group $H_k$ and the Betti number $\beta_k$. In simple terms, the $H_k$ group contains the unique $k$-holes within a complex and $\beta_k$ counts the number of unique $k$-holes in a complex.

\theoremstyle{definition}
\begin{definition}{\textbf{$k$-Homology Group}: The $k$-homology group $H_k$ is given by the quotient group}:
\begin{equation}
    H_k = Z_k/B_{k}.
\end{equation}

\end{definition}

\begin{figure}[!htp]
\begin{subfigure}{.5\textwidth}
  \centering
  \includegraphics[width=.5\linewidth]{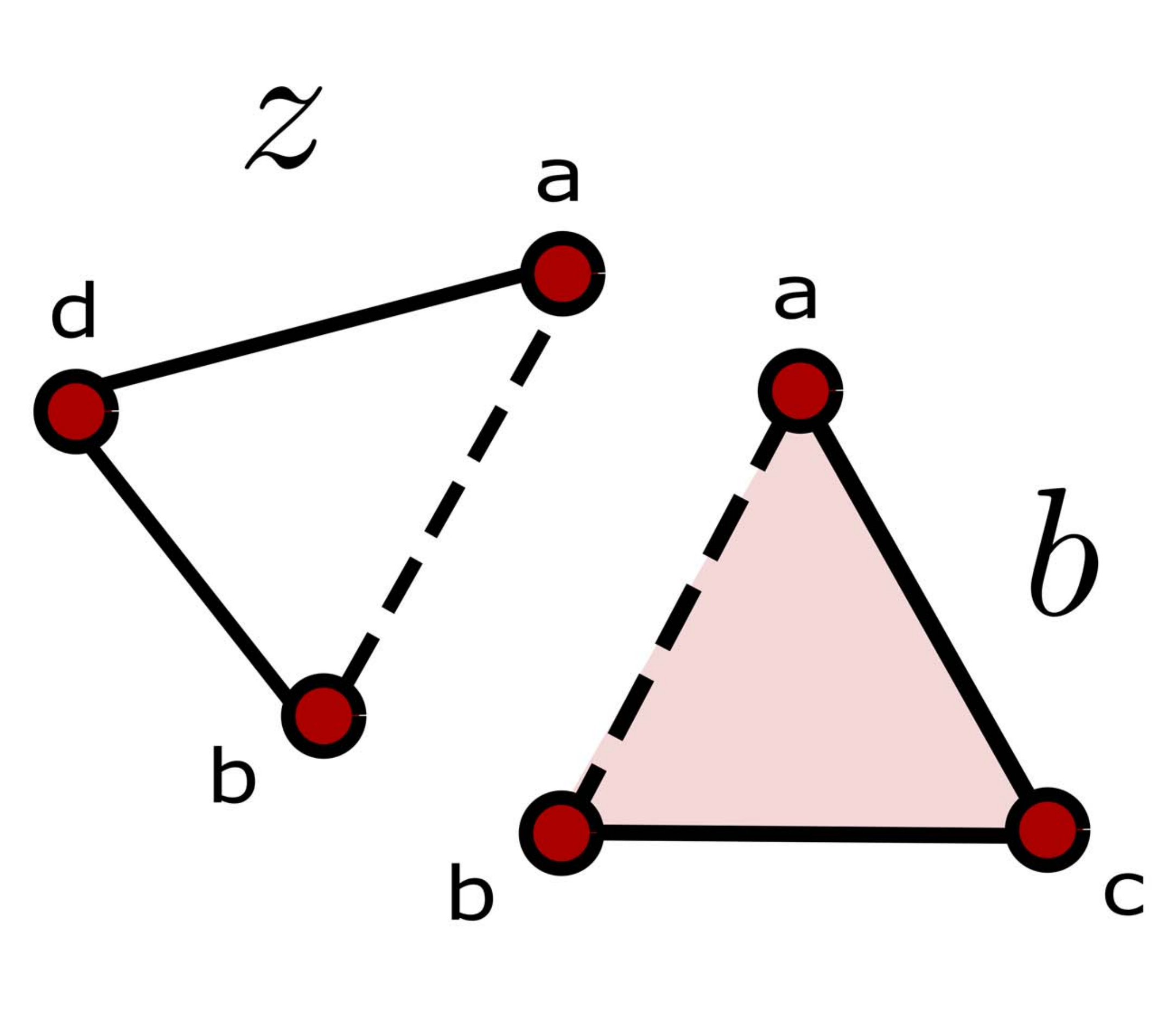}  
  \caption{}
  \label{fig:sub-first}
\end{subfigure}
\begin{subfigure}{.5\textwidth}
  \centering
  \includegraphics[width=.5\linewidth]{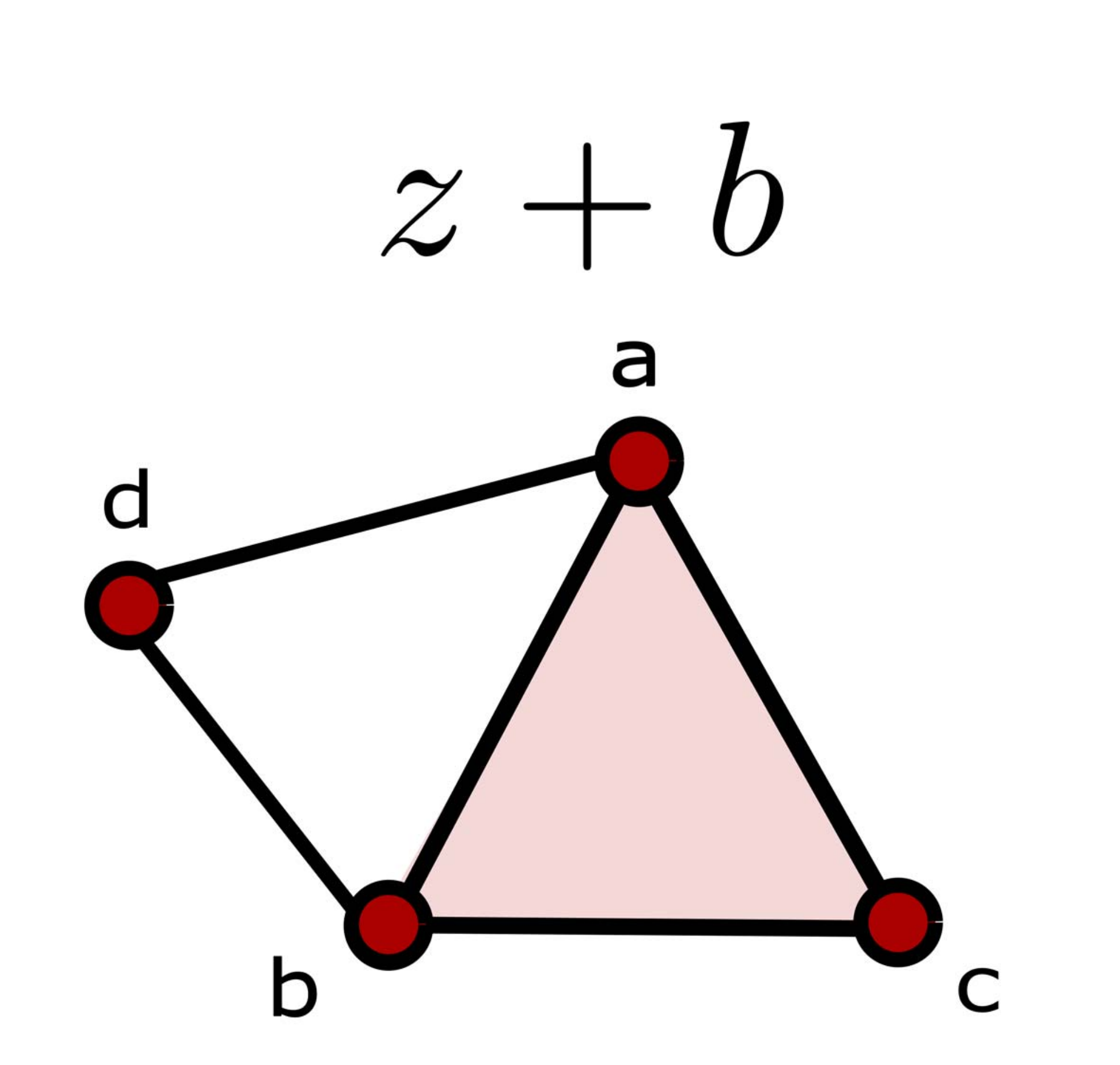}  
  \caption{}
  \label{fig:sub-second}
\end{subfigure}
\caption{An illustration of the homology group $H_1$. (a) $z$ represents a cycle $z\in Z_1$ that is not a boundary ($z \notin B_1$), whereas $b \in Z_1$ is a cycle that bounds a $2$-simplex ($b \in B_1$). (b) The cycle $z$ is homotopy equivalent to $z + b$ ($z \simeq z+b$) and should not be counted as a separate hole.}
\label{fig:Hgroup}
\end{figure}

We illustrate $H_k$ in Figure \ref{fig:Hgroup}; here, we have two simplicial representations $z$ and $b$, where $z$ is a cycle that does not bound a higher dimensional simplex (we have $z \in Z_{1}$ and $z \notin B_{1}$) while $b$ does bound a higher dimensional simplex (we have $b \in B_{1}$ and $b \in Z_1$ as $B_1 \subseteq Z_1$). We can also see that both $z$ and $z+b$ contain the same hole and thus they are homologically equivalent (their difference is a boundary). The homology group $H_k$ formally defines this concept and states that, if a cycle $z_1 = z_2 + B_k$ and $z_1,z_2 \in Z_k$, then the two cycles are equivalent $z_1 \simeq z_2$ and are not independent elements of the group ($H_k$). This logic prevents counting the same topological feature multiple times.

Before we can formally define Betti numbers, we must define the \textit{rank} of a group.

\begin{definition}{\textbf{Group Rank}: We define the $\text{rank}(Z)$ of a group $Z$ as:}
\begin{equation}
  \textrm{rank}(Z) = \min\{|Y|:Y \subseteq Z, \langle Y \rangle = Z\}
\end{equation}
\hl{where $|Y|$ represents the cardinality of set $Y$ and $\langle Y \rangle$ represents the subgroup of $Z$ generated by element of $Y$.}

\end{definition}

We will see that the rank of a group is analogous to the notion of the rank of a matrix (or to the dimension of a vector space). Specifically, it identifies the number of independent basis elements (known as generators) of a group. With this in mind, we define the $k^{th}$-Betti number as follows.

\begin{definition}{\textbf{$k^{th}$-Betti Number}: The $k^{th}$-Betti number $\beta_k$ is the \textrm{rank} of $H_k$ and is given by:}
\begin{equation}
    \beta_k = \text{rank}(H_k) = \text{rank}(Z_k) - \text{rank}(B_{k})
\end{equation}
\end{definition}
We illustrate these concepts using the complex $\mathcal{K}_1$ presented in Figure \ref{fig:loopholeex}. This complex presents a hole; performing the necessary computations to obtain $Z_1$ and $B_1$ for $\mathcal{K}_1$ we find the that:
\begin{subequations}
\begin{align}
    Z_1(\mathcal{K}_1) &= \{(a,b) + (b,c) - (a,c)\} \\
    B_1(\mathcal{K}_1) &= \{0\}
\end{align}
\end{subequations}
We now perform the quotient operation to obtain $H_1$ and $\beta_1$ for $\mathcal{K}_1$:
\begin{subequations}
\begin{align}
    H_1(\mathcal{K}_1) &= \{(a,b) + (b,c) - (a,c)\}/\{0\} = \{(a,b) + (b,c) - (a,c)\}\\
    \beta_1(\mathcal{K}_1)  &= \text{rank}(\{(a,b) + (b,c) - (a,c)\}) = 1
\end{align}
\end{subequations}
We see that the $H_1$ group identifies the $1$-dimensional holes within the complex, and $\beta_1$ counts the number of $1$-dimensional holes in the complex. The same is true for all other dimensions $k \geq 0$ as long as $\dim(\mathcal{K}) \geq k$ because $H_k = 0$ for all $k > \dim(\mathcal{K})$.  \hl{This concept becomes familiar when viewed from the perspective of linear algebra and  matrices. One can consider the $H_k$ group similar to the basis vectors for a given matrix, where in this case $H_k$ defines the \textit{topological} basis for a given shape and the Betti numbers ($\beta_k$) correspond to the rank of this bases and can be seen as the total number of unique topological features. The main difference is that for a given shape there can be multiple sets of topological bases for each dimension of the shape. The goal of TDA is to identify and compare the topological bases for each shape, similar to how one might compare the structure of matrices based upon their basis vectors and corresponding rank.}

\subsubsection{The $0^{th}$ Homology Group ($H_0$)}

The $0^{th}$ homology group $H_0$ plays an important role in topological analysis. $H_0$ is the measure of the number of \textit{connected components} in a complex. A \textit{component} (or \textit{subcomplex}) is defined as follows. 

\theoremstyle{definition}
\begin{definition}{\textbf{Subcomplex:} A subcomplex is a subset $\mathcal{S}$ of a complex $\mathcal{K}$ such that $\mathcal{S}$ is also a complex.}  

\end{definition}

\theoremstyle{definition}
\begin{definition}{\textbf{Connected complex:} A complex is connected if there exists a path made of $1$-simplices from any vertex of the complex to any other vertex.} 

\end{definition}

\begin{figure}[!htp]
\begin{subfigure}{.5\textwidth}
  \centering
  \includegraphics[width=.5\linewidth]{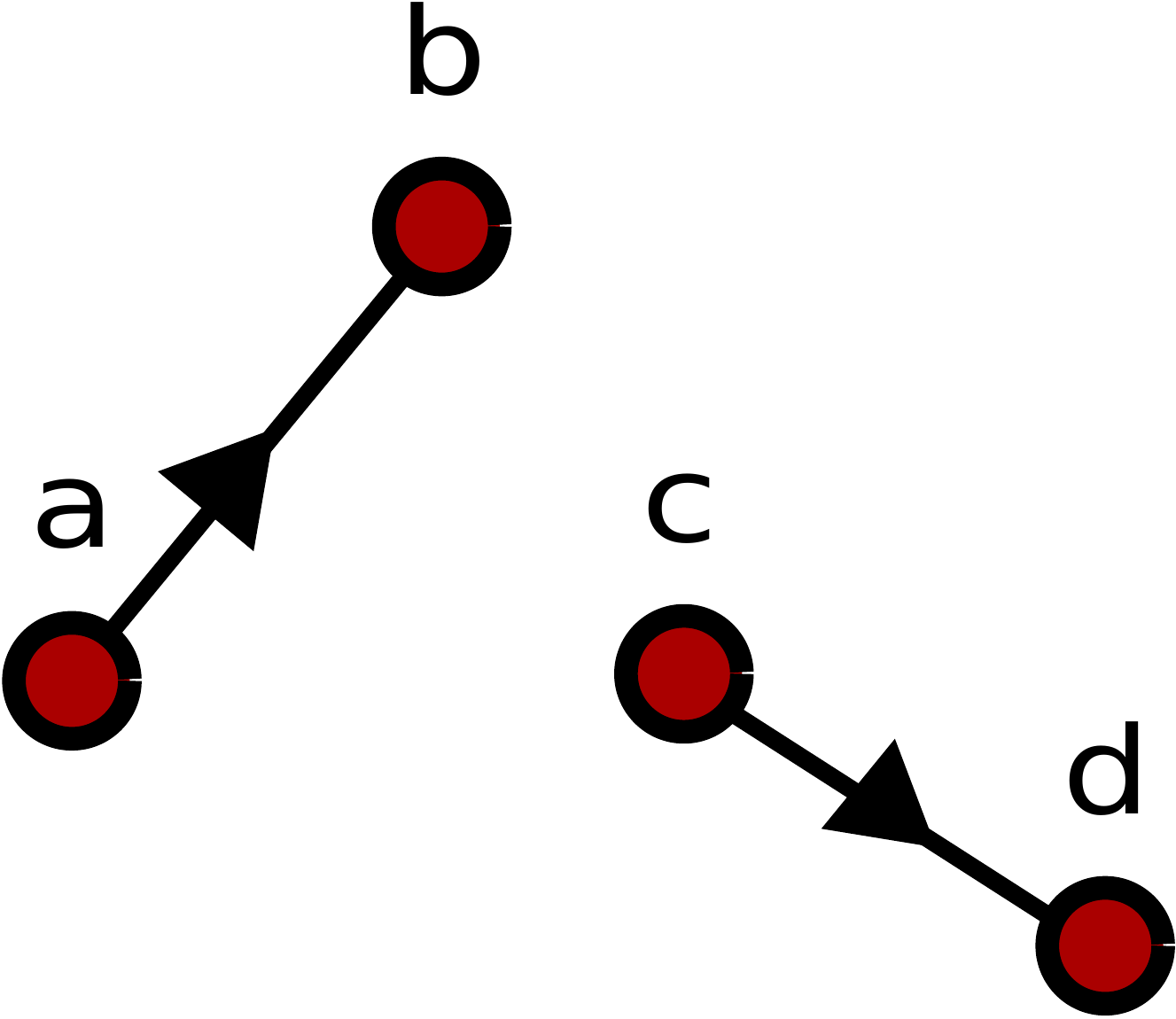}  
  \caption{$1$-complex $\mathcal{K}_3$}
  \label{fig:sub-first}
\end{subfigure}
\begin{subfigure}{.5\textwidth}
  \centering
  \includegraphics[width=.5\linewidth]{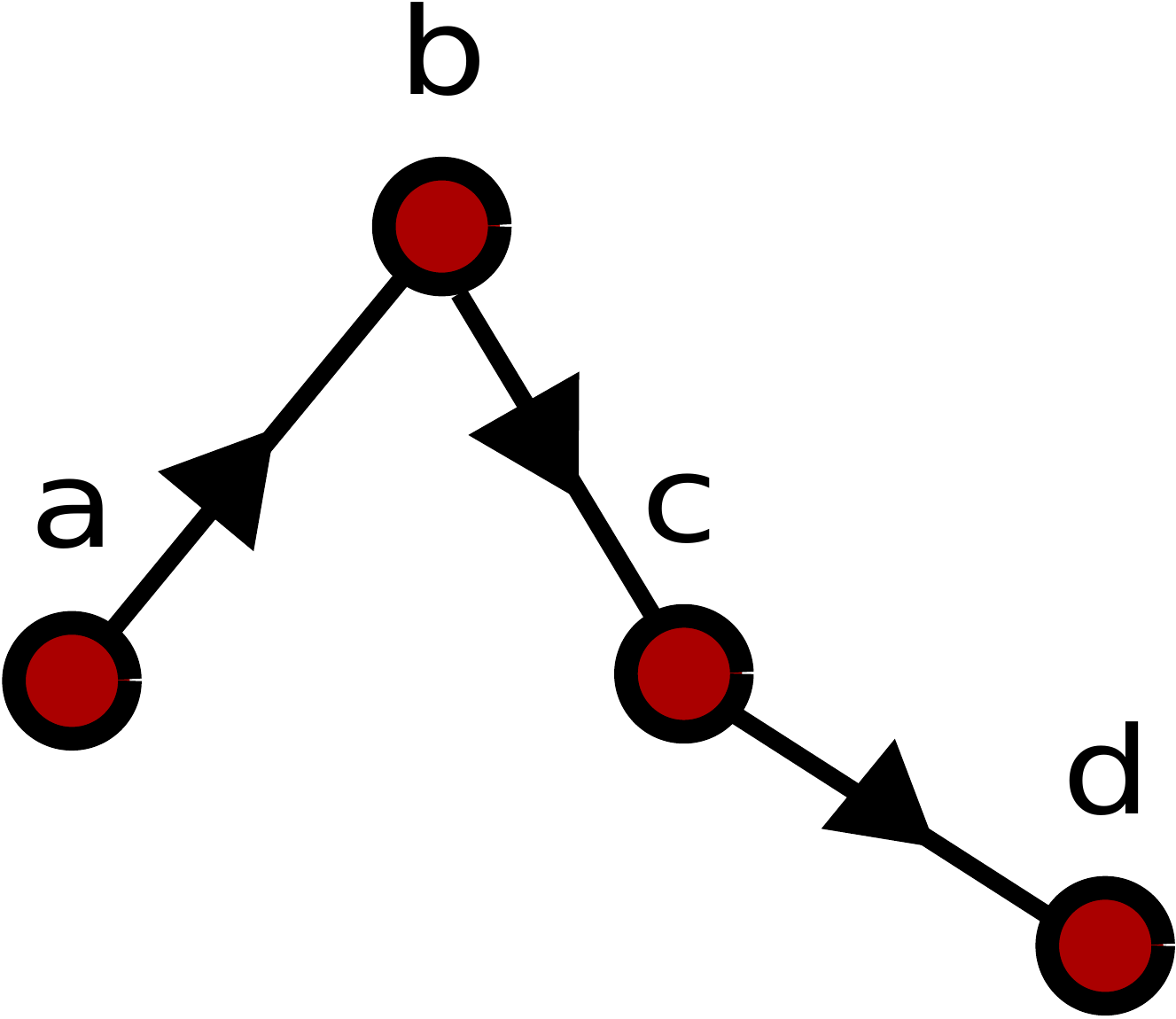}  
  \caption{$1$-complex $\mathcal{K}_4$}
  \label{fig:sub-second}
\end{subfigure}
\caption{Example complexes with different number of connected components. $\mathcal{K}_3$ has two connected components and $\mathcal{K}_4$ has a single connected component.}\label{fig:concompex}
\end{figure}

The group $H_0$ describes how many disconnected subcomplexes $\mathcal{S}$ there are within a given complex $\mathcal{K}$. A simple example is shown in Figure \ref{fig:concompex}; here, we note that the first complex $\mathcal{K}_3$ has two disjoint subcomplexes, and the second complex $\mathcal{K}_4$ has a single connected component.  The calculations for $\mathcal{K}_3$  are given by:
\begin{subequations}
\begin{align}
    H_0(\mathcal{K}_3) &= \{a,b,c,d\}/\{(a-b),(c-d)\}\\
    \beta_0(\mathcal{K}_3) &= 2
\end{align}
\end{subequations}
and for $\mathcal{K}_3$ are:
\begin{subequations}
\begin{align}
    H_0(\mathcal{K}_4) &= \{a,b,c,d\}/\{(a-b),(b-c),(c-d)\}\\
    \beta_0(\mathcal{K}_4) &= 1.
\end{align}
\end{subequations}

\section{Computational Methods for TDA}
\label{section:matrix}

Now that we have developed a basic understanding of simplicial homology, we can begin to streamline computations through the tools of numerical linear algebra. \hl{In order to simplify our discussion, we will define a new $k$-chain where the coefficients $c_i \in \mathbb{Z}_2$ where $\mathbb{Z}_2$ is the set of binaries $\{0,1\}$ (rather than $c_i \in \mathbb{Z}$), for which $1+1 = 0$}. With this new definition, we can remove the need for defining an orientation on a simplicial complex. In \eqref{eq:step1}-\eqref{eq:step3}, we see that computational results are the same as those for the example provided in \eqref{eq:bound1}.
\begin{subequations}
\begin{align}
    C_1(\mathcal{K}_1) &= [a,b] + [b,c] + [a,c]
    \label{eq:step1}\\
    \partial_1(\mathcal{K}_1) &= \partial([a,b]) + \partial([b,c]) + \partial([a,c]) = (a + b) + (b + c) + (a + c)
    \label{eq:step2}\\
    \partial_1(\mathcal{K}_1) &= (1+1)a + (1+1)b + (1+1)c = 0 + 0 + 0 = 0 
    \label{eq:step3}
    \end{align}
\end{subequations}

We can now couple the newly defined $k$-chain with the boundary operator to create what is known as a \textit{boundary matrix} $\mathbf{B} \in \mathbb{Z}_2^{n\times m}$ where $n$ represents the \hl{simplices} of dimension $(k-1)$ in $\mathcal{K}$ and $m$ represents the simplices of dimension $k$ in $\mathcal{K}$. The boundary matrix is built based on the following rules:

\begin{itemize}
    \item The face of a simplex precedes the simplex in column (row) index.
    \item An entry of one is placed in position ($j,i$) if $\sigma_i$ is a face of $\sigma_j$, otherwise it is zero. 
    \item If two simplices are of the same dimension, lexicographic ordering is used (i.e., $[a,b] < [a,c] < [b,c])$.
\end{itemize}

For example, we take complex in Figure \ref{fig:simplicial_complex} and compute the $H_1$ group. We show matrix $\mathbf{B}$ for $\partial_0$ in Table \ref{tab:boundmat0}, for $\partial_1$ in Table \ref{tab:boundmat1}, and for $\partial_2$ in Table \ref{tab:boundmat2}. Again, the boundary matrix associated with $\partial_1$ is mapping the $1$-simplices of the complex to its boundaries and similarly for $\partial_2$.

\begin{table}[!htp]
    \centering
    \begin{tabular}{c|cccccc}
        $\partial_{0}$ & [a] & [b] & [c] & [d] & [e]\\
        \hline
        [0] & 0 & 0 & 0 & 0 & 0 \\
        
    \end{tabular}
    \caption{Boundary matrix $\mathbf{B_{0}}$ for $\partial_0$ of complex in Figure \ref{fig:simplicial_complex}}
    \label{tab:boundmat0}
\end{table}

\begin{table}[!htp]
    \centering
    \begin{tabular}{c|cccccc}
        $\partial_{1}$ & [a,b] & [a,c] & [a,d] & [b,c] & [b,d] & [d,e]\\
        \hline
        [a] & 1 & 1 & 1 & 0 & 0 & 0 \\
        
        [b] & 1 & 0 & 0 & 1 & 1 & 0\\
        
        [c] & 0 & 1 & 0 & 1 & 0 & 0  \\
        
        [d] & 0 & 0 & 1 & 0 & 1 & 1 \\
        
        [e] & 0 & 0 & 0 & 0 & 0 & 1 \\
        
    \end{tabular}
    \caption{Boundary matrix $\mathbf{B_1}$ for $\partial_1$ of complex in Figure \ref{fig:simplicial_complex}}
    \label{tab:boundmat1}
\end{table}

\begin{table}[!htp]
    \centering
    \begin{tabular}{c|cccccc}
        $\partial_{2}$ & [a,b,c]\\
        \hline
        [a,b] & 1 \\
        
        [a,c] & 1 \\
        
        [b,c] & 1 \\
        
        [b,d] & 0 \\
        
        [d,e] & 0 \\
        
    \end{tabular}
    \caption{Boundary matrix $\mathbf{B_2}$ for $\partial_2$ of complex in Figure \ref{fig:simplicial_complex}}
    \label{tab:boundmat2}
\end{table}

In order to compute $\text{rank}(Z_j)$ and $\text{rank}(B_j)$ and the number of $k$-holes in the complex (the Betti number $\beta_k$), we must first reduce the matrices to a canonical form known as the \textit{Smith normal form} (SNF).

\begin{definition}{\textbf{Smith Normal Form}}: A matrix $\mathbf{M} \in \mathbb{Z}_2^{n \times m}$ is in \textit{Smith normal form} if it is diagonal and if it can be obtained by multiplying $\mathbf{M}$ by invertible matrices $\mathbf{S} \in \mathbb{Z}_2^{n \times n}$ and $\mathbf{T} \in \mathbb{Z}_2^{m \times m}$ as  $ \mathbf{M_{SNF}} = \mathbf{S} \mathbf{M} \mathbf{T}$. 
\end{definition}

The reduced matrix $\mathbf{B_{1_{SNF}}}$ is shown in Table \ref{tab:redboundmat1} and $\mathbf{B_{2_{SNF}}}$ is shown in Table \ref{tab:redboundmat2}. The $\partial_0$ matrix cannot be further reduced and thus it is not shown. The SNF matrices contain all the information required to find $\textrm{rank}(Z_k)$ and $\textrm{rank}(B_k)$ for $k \in \{1,2\}$:

\begin{subequations} 
\begin{align}
    \textrm{rank}(Z_k) &= m - \textrm{rank}(\mathbf{B_{(k)_{SNF}}})
    \label{eq:rankz}\\
    \textrm{rank}(B_k) &= \textrm{rank}(\mathbf{B_{(k+1)_{SNF}}})
    \label{eq:rankb}
    \end{align}
\end{subequations}

These simple calculations show us that $\textrm{rank}(Z_1) = 2$ and $\textrm{rank}(B_1) = 1$. From this we can see that the number of $1$-holes in our complex is $\beta_1 = \textrm{rank}(Z_1) - \textrm{rank}(B_1) = 1$, as expected.
\begin{table}[!htp]
    \centering
    \begin{tabular}{c|cccccc}
        $\partial_{1}$ & [a,b] & [a,c] & [a,d] & [b,c] & [b,d] & [d,e]\\
        \hline
        [a] & 1 & 0 & 0 & 0 & 0 & 0 \\
        
        [b] & 0 & 1 & 0 & 0 & 0 & 0\\
        
        [c] & 0 & 0 & 1 & 0 & 0 & 0  \\
        
        [d] & 0 & 0 & 0 & 0 & 0 & 1 \\
        
        [e] & 0 & 0 & 0 & 0 & 0 & 0 \\
        
    \end{tabular}
    \caption{Reduced boundary matrix $\mathbf{B_{1_{SNF}}}$ for $\partial_1$ of the simplicial complex in Figure \ref{fig:simplicial_complex}. \hl{We note that the matrix is not diagonal in order to retain lexicoraphical ordering, but can be easily made diagonal.}}
    \label{tab:redboundmat1}
\end{table}

\begin{table}[!htp]
    \centering
    \begin{tabular}{c|cccccc}
        $\partial_{2}$ & [a,b,c]\\
        \hline
        [a,b] & 1 \\
        
        [a,c] & 0 \\
        
        [b,c] & 0 \\
        
        [b,d] & 0 \\
        
        [d,e] & 0 \\
        
    \end{tabular}
    \caption{Reduced boundary matrix $\mathbf{B_{2_{SNF}}}$ for $\partial_2$ of the simplicial complex in Figure \ref{fig:simplicial_complex}.}
    \label{tab:redboundmat2}
\end{table}

\section{Persistent Homology}
\label{seq:mathend}
Persistent homology is a methodology originally proposed by Edelsbrunner, Letscher, and Zomorodian and further developed by many others for extracting and quantifying topological information from data \cite{edelsbrunner2000topological}. This methodology is discussed in detail in   \cite{zomorodian2005computing,bauer2014induced,chazal2012structure,carlsson2009topology,cohen2007stability,ghrist2008barcodes}. The homology of data provides deep insight into the structure of the data and quantification capabilities of their geometric features \cite{ghrist2008barcodes,cohen2007stability}.

\begin{figure}[!htp]
    \centering
    \includegraphics[width = .7\textwidth]{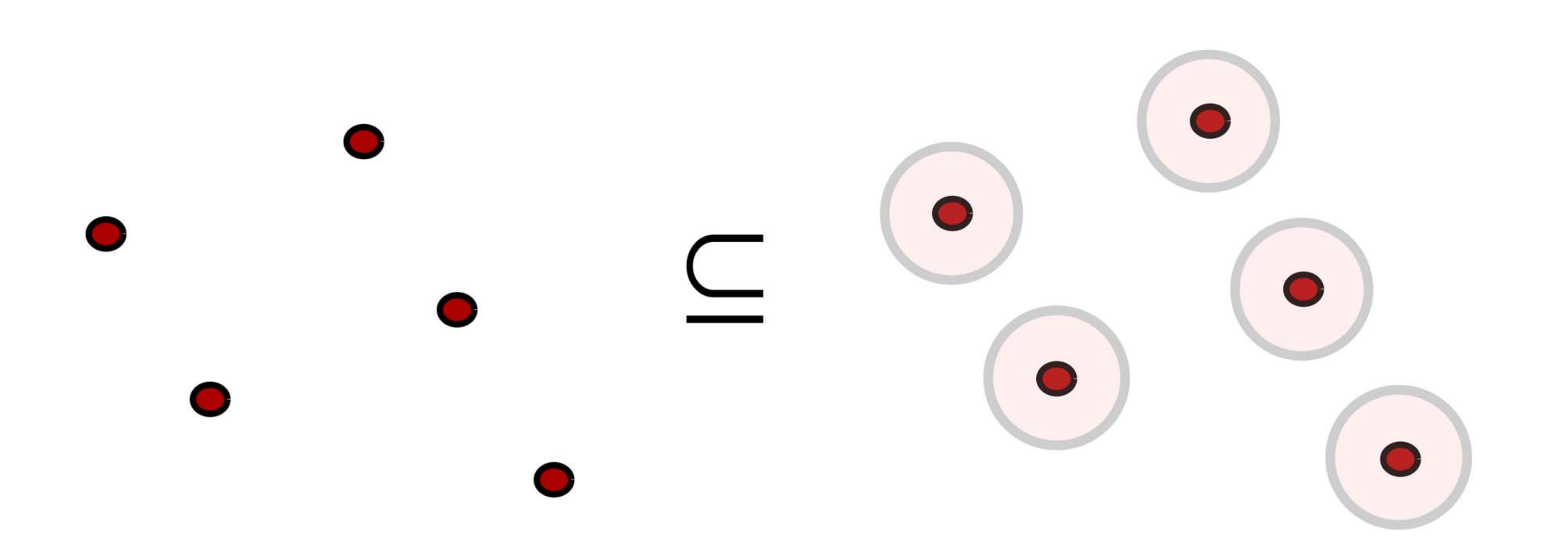}
    \caption{A cover of points $x_i \in \mathcal{X}$ is defined by a set of balls $B(x_i,\epsilon)$ expanded around each point.}
    \label{fig:coverex}
\end{figure}

\subsection{Building Simplicial Complexes from Data}

Direct computation of the homology of an arbitrary set defined by the space $\mathcal{X} \subset \mathbb{R}^n$ is a complex task. To simplify this task, we leverage simplicial homology; here, we identify a simplicial complex $\mathcal{K}$ such that its homology is the same or similar to that of $\mathcal{X}$. We can define such a complex $\mathcal{K}$ by creating a geometric object known as the \textit{cover} ($\mathcal{U}$) of $\mathcal{X}$. 

\begin{definition}{\textbf{Cover}:} $\mathcal{U} = \{U_i\}_{i\in I}$ is a cover of a metric space $\mathcal{X}$ if $\mathcal{X} \subseteq \bigcup\limits_{i\in I} U_i.$
\end{definition}

If we define our metric space ($\mathcal{X}$) to be a finite set of points in $\mathbb{R}^n$, then we can imagine each set $U_i$ to be a ball $\{B(x_i,\epsilon):x_i \in \mathcal{X}, \epsilon \in \mathbb{R}^+ \}$ centered around each point $x_i \in \mathcal{X}$,  where $\epsilon$ is the ball's radius. An example of such a cover is presented in Figure \ref{fig:coverex}. We now utilize this cover to develop a simplicial complex, known as a \textit{\v{C}ech complex}.

\begin{definition}{\textbf{\v{C}ech Complex}: The \v{C}ech Complex ($\check{C}$) is a simplicial complex built from $k$-simplicies that are the non-empty intersection of $k$+1 sets of a cover $\mathcal{U}$.}
\end{definition}

A \v{C}ech complex is also known as the ``nerve" of the cover $\mathcal{U}$.

\begin{definition}{\textbf{Nerve}: The nerve of collection $\mathcal{U} = \cup_{i \in I}\{U_i\}$ is the simplicial complex with vertices $I$ and $k$-simplices built from $\{i_0,i_1,...,i_k\}$ if and only if $U_{i_0} \cap U_{i_1} \cap ... \cap U_{i_k} \neq \varnothing$.}
\end{definition}

\begin{figure}[!htp]
    \centering
    \includegraphics[width = 1\textwidth]{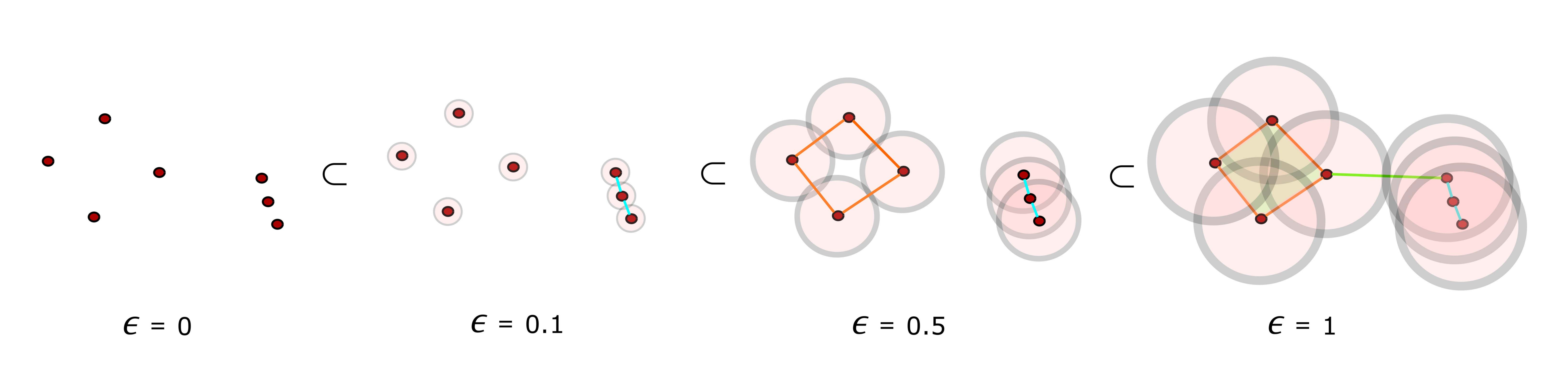}
    \caption{Filtration of points $x_i \in \mathcal{X}$ by a set of balls $B(x_i,\epsilon)$ with expanding $\epsilon$. As $\epsilon$ is increased the topology of the \v{C}ech complex evolves and this introduces holes and higher dimensional simplices. This filtration builds a filtered simplicial complex ($\mathcal{K}_{\epsilon = 0} \subset \mathcal{K}_{\epsilon = 0.1} \subset \mathcal{K}_{\epsilon = 0.5} \subset \mathcal{K}_{\epsilon = 1}$).}
    \label{fig:filterex}
\end{figure}

\hl{By using the so-called ``Nerve Theorem" construction, we can build a simplicial complex $\mathcal{K}$ for space $\mathcal{X} = \cup_{i \in I}\{U_i\}$. Under certain assumptions, the simplicial complex $\mathcal{K}$ and the space $\mathcal{X}$ are homotopy equivalent \cite{alexandroff1928allgemeinen,ghrist2014elementary}.} With this, we can apply the calculations and analysis of simplicial homology directly to our data. However, there is one caveat that is important to note, which is the selection of the distance $\epsilon$. An example of the the nerve of a dataset with varying levels of $\epsilon$ is shown in Figure \ref{fig:filterex}. We can see that, as we adjust $\epsilon$ of the cover $\mathcal{U}$, we obtain different \v{C}ech complexes, each with a different homology. We want to ensure that the homology captures the most interesting features of the data. In the complex shown in  Figure \ref{fig:filterex}, these features are the two clusters of points; here, one cluster forms a loop and the other does not. It is easy to see in this example what range of $\epsilon$ values would be most effective at capturing this information. However, if our dataset is of much higher dimension, finding the correct $\epsilon$ is much more difficult. Consequently, we characterize the dataset for multiple values of $\epsilon$. This information is captured by a \textit{filtered simplicial complex} \cite{carlsson2009topology}.

\begin{definition}{\textbf{Filtered Simplicial Complex}: A filtered simplicial complex $\mathcal{K} \in \mathbb{R}^m$ is a simplicial complex for which there is a series of nested simplicial subcomplexes $\mathcal{K}_\epsilon \in \mathbb{R}^m$ such that:}
\begin{equation}
    \mathcal{K}_{\epsilon_0} \subset \mathcal{K}_{\epsilon_1} \subset ... \subset \mathcal{K}_{\epsilon_n} 
\end{equation}
where $\epsilon_0<\epsilon_1<\cdots <\epsilon_n$. 

\end{definition}

Referring back to the example in which we define the cover $\mathcal{U}$ as a set of balls $B(x_i,\epsilon)$ of radius $\epsilon$, we can view the filtered complex as the set of \hl{nerve complexes} that are formed as we expand $\epsilon$. Figure \ref{fig:filterex} demonstrates that the filtration of nerve complexes $\mathcal{K}_{\epsilon = 0.1} \subset \mathcal{K}_{\epsilon = 0.3} \subset \mathcal{K}_{\epsilon = 0.5} \subset \mathcal{K}_{\epsilon = 1}$. Note that the interesting features of the data (the two clusters and one loop) are present in the homology of the subcomplexes and persist during a large portion of the filtration. The main goal of this analysis is to identify \textit{persistence intervals} in the complex filtered by $\epsilon$. Given a topological feature present in a filtration, we identify the value of $\epsilon$ where the feature is \textit{born} (appears) and \textit{dies} (disappears).

\begin{definition}{\textbf{Birth}: For a filtered complex $\mathcal{K}$ and subcomplexes $\mathcal{K}_i,\mathcal{K}_j$ where $i < j$. A topological feature $x \in H_p(\mathcal{K}_j)$ is \textit{born} at $j$ if $x \notin H_p(\mathcal{K}_i).$} 
\end{definition}

\begin{definition}{\textbf{Death}: For a filtered complex $\mathcal{K}$ and subcomplexes $\mathcal{K}_i,\mathcal{K}_j$ where $i < j$. A topological feature $x \in H_p(\mathcal{K}_i)$ \textit{dies} at $j$ if $x \notin H_p(\mathcal{K}_j)$. A feature will also die if the feature merges with a feature born earlier in the filtration, this is known as the \textit{elder rule}.}
\end{definition}

\begin{definition}{\textbf{Persistence Interval}: For a given topological feature $x$ with birth point $i$ and death point $j$, the \textit{persistence interval} ($\textrm{Int}$) for the feature is given by:}
\begin{equation}
    \textrm{Int} = [i,j) : i,j \in \ \bar{\mathbb{R}}
\end{equation}
If $j = \infty$ then the component does not die during the filtration (persists forever).

\end{definition}

With a filtration we are identifying the appearance and disappearance of topologically interesting features in our dataset. The filtered complex in Figure \ref{fig:filterex} demonstrates this concept. We can see that the hole in the dataset is born at $\epsilon=0.5$ and is completely filled in at $\epsilon=1$, persisting for a majority of the filtration. We can also see that the individual points (at $\epsilon = 0$) become two connected components at $\epsilon=0.5$ and then become a single component a $\epsilon=1$. Thus, the longest persistence intervals in both $H_1$ and $H_0$ capture the defining topological characteristics of the filtered complex.

\subsection{Persistence Diagrams}
\label{sec:pds}
\begin{figure}[!htp]
    \centering
    \includegraphics[width = 0.8\textwidth]{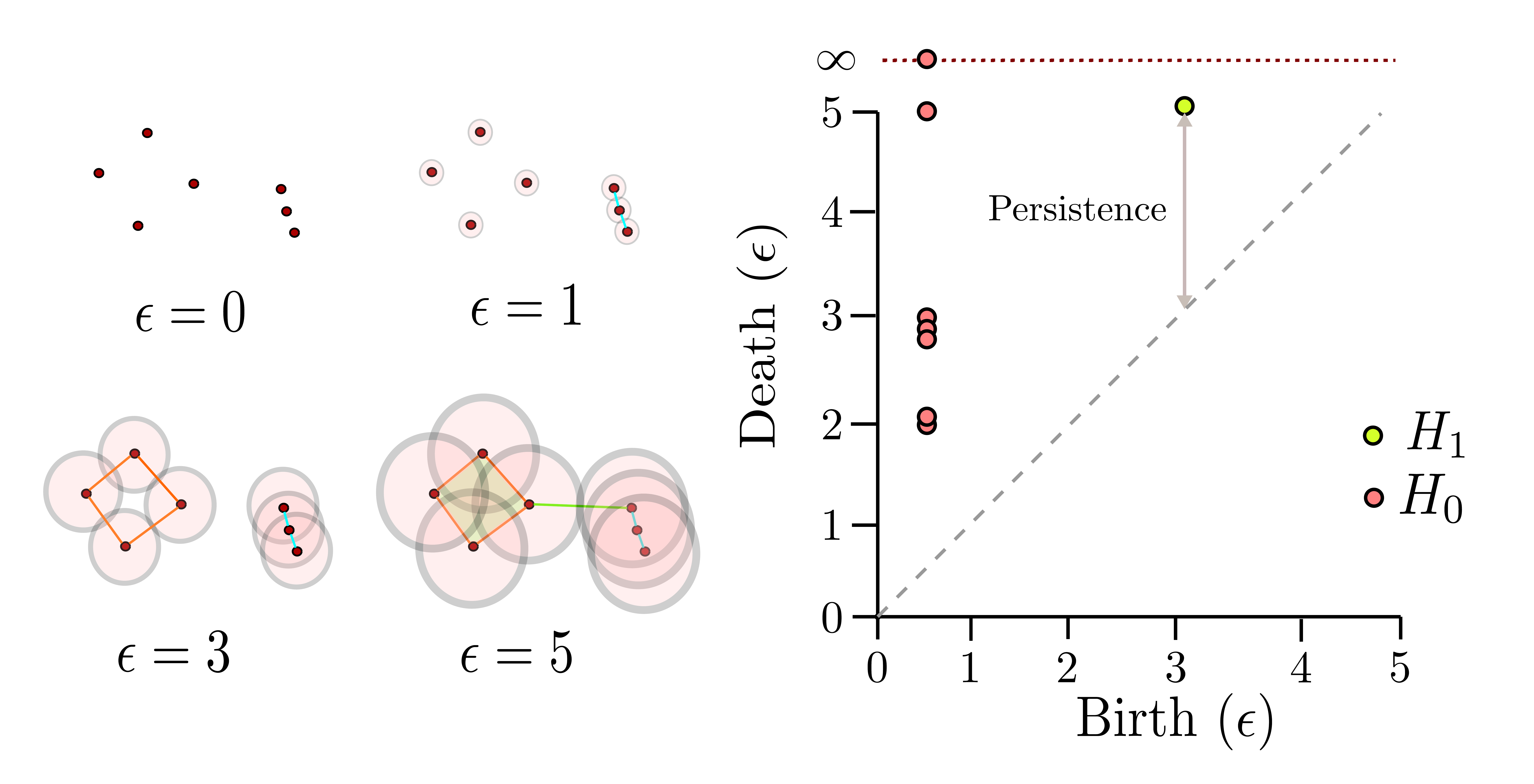}
    \caption{Filtration of points $x_i \in \mathcal{X}$ by a set of balls $B(x_i,\epsilon)$ with expanding $\epsilon$ and its corresponding persistence diagram. The PD records the $\epsilon$ value at which topological features are born and the $\epsilon$ value of their death during the filtration. For example, the cycle born at $\epsilon = 3$ ($x = 3$) dies at $\epsilon = 5$ ($y = 5$) when it is filled in, with a total persistence of $5-3 = 2$, which is seen in the PD.}
    \label{fig:points_pers}
\end{figure}

The information about topological features contained within a filtered complex is summarized into what is known as a \textit{persistence diagram}  (PD)  \cite{edelsbrunner2000topological}. This can be computed via an extension of the matrix methods presented in Section \ref{section:matrix} \cite{edelsbrunner2010computational}. The PD is a visual method that represents the \textit{birth} ($x$) and \textit{death} ($y$) of topological features as a set of points  in ${\mathbb{R}}^2$. The persistence diagram associated with the filtration in Figure \ref{fig:filterex} is shown in Figure \ref{fig:points_pers}. \hl{This diagram represents the \textit{birth} and \textit{death} of the features of the $H_0$ and $H_1$ homology groups for the filtered complex. The \textit{persistence interval} associated with each feature is the vertical line segment between the persitence point and the diagonal. This information allows for a direct visual understanding of the topology of the dataset.}

An important property of PDs is that they can be \textit{vectorized} to enable quantification and these vectors can be used to perform tasks such as regression, classification, or clustering. For instance, as we will see in Section (\ref{section:examples}), one can apply PCA to the vectorized PDs to identify clusters defined by different topological features. 

\begin{definition}{\textbf{Vectorization}}: A persistence diagram $PD_{X_i}$  of a dataset $X_i$ is vectorized by mapping the $PD_{X_i}$ to a vector $\overrightarrow{PD}_{X_i}  \in \mathbb{R}^q$.

\end{definition}

There are multiple ways to vectorize a PD, the most widely used mappings are the \textit{persistence landscape} \cite{bubenik2017persistence,bubenik2012statistical} and the \textit{persistence image} \cite{adams2017persistence} (we will focus on the latter). The \textit{persistence image} is a smoothed representation of the points ($x,y$) in a persistence diagram $(x,y) \in PD$. Typically, the smoothing is done by applying a Gaussian kernel (with mean $u$ and variance $\sigma^2$) to each of the points $\{x,y\} \in PD$: 
\begin{equation}
    \phi(x,y) = \frac{1}{2\pi\sigma^2}\exp -\frac{(x - u_x)^2 + (y - u_y)^2}{\sigma^2}
	\label{eq:kernel}
\end{equation}

\begin{definition}{\textbf{Persistence Surface}}: The \textit{persistence surface} is a scalar mapping $\rho: \mathbb{R}^2 \rightarrow \mathbb{R}$ with weighting function $w: \mathbb{R}^2 \rightarrow \mathbb{R}$ and is defined as:

\begin{equation}
    \rho(PD) = \sum_{u \in PD}w(u)\phi(u)
\end{equation}

\end{definition}

A \textit{persistence image} is the discretization of the \textit{persistence surface}.

\begin{definition}{\textbf{Persistence Image}: \hl{For a given PD, a \emph{persistence image} of size $n$ by $m$ is a collection of pixels supported in a rectangle $R = [a,b]\times [c,d]$. The ($i,j$) pixel $p_{i,j}$, is the area $[a + \frac{i-1}{n}(b - a), a+\frac{i}{n}(b-a)] \times [c + \frac{j-1}{m}(d-c), c + \frac{j}{m}(d - c)]$.}}

\begin{equation}
    \overrightarrow{\rho}[i,j] = \iint_{p_{i,j}}\rho dx dy
\end{equation}

\end{definition}

The ultimate goal of vectorization methods is to create a \textit{stable} representation of the PD.

\begin{definition}{\textbf{Stability}}: A vector representation $\overrightarrow{PD}_{X_i}$ of a persistence diagram $PD_{X_i}$ is said to be \textit{stable} if small perturbations in $PD_{X_i}$ (represented as $PD_{X_i}'$) results \hl{in a bounded change} in $\overrightarrow{PD}_{X_i}'$. 
\end{definition}

Mathematically, we establish the stability property as:
\begin{equation}
    \textrm{dist}(\overrightarrow{PD}_{X_i},\overrightarrow{PD}_{X_i}') \leq L\cdot \textrm{dist}(PD_{X_i},PD_{X_i}') 
	\label{eq:imagestable}
\end{equation}

where $dist(\cdot,\cdot)$ represents a distance metric and $L$ represents a scalar constant. \hl{The distance between PDs is commonly measured using the \textit{Wasserstein distance}  or the \textit{bottleneck distance}. Whereas the distance between the vectorized PDs ($\overrightarrow{PD}$) can be expressed in terms of $l_p$ norms.}

\begin{definition}{\textbf{Wasserstein distance}: The $p^{th}$-\textit{Wasserstein distance} between persistence diagrams $PD_1$ and $PD_2$ is defined as:}

\begin{equation}
    d_{W_p}(PD_1,PD_2) = \left(\inf_{\gamma}\sum_{x\in PD_1}||x - \gamma(x)||_{\infty}^{p}\right)^{1/p}
\end{equation}

where $\gamma$ ranges over all possible bijections from $PD_1$ to $PD_2$. \hlp{By convention, we add an infinite number of points in the diagonal to allow for bijections between PD's containing difference numbers of points.} 

\end{definition}

\begin{definition}{\textbf{Bottleneck distance}: The \textit{bottleneck distance} between persistence diagrams, $PD_1$ and $PD_2$, is defined as:}

\begin{equation}
    d_{B}(PD_1,PD_2) = \inf_{\gamma}\sup_{x}||x - \gamma(x)||_{\infty}
\end{equation}

where $\gamma$ ranges over all possible bijections from $PD_1$ to $PD_2$ \hlp{with the same considerations made in the Wasserstein distance definition.}

\end{definition}

Notably, it has been proven that \textit{persistence landscapes} and \textit{persistence images} are stable under the \hl{appropriate distances} \cite{bubenik2017persistence,bubenik2012statistical,adams2017persistence}. Intuitively, stability indicates that topology changes in a {\em continuous} manner under perturbations. This makes them excellent representations of PD and amenable to use in diverse tasks such as regression and classification, as we  demonstrate in Section \ref{section:examples}. 

\subsection{Topology of Continuous Functions}
\label{sec:contfunc}
\begin{figure}
    \centering
    \includegraphics[width = 0.8\textwidth]{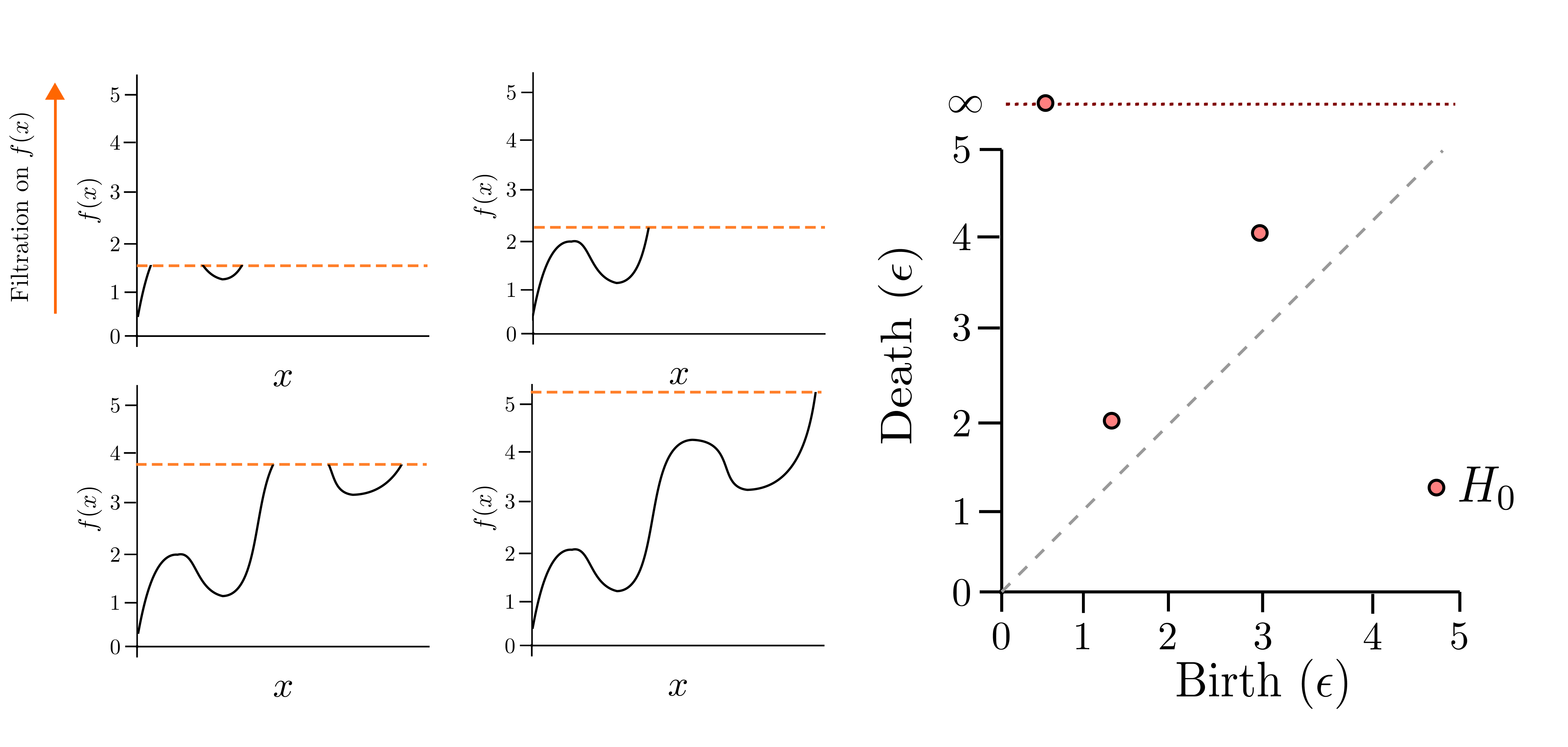}
    \caption{\hl{Morse filtration of a function $f: \mathbb{R} \rightarrow \mathbb{R}$. Consider a sublevel set $f^{-1}(-\infty,b)$ for increasing values of $b$. The topology of these sublevel sets changes at the critical points of $f$. As $b$ increases the value of the sublevel set, the persistence diagram records the topological changes in the function.}}
    \label{fig:func_pers}
\end{figure}

In the previous sections we focus on understanding the topology of datasets that are made of \hl{point clouds}. We now discuss how to quantify the topology of continuous functions. An example of this type of object is the scalar function shown Figure \ref{fig:func_pers}. Topologically, the interesting features of this function are its critical points (min and max points). In order to characterize these critical points we utilize a new form of filtration known as a \textit{Morse filtration} or \textit{sublevel set} filtration \cite{mischaikow2013morse}. The Morse filtration is derived from ideas of Morse theory, which is the study of the topology of manifolds through differential functions and the analysis of the critical points of these functions \cite{milnor1969morse}. Here, we consider the graph of a continuous function as a differentiable manifold \cite{poincare1895analysis}.

\begin{definition}{\textbf{Level Set}: Given a differentiable manifold $M$ and function $f:M \rightarrow \mathbb{R}$, the \textit{level set} $M^{a}$ at a point $a$ is defined as the pre-image:}
\begin{equation}
    M^a = f^{-1}(a) = \{x \in M : f(x) = a\}.
\end{equation}

\end{definition}

The level set contains all points of the manifold that have the same function value. In order to create a filtration, we use a \emph{sublevel set} defined below.

\begin{definition}{\textbf{Sublevel Set}: A \textit{sublevel set} $M^{(a,b)}$, where $a,b \in \bar{R}$ and $a = -\infty$ is defined as:}

\begin{equation}
    M^{(-\infty,b)} = \{x\in M: f(x) \leq b\}.
\end{equation}

\end{definition}

As we pass through the Morse filtration and build the sublevel set of the function, we are creating a well-defined filtered complex. The topology of the function will change as the filtration passes through critical points in the manifold \cite{mischaikow2013morse}. These topological changes are quantified in a persistence diagram \hl{which is subsequently vectorized for analysis}. An example of this type of filtration and the \hl{corresponding} persistence diagram are shown in Figure \ref{fig:func_pers}. The Morse filtration is a good choice for functions that are continuous or that can be approximated as piece-wise \hl{linear} functions (e.g., a time series). The method can be expanded to $k$-dimensional functions  \cite{gunther2012efficient}. This makes it a powerful approach to characterize complex surfaces (landscapes) that have many minima/maxima. We demonstrate this technique on 2D and 3D functions in Section \ref{section:examples}.

\subsubsection{Stability of Persistence Diagrams for Functions}

Persistence diagrams of real valued functions are also \textit{stable} representations of data. The following Theorem \ref{thrm:stab}, established by Cohen-Steiner and co-workers, highlights this result \cite{cohen2007stability}.

\begin{figure}[!htp]
    \centering
    \includegraphics[width = 0.6\textwidth]{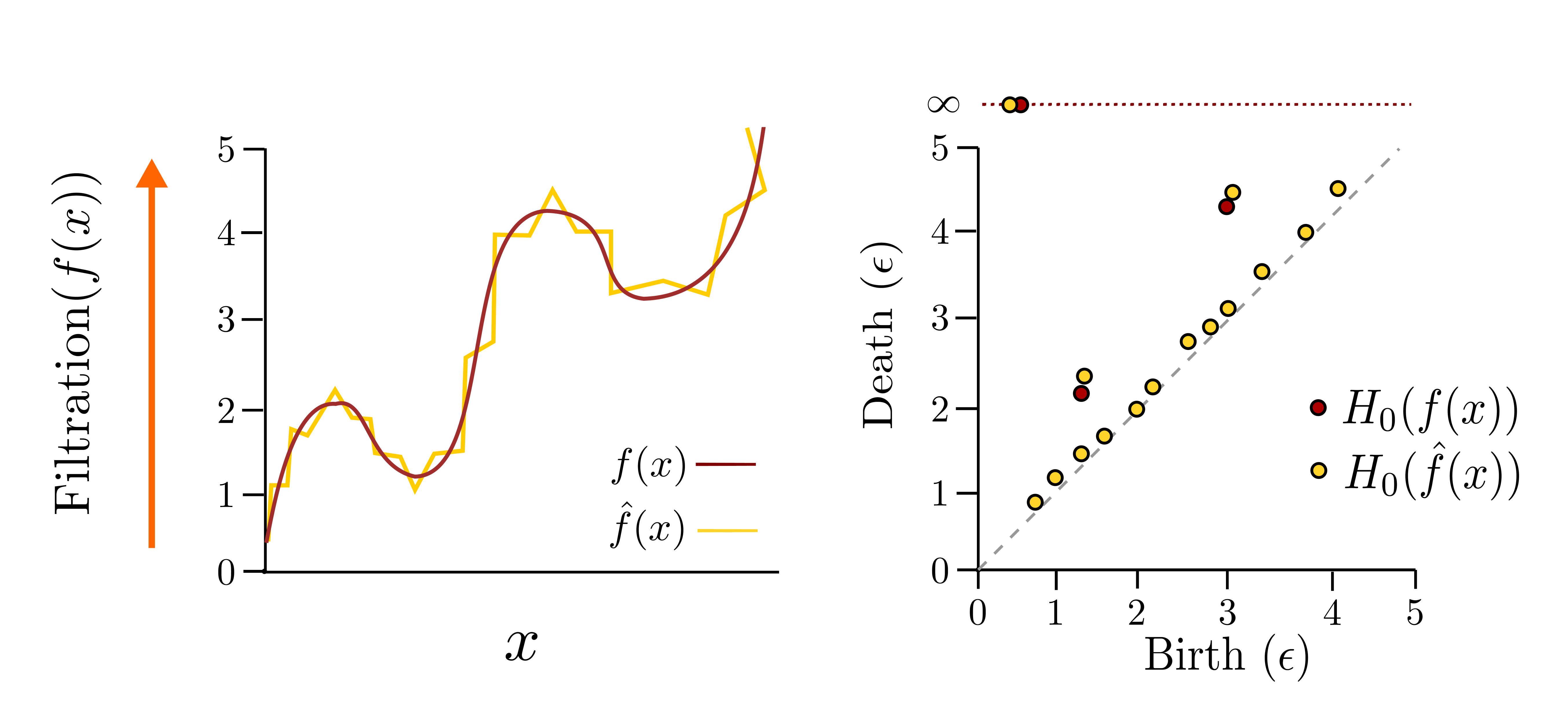}
    \caption{Filtration of functions $f(x)$ and $\hat{f}(x)$. The strong critical points of both functions are captured in the persistence diagram analysis while the weak critical points (arising from noise) remain close to the diagonal as they have minimal persistence.}
    \label{fig:points_pers2}
\end{figure}

\begin{theorem}{Given real valued functions $f,g$ with finitely many critical points and their corresponding persistence diagrams $PD_f,PD_g$, we have that:}

\begin{equation}
    d_{B}(PD_f,PD_g) \leq ||f - g||_\infty
\end{equation}

\hl{where $||\cdot||_{\infty}$ represents the $l_\infty$ distance between two functions $f,g : X \rightarrow \mathbb{R}$:}

\begin{equation}
    ||f - g||_\infty = \sup\{||f(x)-g(x)|| : x\in X\}.
\end{equation}

\label{thrm:stab}
\end{theorem}
An illustration of the stability of persistence diagrams is shown in Figure \ref{fig:points_pers2}. \hl{Here, we can see that the persistence diagram of the two functions are similar. The presence of the strong critical points, with large persistence, are well captured in both diagrams. We can also see that the persistence diagram captures the structure of the weak critical points (arising from noise), which have short persistence.}

\begin{figure}[!htp]
\begin{subfigure}{.5\textwidth}
  \centering
  \includegraphics[width=1\linewidth]{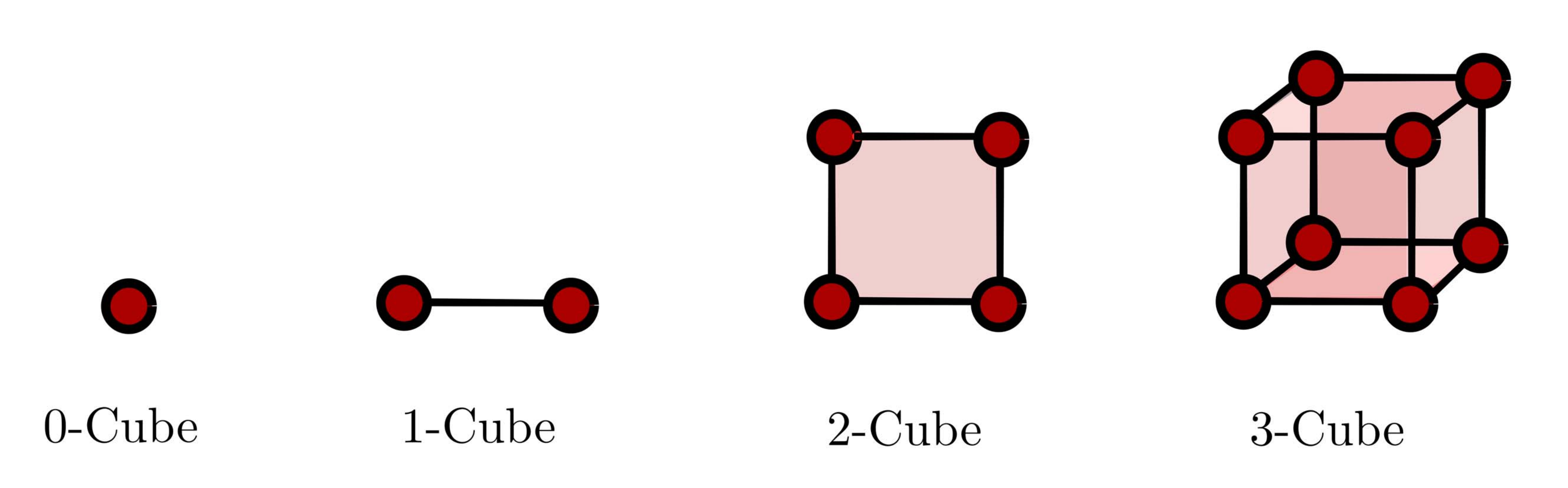}  
  \caption{}
  \label{fig:cubes}
\end{subfigure}
\begin{subfigure}{.5\textwidth}
  \centering
  \includegraphics[width=.6\linewidth]{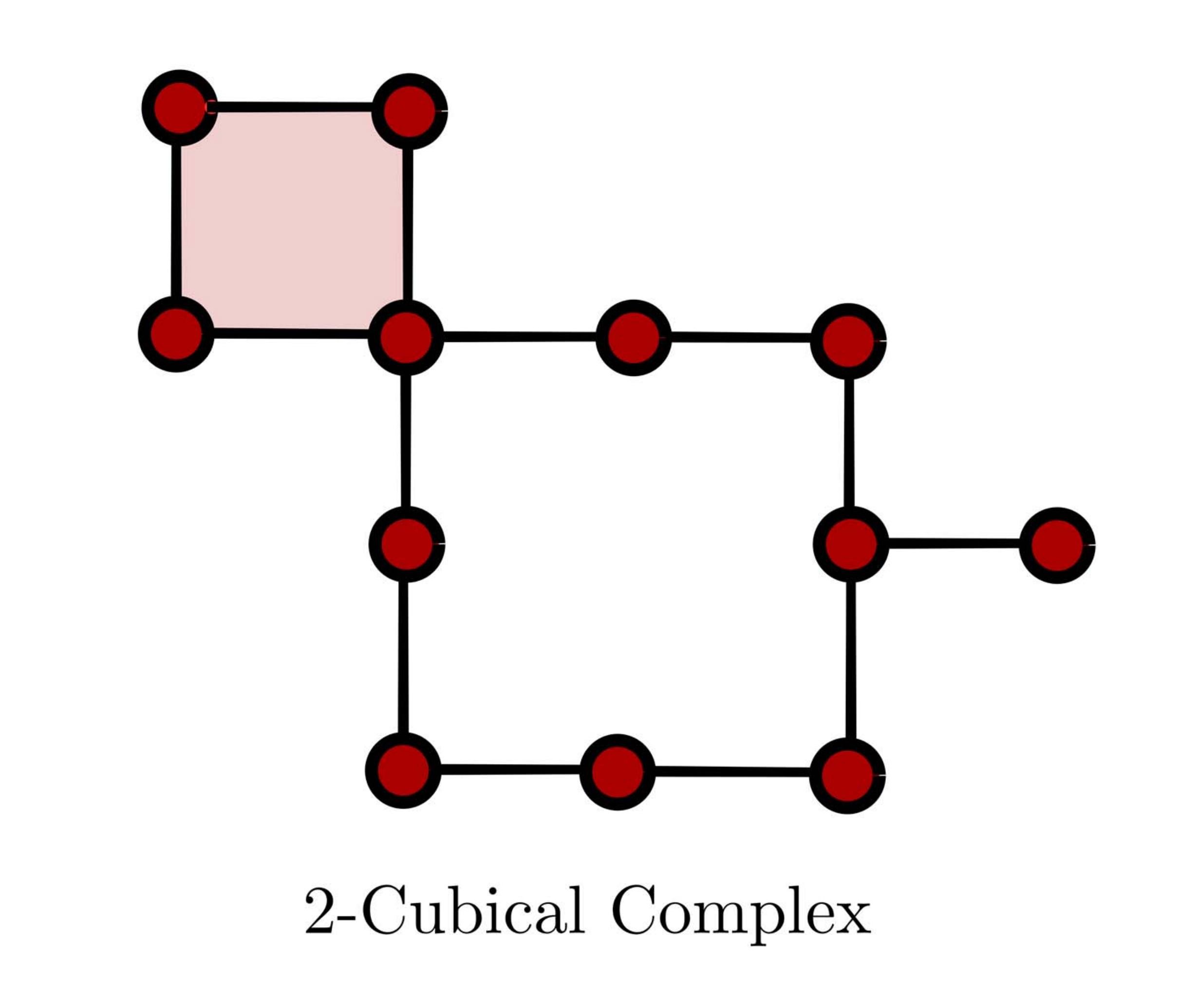}  
  \caption{}
  \label{fig:cubical_comp}
\end{subfigure}
\caption{(a) Representations of the elementary $k$-cubes of dimension $0$ to dimension $3$. (b) A $2$-cubical complex that contains a $2$-cube, along with a hole created by $1$-cubes.}
\label{fig:concompex2}
\end{figure}

\subsection{Cubical Complexes and Images}
\label{sec:cubicalcomp}
We briefly discuss cubical complexes and cubical homology as these are important in understanding data over rectangular domains such as images (represented by pixels in 2D or voxels in 3D) and are used primarily in Morse filtrations \cite{allili2001cubical,niethammer2002analysis}.  Cubical homology is similar to simplicial homology but uses different basis shapes. \hl{This highlights the fact that one can choose different basis shapes in topological analysis. The basis shapes for cubical complexes are $k$-cubes (hypercubes) or elementary cubes \hlp{built from sides of equal length}. In Figure \ref{fig:cubes} we show $k$-cubes of dimension $0$ through $3$. An example cubical complex is found in Figure \ref{fig:cubical_comp}.}

\begin{figure}[!htp]
\begin{subfigure}{.4\textwidth}
  \centering
  \includegraphics[width=.7\linewidth]{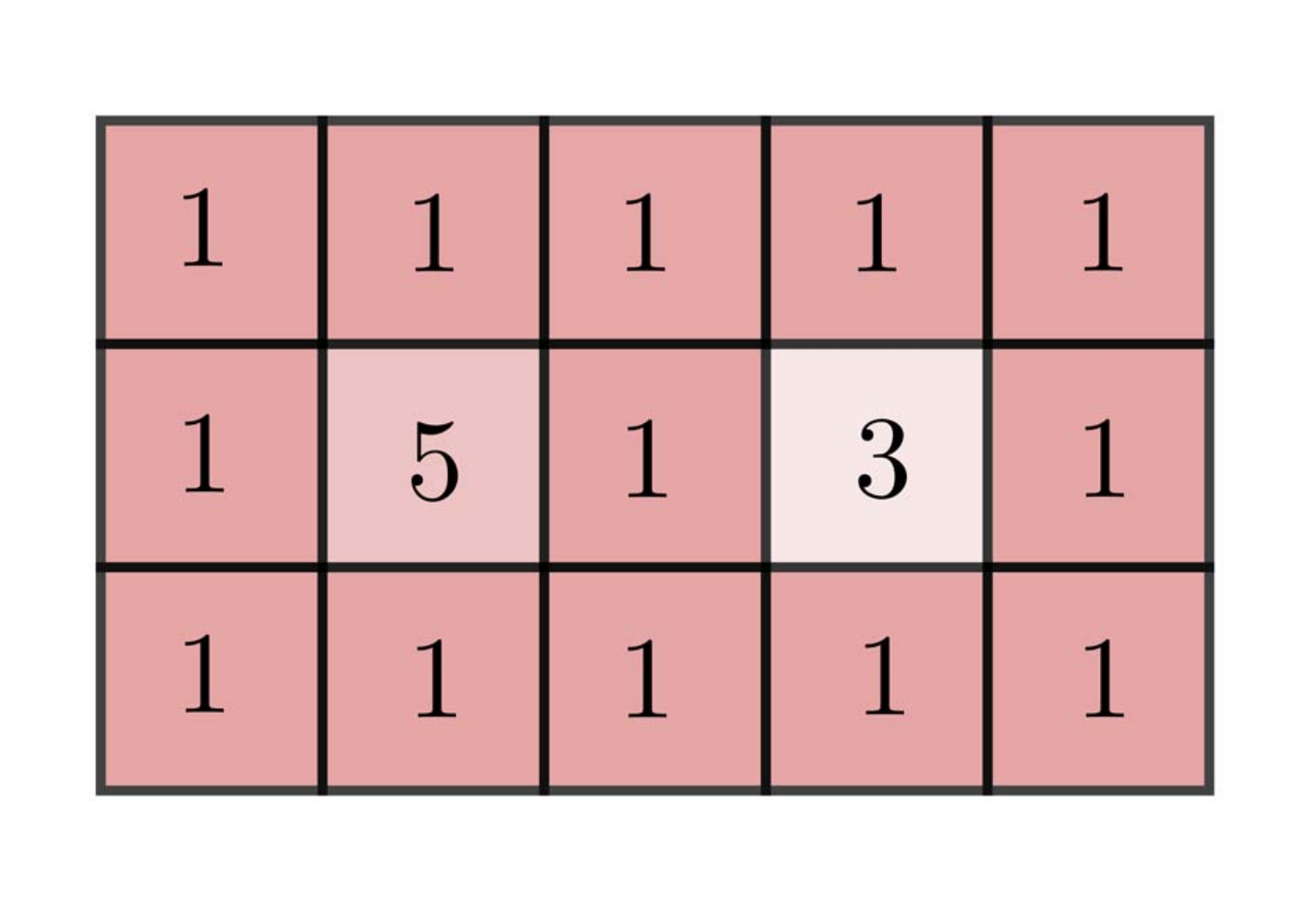}  
  \caption{}
  \label{fig:cubeimg}
\end{subfigure}
\begin{subfigure}{.6\textwidth}
  \centering
  \includegraphics[width=1\linewidth]{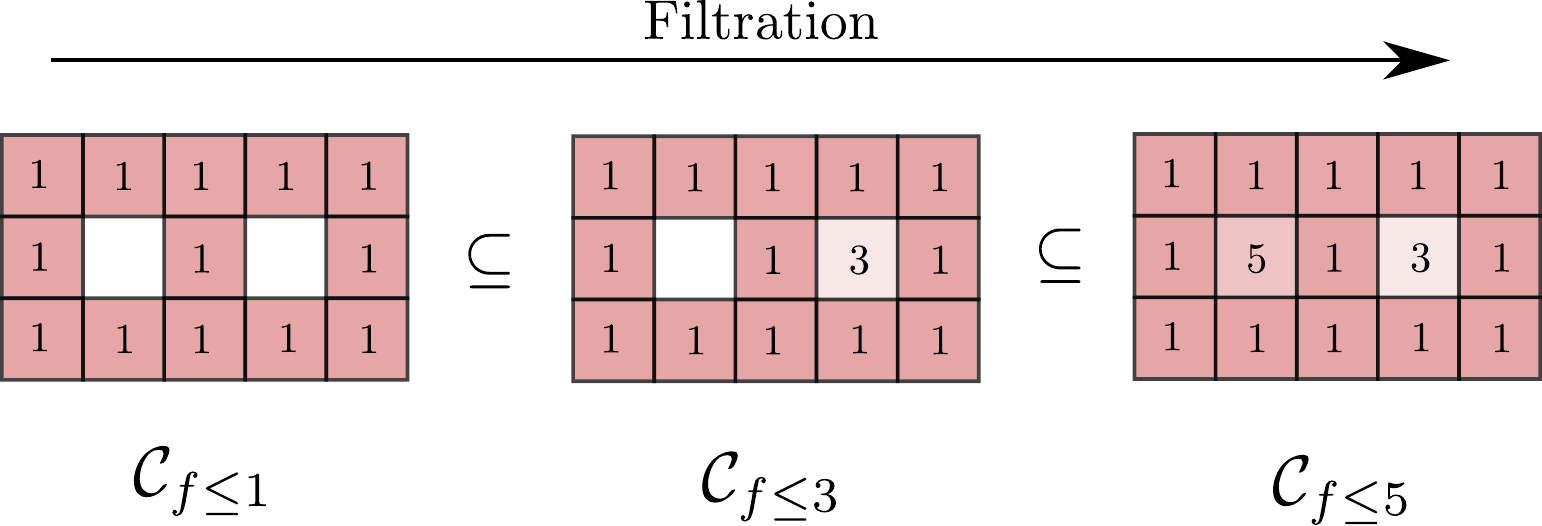}  
  \caption{}
  \label{fig:cube_filt}
\end{subfigure}
\caption{\hl{(a) Representation of an image; the values within each pixel represent intensity. (b) The representative filtration on the image itself and the corresponding sublevel sets.}}
\end{figure}

\begin{figure}[!htp]
    \centering
    \includegraphics[width = 0.9\textwidth]{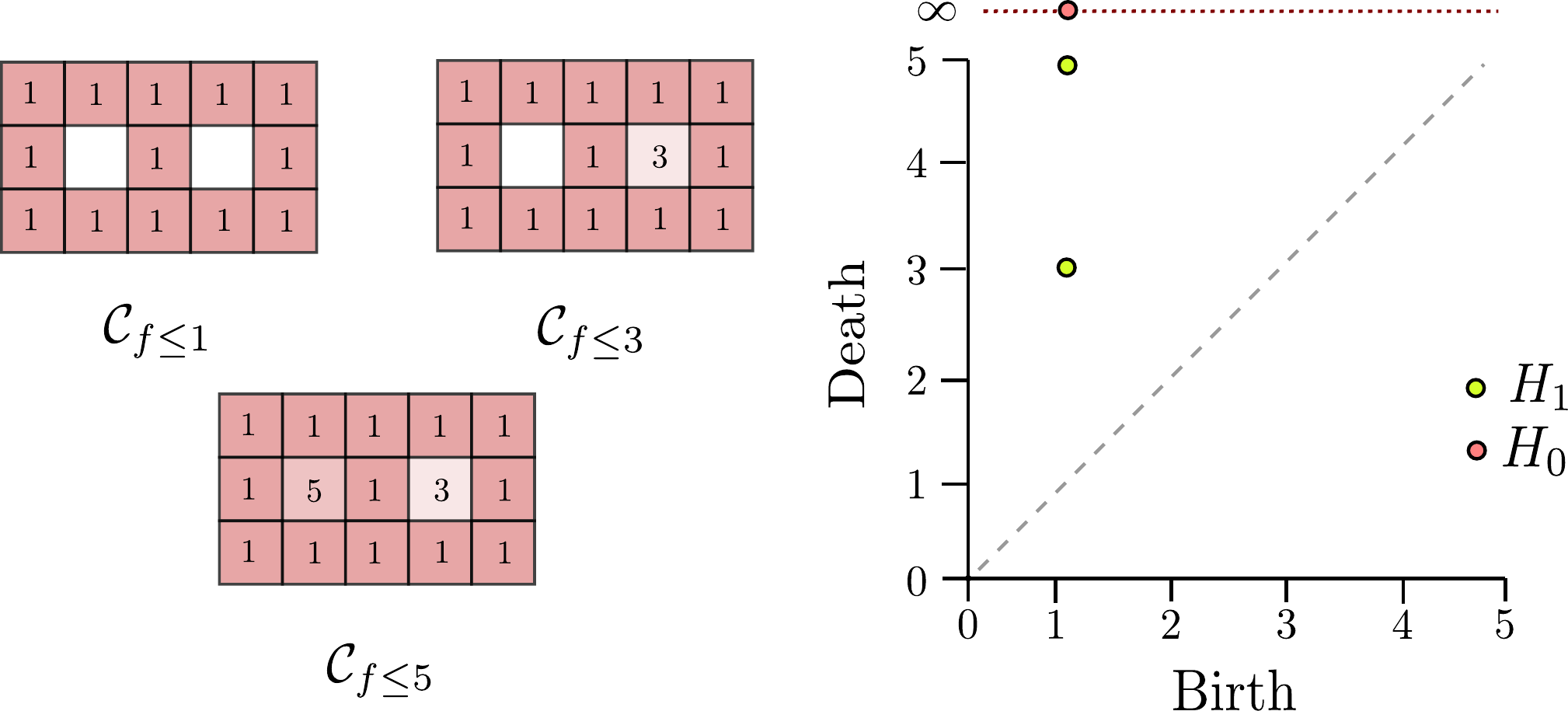}
    \caption{Persistence diagram for the filtration of the image in Figure \ref{fig:cubeimg}. The PD reveals the structure of the image,which changes at the critical points of the image (i.e., $f=3$ and $f=5$.)}
    \label{fig:cube_diagram}
\end{figure}

Filtered cubical complexes can be developed from data and we can perform persistence homology calculations on them, similar to filtered simplicial complexes.  An important application of cubical analysis is the analysis of images. An image can be viewed as a 2D surface embedded in three dimensions where two dimensions are the coordinates of each pixel, and the third dimension is the scalar value associated with each pixel (e.g., intensity). We can also view an image as an approximation of a continuous object over which a Morse filtration can be performed. We demonstrate the concept of filtration on the image in Figure \ref{fig:cubeimg}. We filter through the level sets of this image and develop the filtered complex through the sublevel set. This filtration is demonstrated in Figure \ref{fig:cube_filt}, where  $\mathcal{C}_{f \leq 1} \subseteq \mathcal{C}_{f \leq 3} \subseteq \mathcal{C}_{f \leq 5}$. The persistence diagram is then represented in Figure \ref{fig:cube_diagram}.  The PD  is able to capture the dominant topological features of our image such as the presence of the two critical points in the image and the fact that there is only a single connected component. The persistence diagram generated from a cubical complex filtration and a simplicial complex filtration have the same properties and can be vectorized in the same way. Thus, we can extend our analysis from discrete point data to that of images or other high-dimensional continuous objects.

\subsection{Inverse Analysis}

\begin{figure}[!htp]
    \centering
    \includegraphics[width = 0.7\textwidth]{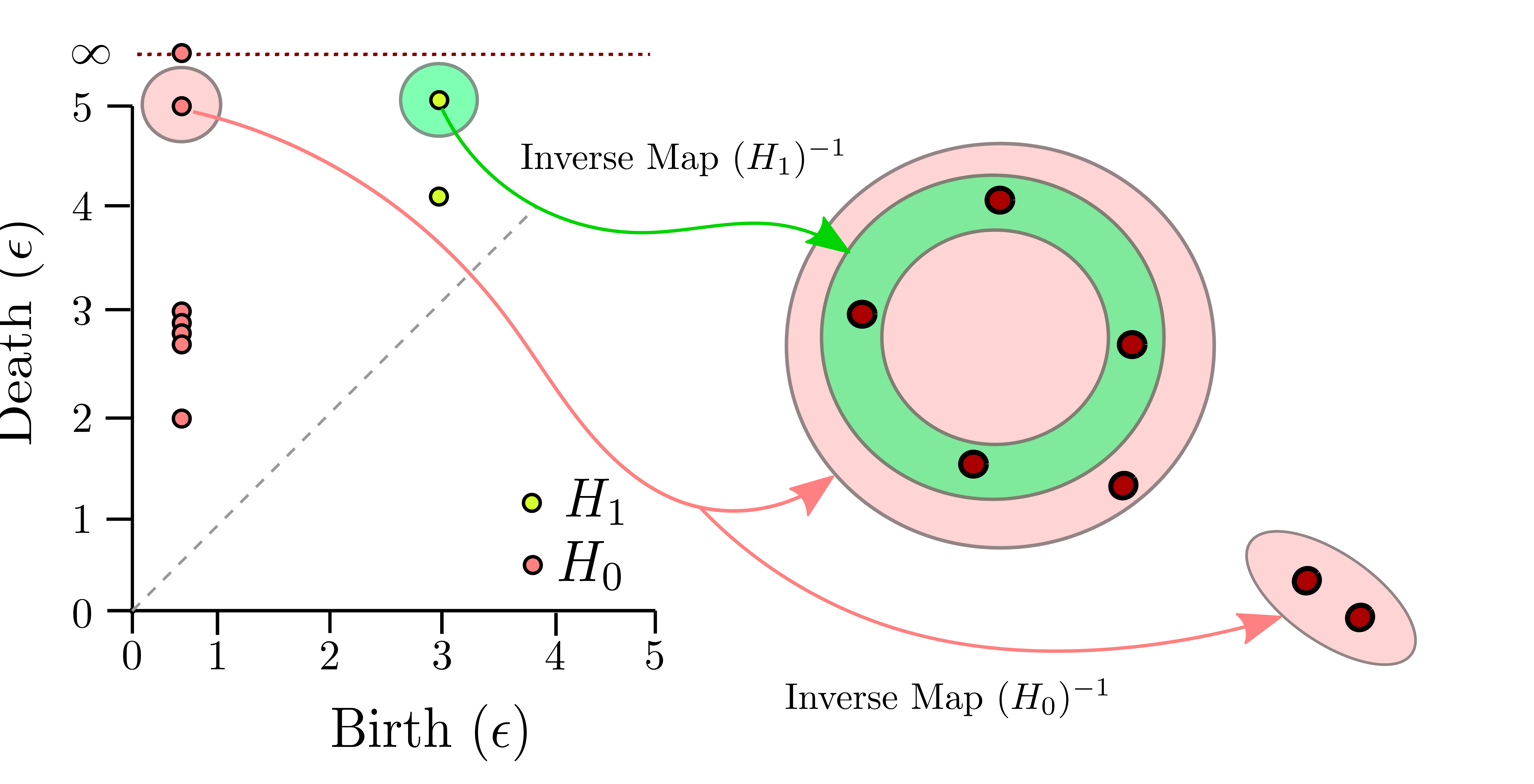}
    \caption{\hl{Inverse mapping of the features of persistence diagram to the features dataset in Figure \ref{fig:points_pers}. The features of interest are picked up within the persistence diagram; In this case they are the highly persistent hole and connected component. The inverse mapping for the selected point in one dimensional persistence identifies the dominant loop in the dataset and the inverse mapping for the selected point in zero dimensional persistence identifies the two most distinct clusters of the data.}}
    \label{fig:inverse}
\end{figure}

\hl{The algorithm to compute persistence homology allows us, for each point in a persistence diagram, to identify representative features in the filtration. The intuition behind this idea is presented in Figure \ref{fig:inverse}. However, before we discuss this procedure along with an appropriate post-processing step, a few issues need to be clarified. The persistence algorithm will provide a set of cycles that \emph{generate} a persistent homology group. However, it will not necessarily be the set that is close to what we may call \emph{geometrically optimal generators}. A similar problem can occur when identifying an optimal basis for a vector space. For example, let us consider the $\mathbb{R}^2$ space and various possible bases for this space presented below.}

\[
B_1= 
\Bigg\{
\begin{bmatrix}
    1     \\
    0       
\end{bmatrix}, 
\begin{bmatrix}
    0     \\
    1       
\end{bmatrix}
\Bigg\},
B_{2}= 
\Bigg\{
\begin{bmatrix}
    1     \\
    1       
\end{bmatrix}, 
\begin{bmatrix}
    0     \\
    10       
\end{bmatrix}
\Bigg\},
B_{3}= 
\Bigg\{
\begin{bmatrix}
    123     \\
    21       
\end{bmatrix}, 
\begin{bmatrix}
    2121     \\
    1243       
\end{bmatrix}
\Bigg\}.
\]

\hl{$B_1, B_2$ and $B_3$ are all valid basis of $\mathbb{R}^2$ and they can be mapped to each other via multiplication by a non-singular matrix. Yet, only $B_1$ can be called \emph{geometrically optimal} as it the natural representation of the space. }

\begin{figure}[!htp]
    \centering
    \includegraphics[width = 0.7\textwidth]{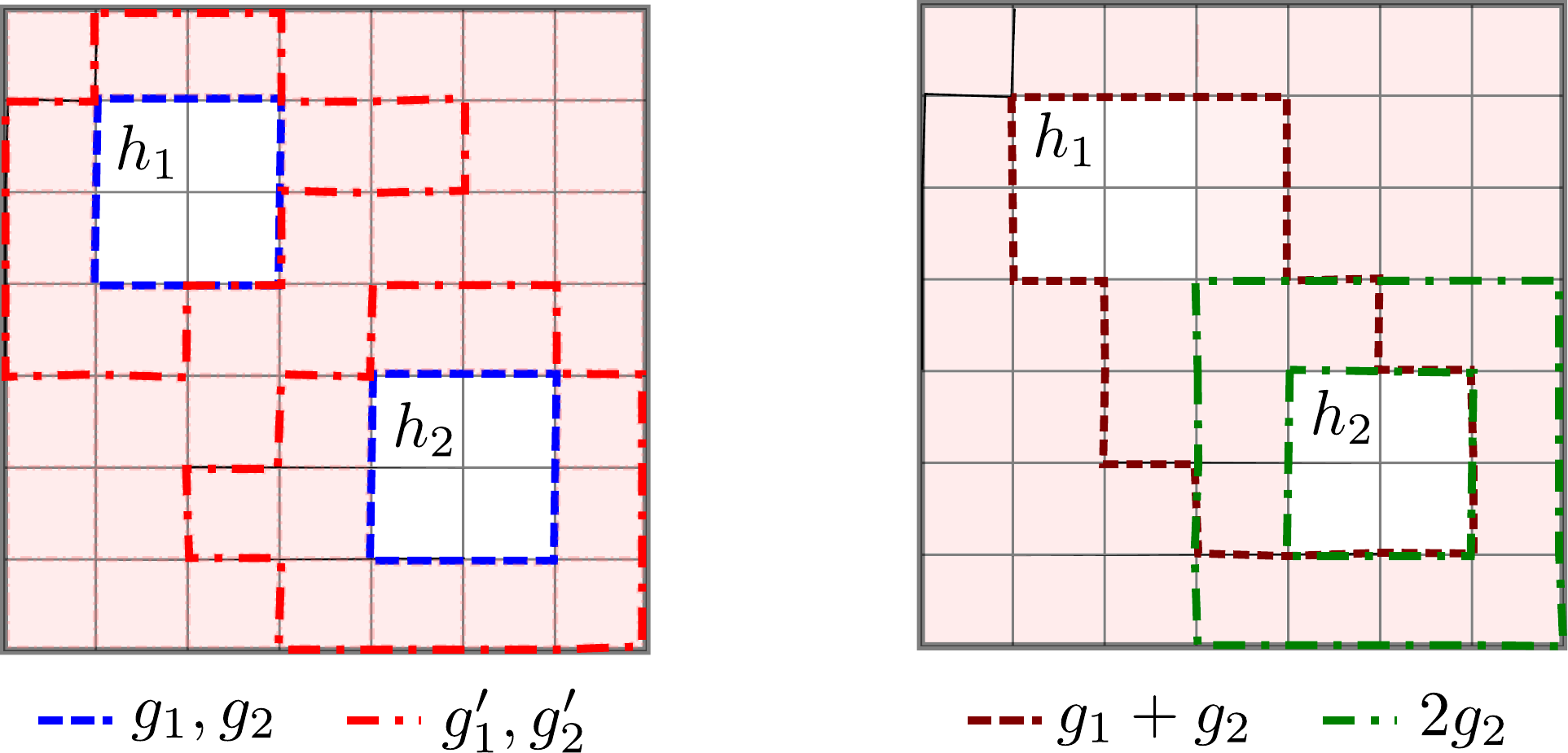}
    \caption{\hl{A representation of a cubical complex with two holes. (left) Cycles $g_1$ and $g_2$ are an optimal representation of holes ($h_1,h_2$) as they trace each hole. Cycles $g'_1$ and $g'_2$  also represent these holes, but are not geometrically optimal. (right) Another set of generators for the holes ($h_1,h_2$) that are not optimal as they represent linear combinations of the cycles $g_1$ and $g_2$. }}
    \label{fig:opt_bases}
\end{figure}

\hl{We now illustrate, with a simple example, how this issue manifests itself when selecting a basis for a persistent homology group. The cubical complex presented in Figure \ref{fig:opt_bases} contains two holes ($h_1,h_2$). Cycles $g_1$ and $g_2$ surrounding $h_1$ and $h_2$ are generators (they form a basis) for the first homology group of this shape. They are the optimal geometric representation of these holes. Let $\mathbf{G_{b1}}$ be the basis made up of these two cycles.}
\hl{The second basis $\mathbf{G_{b2}}$ contains cycles $g_1'$ and $g_2'$. They are a perturbed version of cycles $g_1$ and $g_2$, as $g_1  \simeq g_1'$ and $g_2  \simeq g_2'$, but are not geometrically optimal.}

%
\hl{Lastly, let us consider a basis $\mathbf{G_{b3}}$. It generates the first homology group of our shape, but it is a non-optimal basis for the holes because there is no unique correspondence between cycles and holes. $\mathbf{G_{b3}}$ and $\mathbf{G_{b1}}$ can be transformed into each other by multiplication of a nonsingular matrix (change of basis).}

\[
\mathbf{G_{b1}}= 
\begin{bmatrix}
    g1  &  0      \\
    0  &  g2      
\end{bmatrix}, 
\mathbf{G_{b2}}= 
\begin{bmatrix}
    g1'  &  0      \\
    0  &  g2'      
\end{bmatrix},
\mathbf{G_{b3}}= 
\begin{bmatrix}
    g1  &  0      \\
    g2 &  2g2      
\end{bmatrix} 
\]

\hl{The persistent homology algorithm may return non-optimal bases, like $\mathbf{G_{b3}}$. There is not much that can be done to fix this issue in general. However, there is a way of finding the most optimal generators within their homology class. This will allow, for instance, to simplify $\mathbf{G_{b2}}$ into $\mathbf{G_{b1}}$.}
\hl{In order to address the problem of non-optimal generators, we construct an integer optimization problem to obtain the sparsest representation of a given generator. We begin this optimization with the representative cycle $c$ obtained from the persistence algorithm and identify the sparsest chain of simplices from the homology class of $c$. For example, in Figure \ref{fig:inverse}, we wish to identify the basis (generator) for the persistent cycle identified in green. We construct a simplified representation for the dataset in Figure \ref{fig:inverse} where $\mathcal{K}_{\epsilon}$ represents the associated simplicial complexes with the specified level of filtration. Because we know the cycle is born at $\epsilon = 3$ and dies at $\epsilon = 5$ we only focus on the portion of the filtration up to the level $\epsilon = 5$. Given the simplices below that filtration level, we wish to identify the sparsest set of $1$-simplices $z = \sum_{i}\sigma_i$, where $\sigma_i \in \mathcal{K}$, that is homologous to the generator $c$. We denote $||z||_0$ as the cardinality of $z$. }

\begin{figure}[!htp]
    \centering
    \includegraphics[width = 0.7\textwidth]{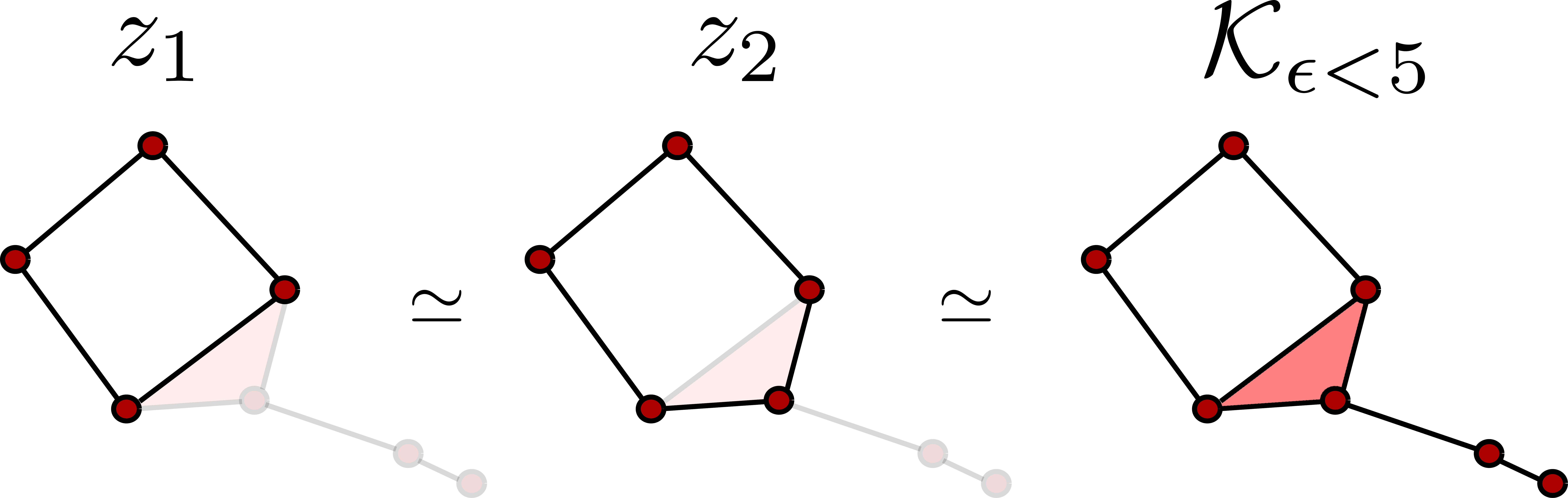}
    \caption{\hl{A representation of potential simplicial complexes that are homotopy equivalent ($\simeq$) to the filtration complex ($\mathcal{K})$. The complex $z_1$ contains the fewest number of $1$-simplices, thus is the optimal representation of the cycle contained in $\mathcal{K}_{\epsilon<5}$.}}
    \label{fig:inv_opt}
\end{figure}

\hl{With this information, we construct an integer optimization formulation, found in \ref{eq:optform}, where we seek to identify the sparsest set of simplices $z \in \mathcal{K}_{\epsilon<5}$ that is homologous to the cycle $c$ obtained from persistent homology computations. In Figure \ref{fig:inv_opt} we identify multiple cycles ($z_1,z_2 \in \mathcal{K}_{\epsilon<5}$), and note that $z_1 = z$ because $||z_1||_0 < ||z_2||_0$.}

\begin{equation}
z = \amin ||z||_0, \ \text{subject to} \  z \simeq c
\label{eq:optform}
\end{equation}

\hl{Inverse analysis can be used to understand dominant features that drive classification and regression results. For instance, the defining characteristics of a given class within a dataset can be identified via the weights of regression/classification models in the space of PDs (see Section \ref{section:examples}). These defining characteristics can then be mapped back to the original data set which can be used to gain physical understanding of differences between datasets}. Applications of these methods in material science are found in \cite{buchet2018persistent,nakamura2015persistent,ichinomiya2017persistent}. Here, inverse analysis is used to identify fracture or degradation sites in materials and to identify pore configurations in granular crystallization. 
\section{Applications}
\label{section:examples}

We now proceed to demonstrate how these methods can be applied to different types of datasets. To do so, we use a couple of illustrative examples and real datasets derived from soft materials and molecular dynamics simulations. All the calculations presented were conducted in \texttt{Python} \cite{oliphant2007python} using the TDA packages \texttt{GUDHI} \cite{maria2014gudhi} and \texttt{Homcloud} \cite{homcloud}. All scripts and data needed to reproduce these results can be found in \hl{\url{https://github.com/zavalab/ML/tree/master/TDApaper}.}

\subsection{Topology of Point Clouds}
\label{sec:ex1}

We illustrate how to use TDA to analyze point clouds; specifically, we seek to extract topological features from the data to perform binary classification of point clouds. The clouds used here are collections of points in two dimensions $x_1,x_2$ (for visualization purposes). In actual applications, one can conduct analysis on a point cloud of any dimension. We represent each cloud as $X_i$ where each cloud can belong to two different types of classes (Class 1 or Class 2). Our goal is to take each point cloud $X_i$ as an input, project this data to their respective $H_1$ persistence diagram $PD_{X_i}$ through an epsilon ball filtration (i.e., extract the topological features), vectorize the PDs ($\overrightarrow{PD}_{X_i} \in \mathbb{R}^{q}$), and perform classification of the cloud based on the vectorized topological features. 

\begin{figure}[!htp]
\begin{subfigure}{.5\textwidth}
  \centering
  \includegraphics[width=0.7\linewidth]{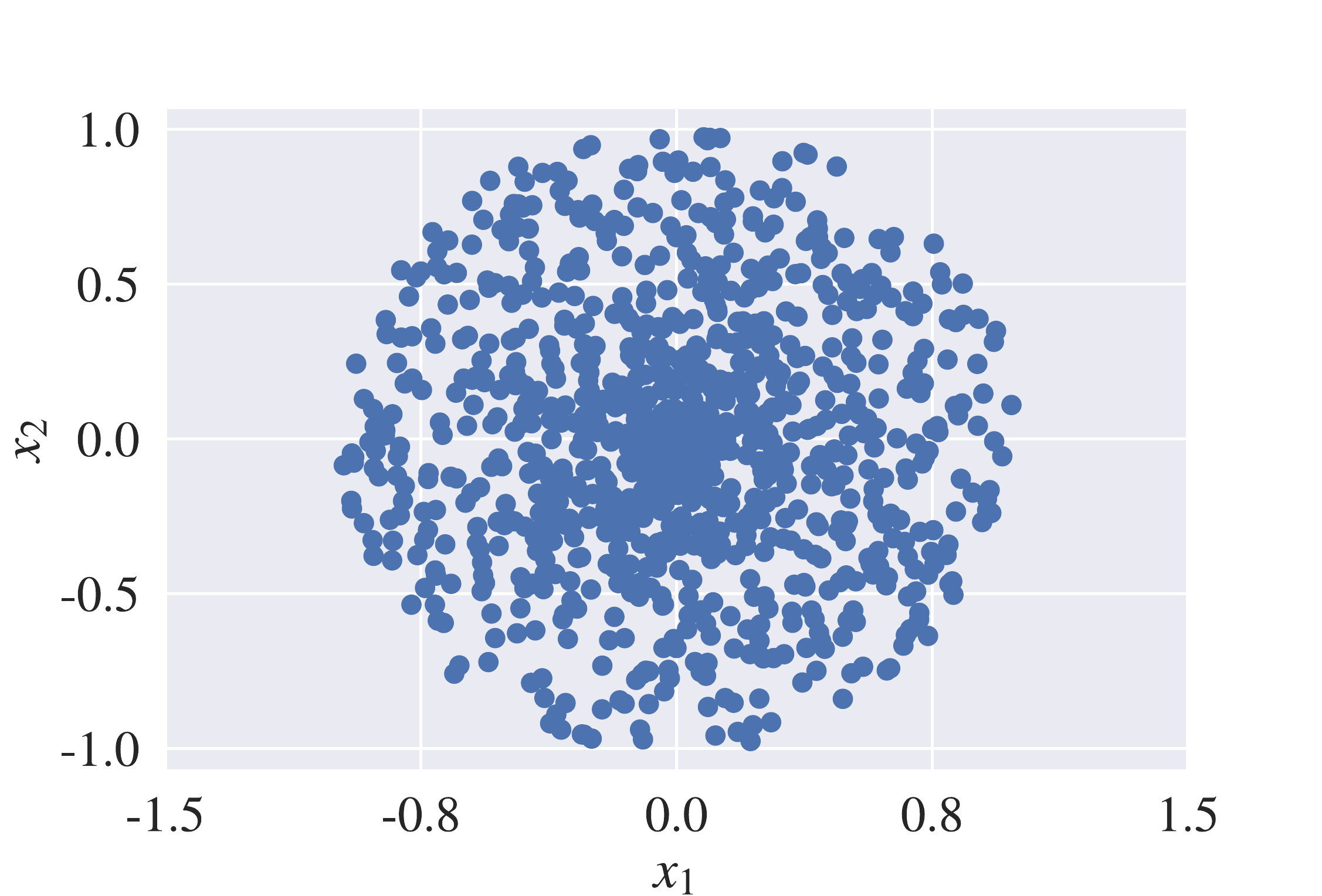}  
  \caption{Sample point cloud of Class 1.}
  \label{fig:sub-first}
\end{subfigure}
\begin{subfigure}{.5\textwidth}
  \centering
  \includegraphics[width=0.7\linewidth]{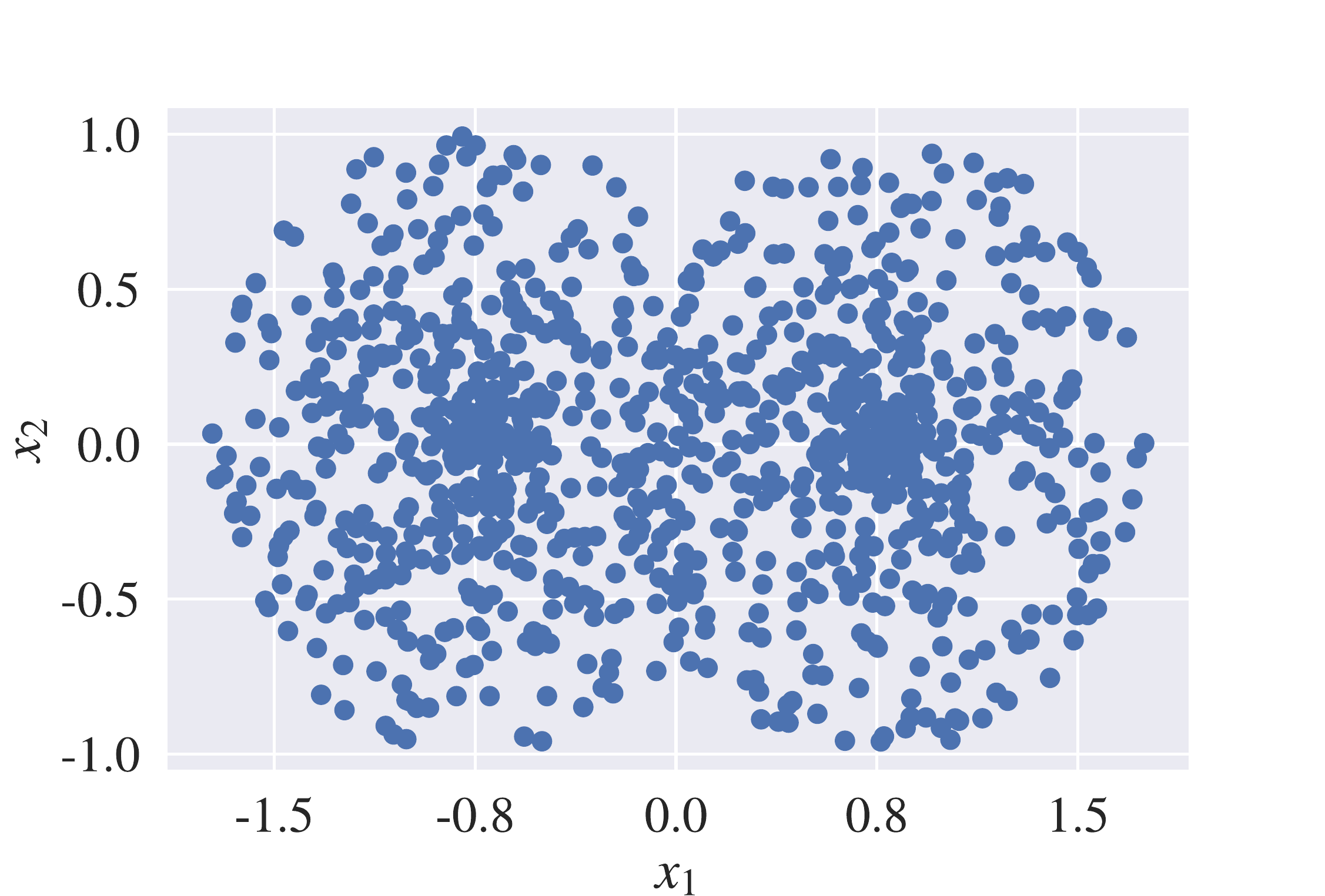}  
  \caption{Sample point cloud of Class 2.}
  \label{fig:sub-second}
\end{subfigure}
\caption{Types of point clouds analyzed using epsilon ball filtration.}
\label{fig:pclouds}
\end{figure}

\begin{figure}[!htp]
\begin{subfigure}{.33\textwidth}
  \centering
  \includegraphics[width=.7\linewidth]{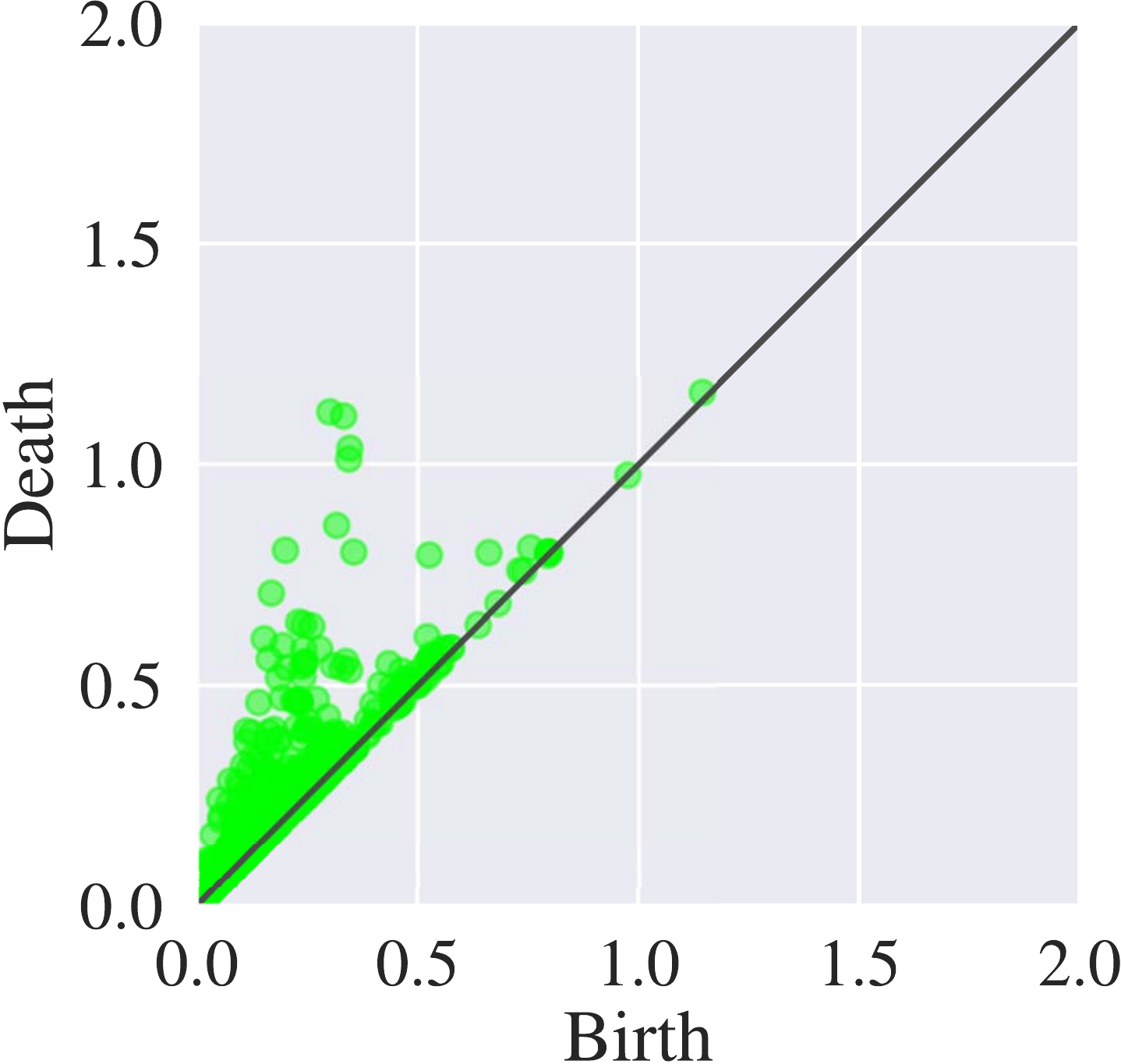}  
  \caption{}
  \label{fig:sub-first}
\end{subfigure}
\begin{subfigure}{.33\textwidth}
  \centering
  \includegraphics[width=.7\linewidth]{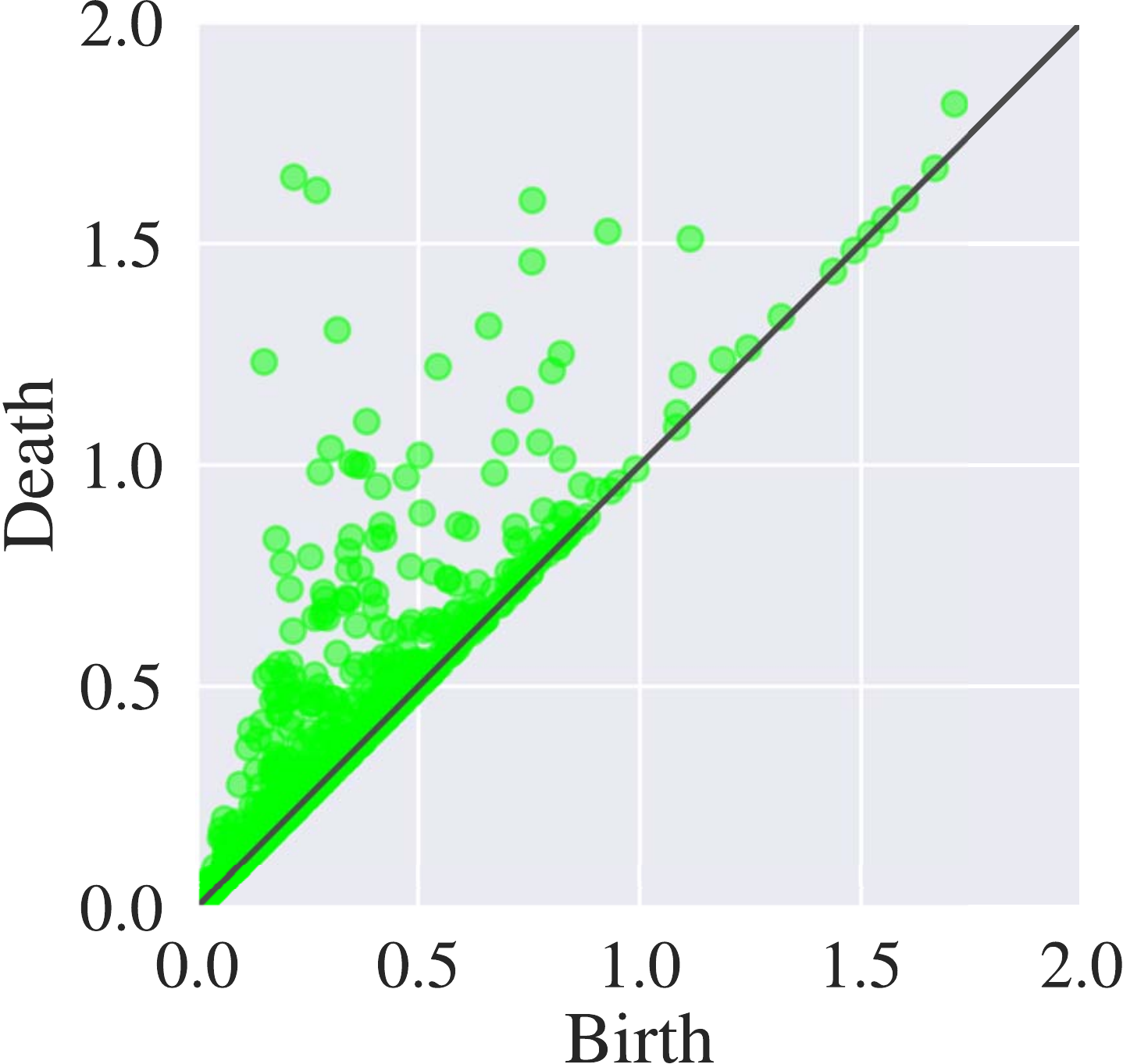}  
  \caption{}
  \label{fig:sub-second}
\end{subfigure}
\begin{subfigure}{.33\textwidth}
  \centering
  \includegraphics[width=1\linewidth]{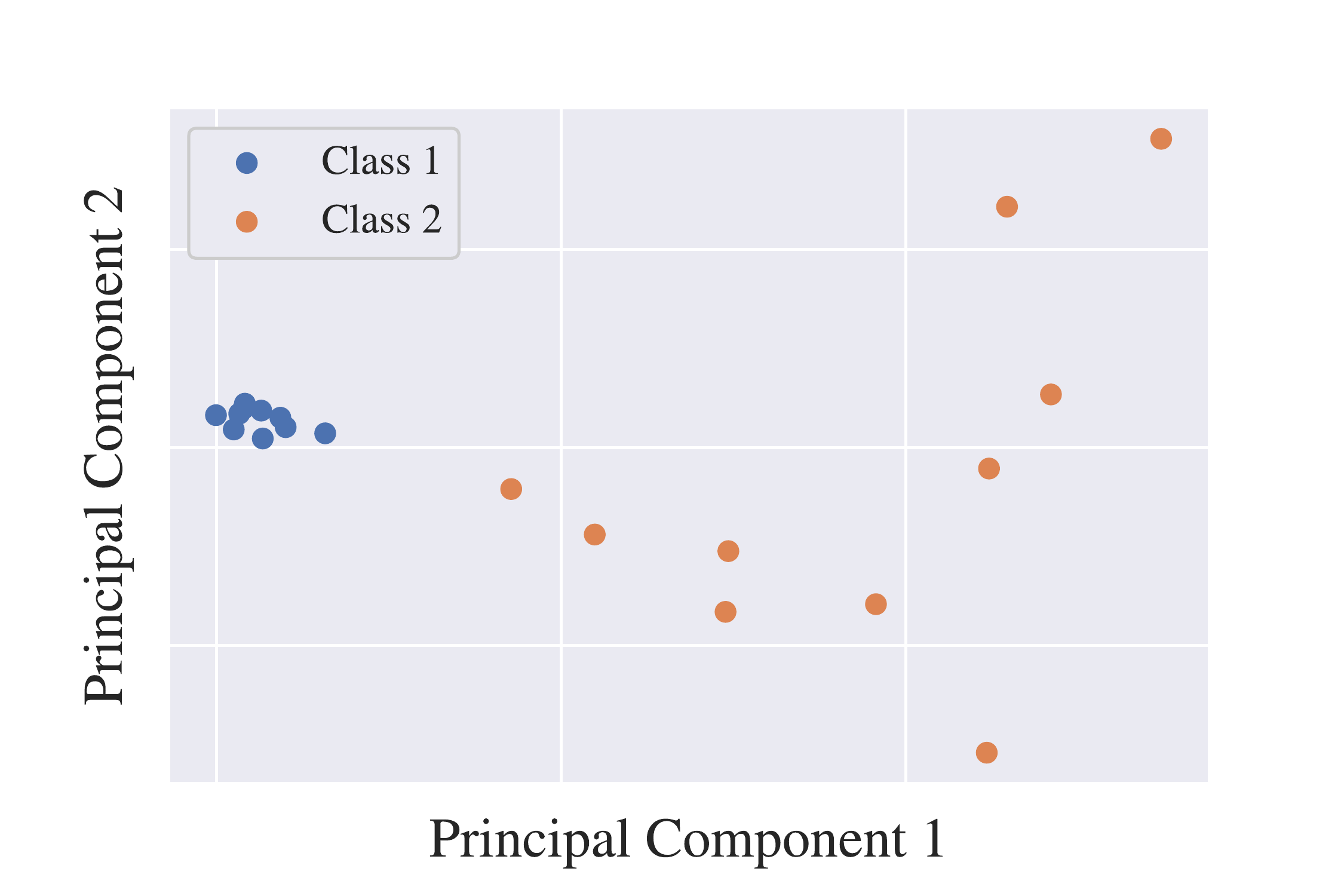}  
  \caption{}
  \label{fig:pcloudpca}
\end{subfigure}
\caption{Persistence diagrams for Class 1 and Class 2 point clouds and the corresponding \hlp{PCA analysis on the set of persistence images created from Class 1 and Class 2}. (a) $H_1$ persistence diagram for Class 1. (b) $H_1$ persistence diagram for Class 2. (c) \hlp{Principal components of persistence images} for Class 1 and Class 2 datasets. It is clear that there is separation of the \hl{persistence diagrams}.}
\label{fig:H1pcloud}
\end{figure}

The point cloud classes are shown in Figure \ref{fig:pclouds} and the $H_1$ PDs are found in Figure \ref{fig:H1pcloud}. Note how the point clouds of Class 1 define a simple object (ellipse) while those of Class 2 define a more complex object (overlapping ellipses). We can see that the \hl{persistence diagrams} are visually distinct; specifically, point clouds of Class 2 have features that persist over a longer range of the filtration. 

We utilize the persistence image method to vectorize the PDs. We apply principal component analysis (PCA) to the vectorized PDs to verify if there is a separation between Class 1 and Class 2. We emphasize that the PCA projection is not applied to the original datasets $X_i$ but to the transformed datasets $\overrightarrow{PD}_{X_i}$ obtained from TDA.  The PCA projection on first and second principal components is shown in Figure \ref{fig:pcloudpca}. This shows an obvious separation between Class 1 and Class 2, which means that the topological features extracted from the data (contained in $H_1$) are informative. This also suggests that a simple linear classifier using vectorized PDs as features should work well. To test this hypothesis, we apply a linear support vector machine (SVM) classifier using the $\overrightarrow{PD}_{X_i}$ as features; we find that we can {\em perfectly classify} the datasets (we use a $5$-fold cross validation scheme). This again indicates that the topological features extracted with TDA are highly informative. 

\begin{figure}[!htp]
\begin{subfigure}{.33\textwidth}
  \centering
  \includegraphics[width=1\linewidth]{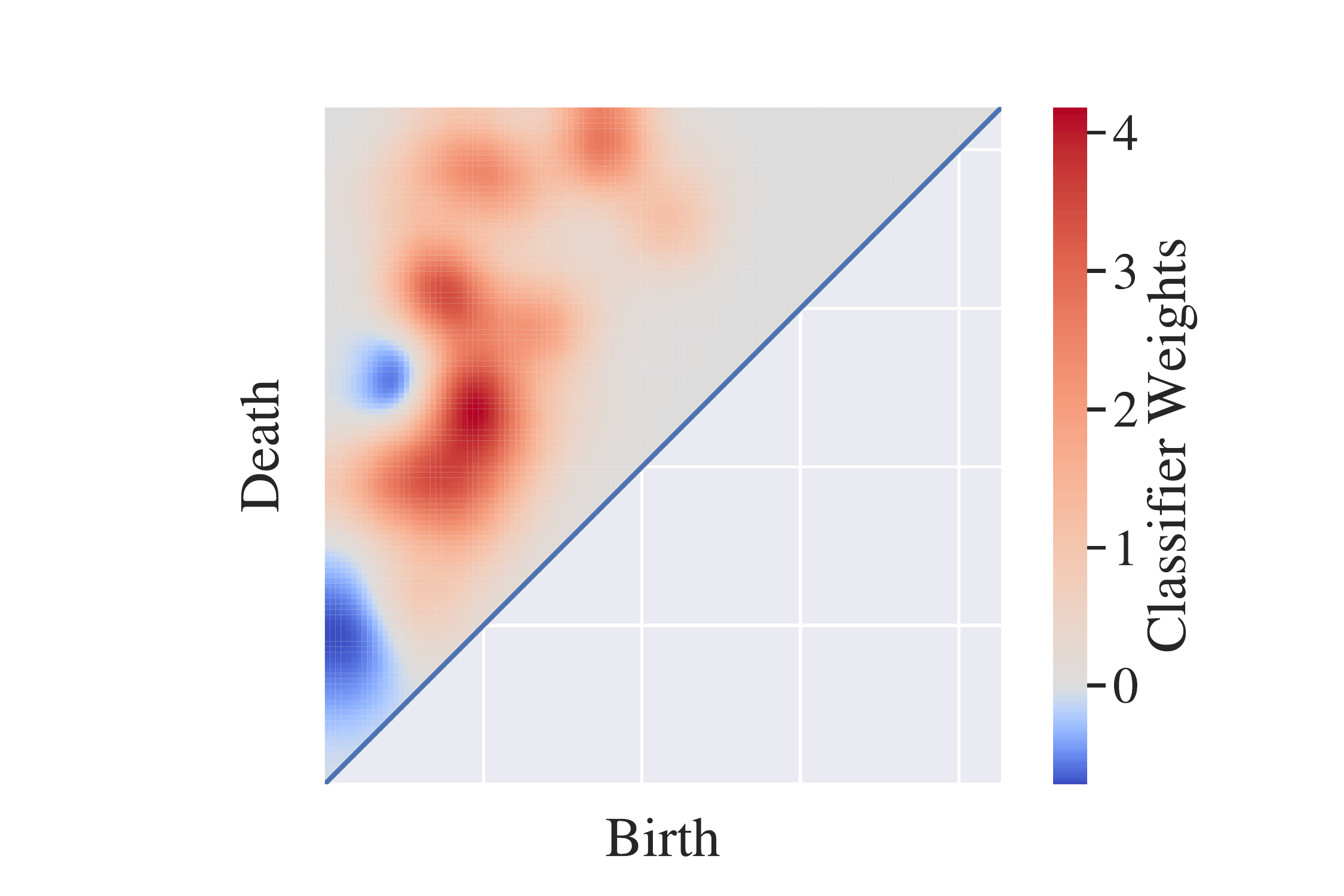}  
  \caption{}
  \label{fig:pcloudweights}
\end{subfigure}
\begin{subfigure}{.33\textwidth}
  \centering
  \includegraphics[width=1\linewidth]{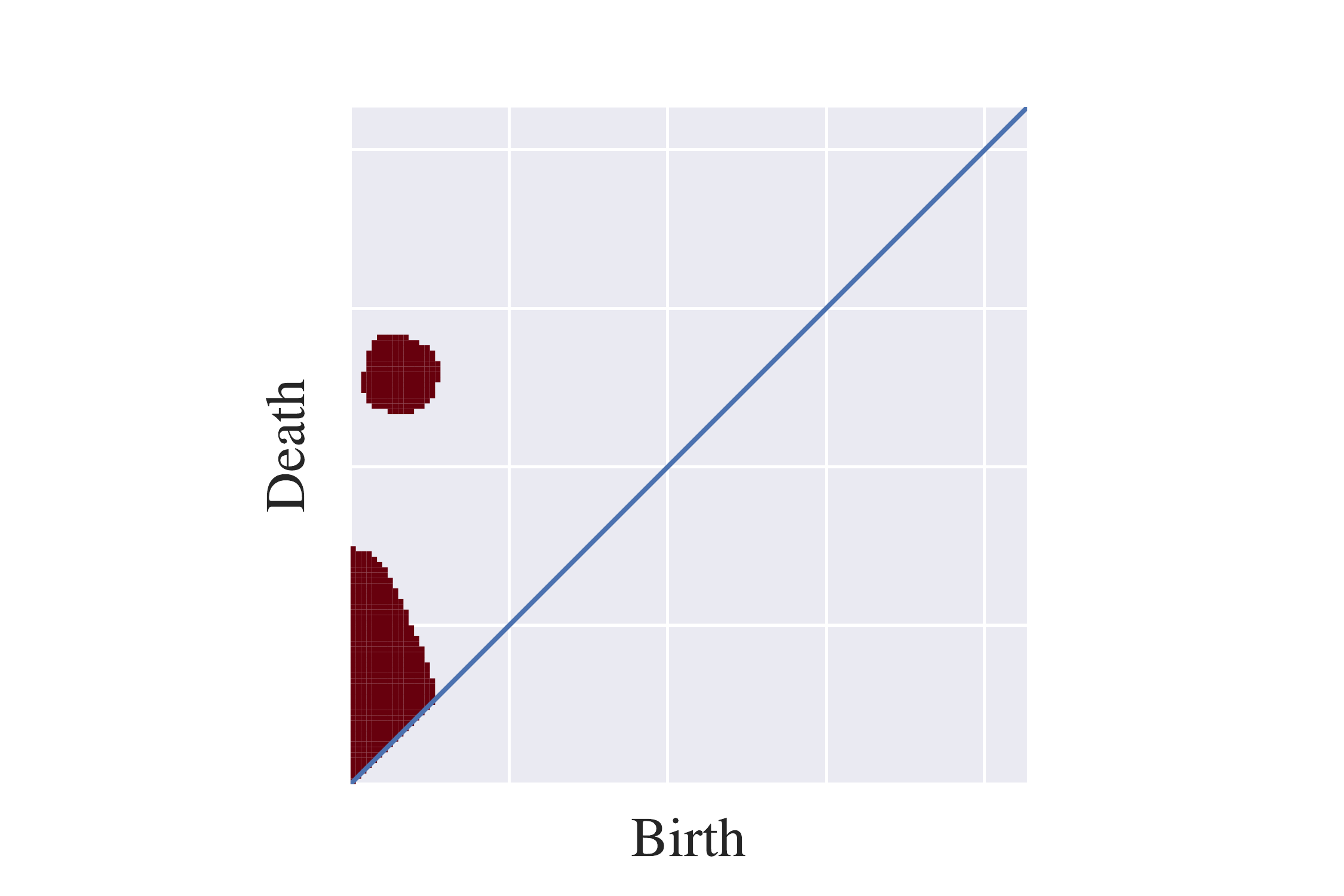}  
  \caption{}
  \label{fig:sub-first}
\end{subfigure}
\begin{subfigure}{.33\textwidth}
  \centering
  \includegraphics[width=1\linewidth]{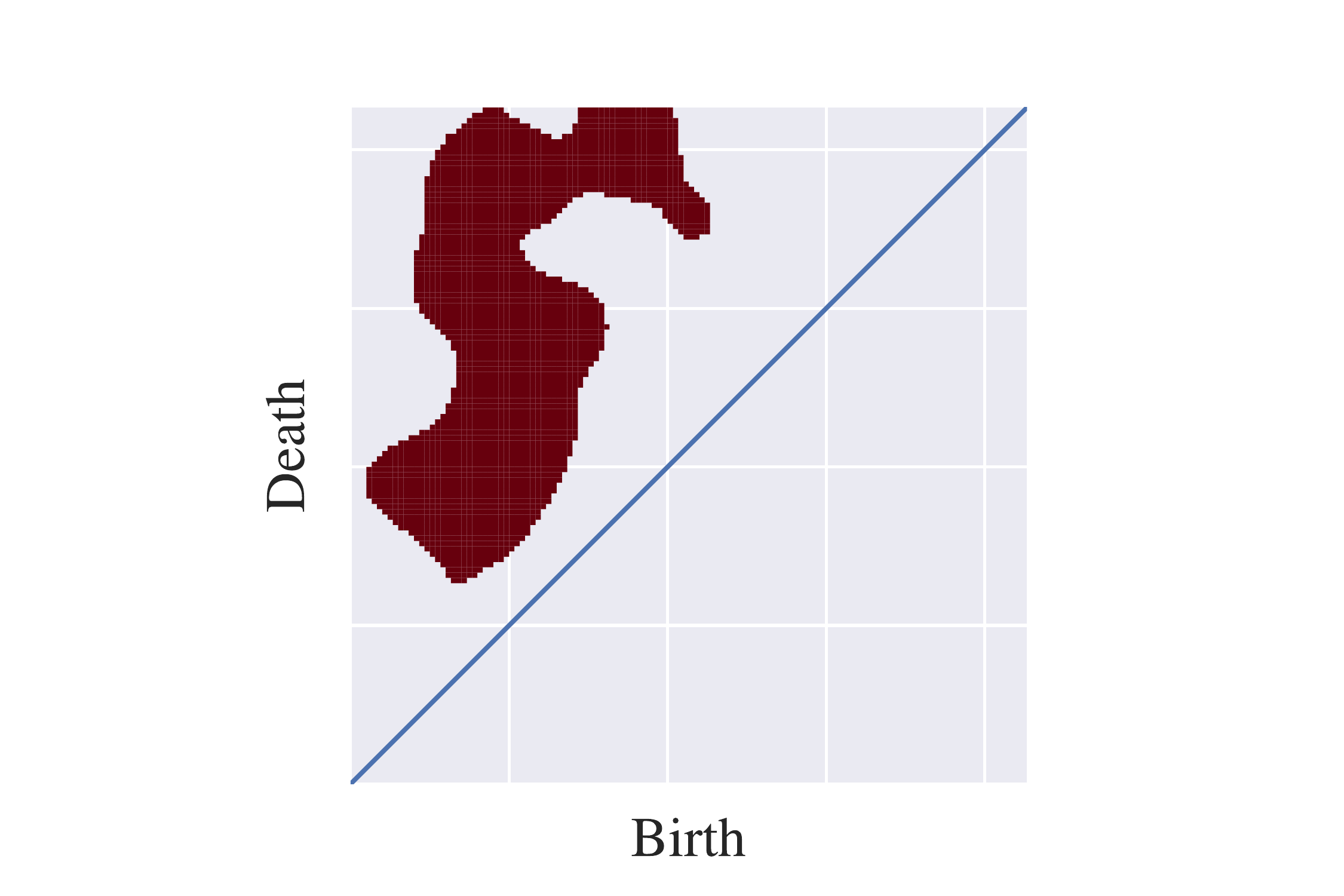}  
  \caption{}
  \label{fig:sub-second}
\end{subfigure}
\caption{Masks highlight the areas of the PD that are important in distinguishing Class 1 from Class 2. We perform inverse analysis on these areas to visualize what features of the original data distinguishing classes. (a) Weights from SVM classification in the space of PDs. The areas of the diagram that distinguish Class 2 are in red and the areas of the diagram that distinguish Class 1 are in blue. (b)  PD mask for Class 1. (c) PD mask for Class 2.}
\label{fig:impweightpcloud}
\end{figure}

An advantage of utilizing a linear classifier is the ability to extract which features are the ones driving classification. Specifically, the magnitude of the weights of the SVM classifier $w \in \mathbb{R}^q$ can be directly associated with the importance of each feature of the vectorized PDs \cite{smith2020convolutional}. Weights with a large negative value are characteristic of point clouds in Class 1 and weights with a large positive value are characteristic of point clouds in Class 2. A visualization of the weights is shown in Figure \ref{fig:pcloudweights}. From this \hl{representation, we} can see that Class 1 is characterized by 1-D holes that are born and die in the early stages of the filtration suggesting a higher number of small radius holes. Class 2 is characterized by 1-D holes that persist over multiple stages of the filtration, suggesting the presence of holes with large radius. These results highlight how one can exploit topological information obtained with TDA to perform statistical (PCA) or machine learning tasks (SVM classification). 

\begin{figure}[!htp]
\begin{subfigure}{.5\textwidth}
  \centering
  \includegraphics[width=0.7\linewidth]{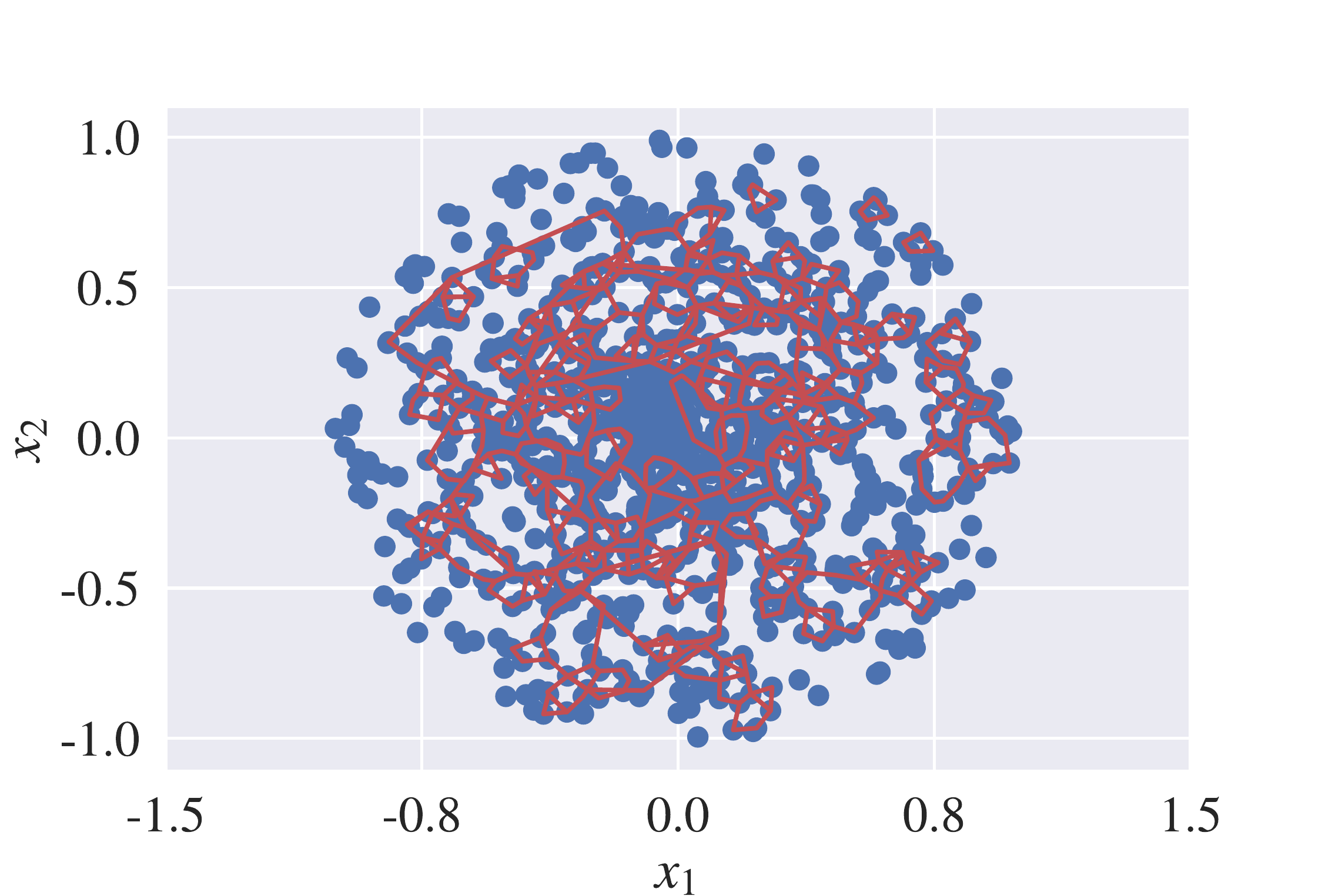}  
  \caption{Inverse $H_1$ analysis for Class 1.}
  \label{fig:sub-first}
\end{subfigure}
\begin{subfigure}{.5\textwidth}
  \centering
  \includegraphics[width=0.7\linewidth]{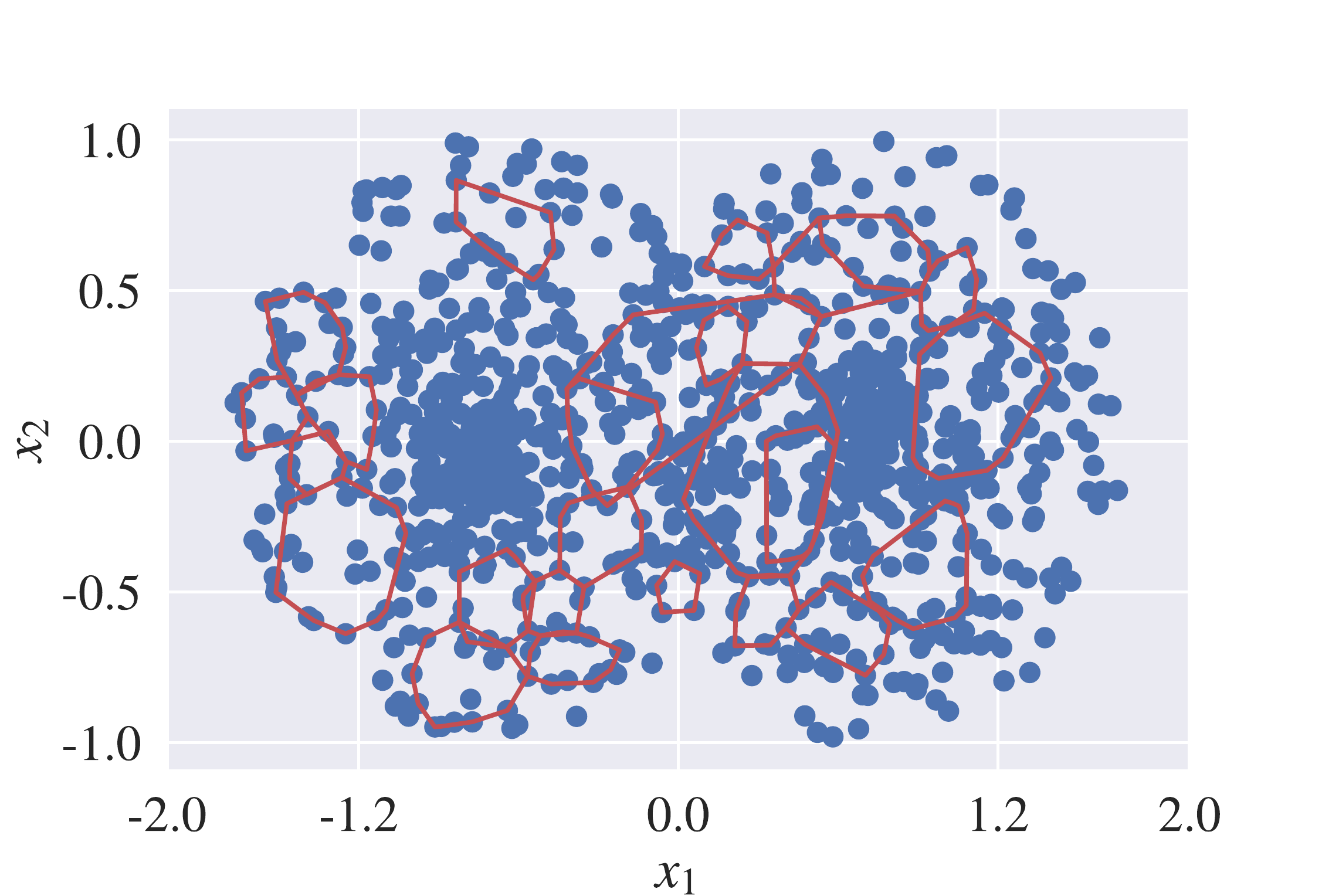}  
  \caption{Inverse $H_1$ analysis for Class 2.}
  \label{fig:sub-second}
\end{subfigure}
\caption{Inverse analysis based on classification weights for Class 1 and Class 2. The analysis reveals that the classifier is separating the classes based on the presence of large cycles in Class 2 and the \hlp{larger number} of smaller cycles in Class 1.  }
\label{fig:invcloud}
\end{figure}

Work conducted in \cite{obayashi2018volume} has led to the development a set of techniques that are useful in the interpretation of PDs. These techniques allow for the identification of \textit{volume-optimal} cycles, which are cycles that correspond to the sparsest representation of the topological features identified in a PD. This technique has been implemented in the \texttt{Homcloud} software, which can be used to identify the 1D holes that are responsible for differences between the classes. Inverse analysis identifies features of the data corresponding to areas of the PD associated with the masks found in Figure \ref{fig:impweightpcloud}. The inverse analysis for a sample of Class 1 and Class 2 is shown in Figure \ref{fig:invcloud}. For Class 2 we see larger separation between points and larger holes that persist for a longer period of the filtration. In Class 1 we see that the separation between points is smaller, resulting in smaller loops that are formed early and die quickly. 

\subsection{Topology of Time-Series and Phase-Planes}
\label{sec:ex2}

\begin{figure}[h!]
\begin{subfigure}{.5\textwidth}
  \centering
  \includegraphics[width=0.7\linewidth]{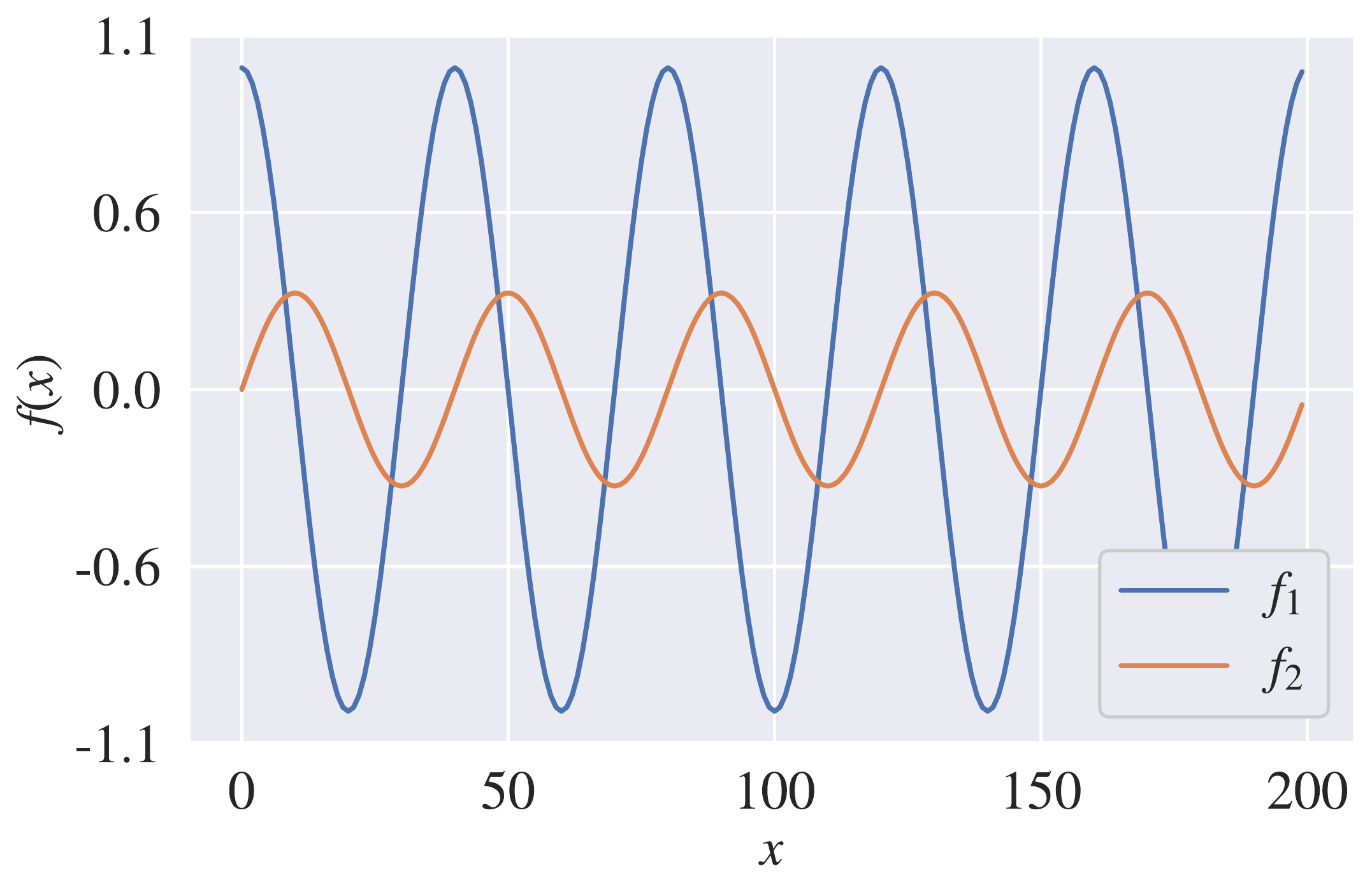}  
  \caption{Time series for functions $f_1$ and $f_2$.}
  \label{fig:sub-first}
\end{subfigure}
\begin{subfigure}{.5\textwidth}
  \centering
  \includegraphics[width=0.7\linewidth]{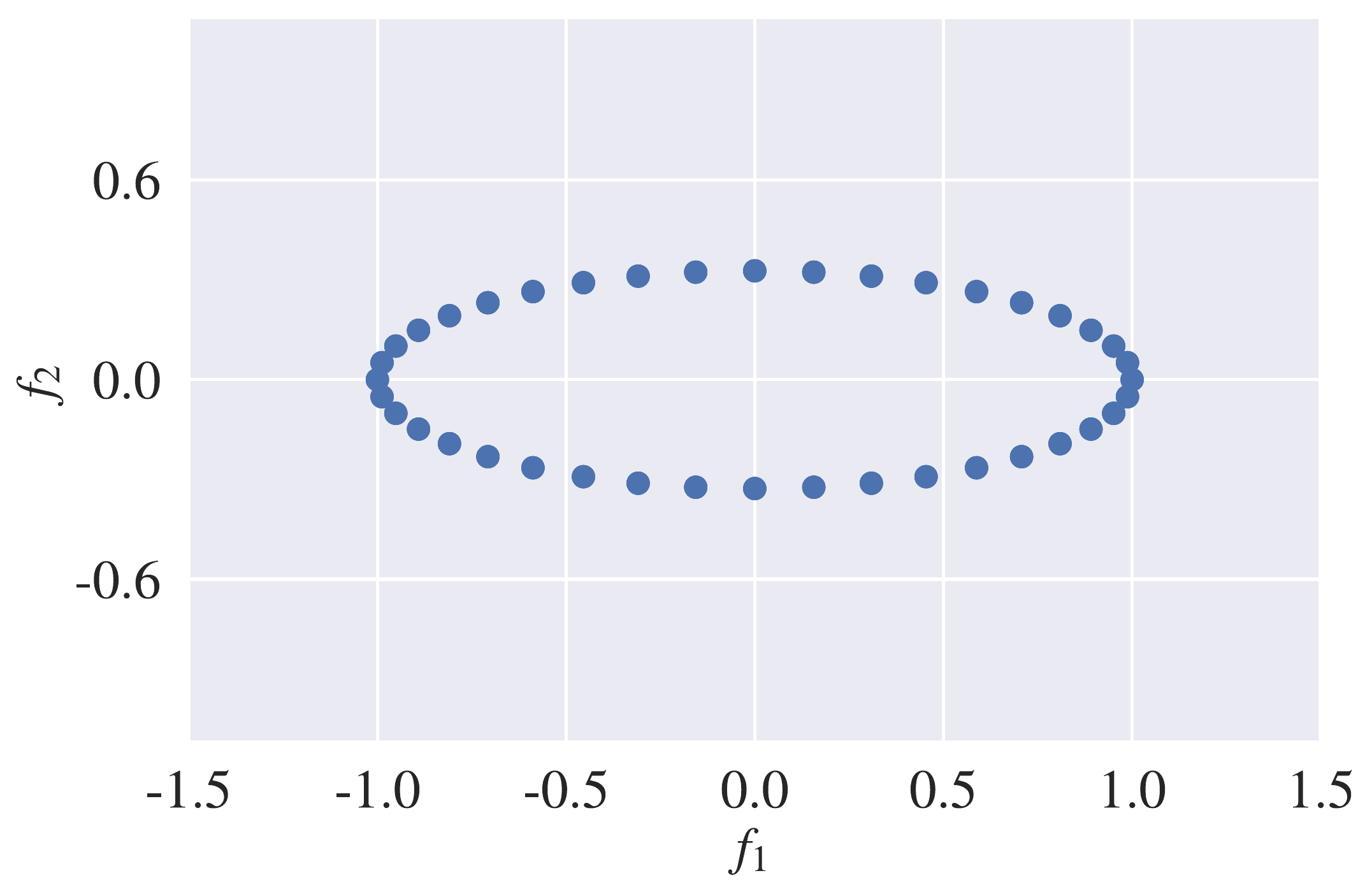}  
  \caption{Phase plane for functions $f_1$ and $f_2$ }
  \label{fig:SS1b}
\end{subfigure}
\caption{Phase plane for periodic orbit with two states. The plane is represented as cloud points from the edge of an ellipse and is ideal for a geometric analysis.}
\label{fig:SS1}
\end{figure}

\begin{figure}[h!]
\begin{subfigure}{.5\textwidth}
  \centering
  \includegraphics[width=0.7\linewidth]{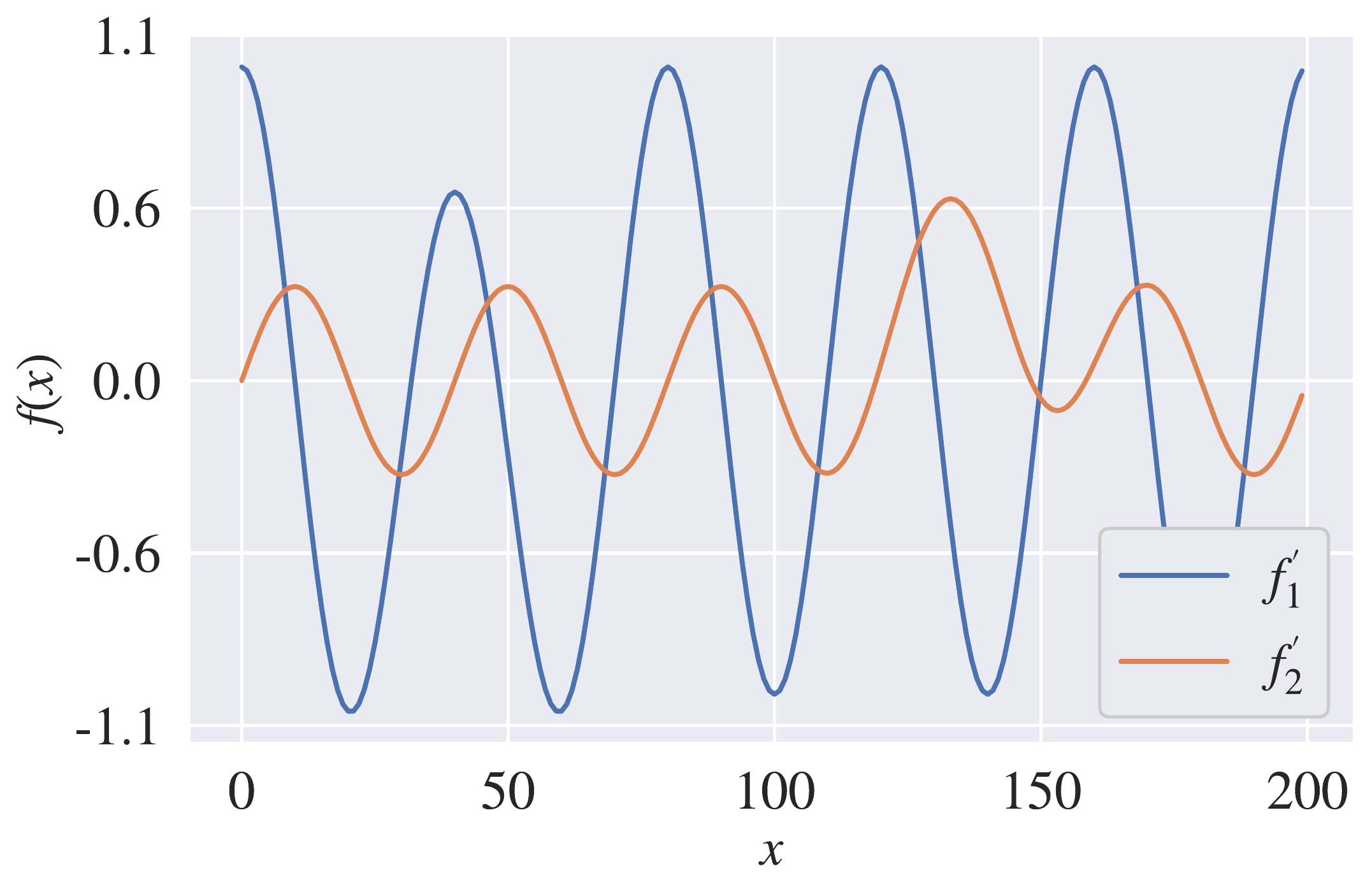}  
  \caption{Timer series for functions $f_1'$ and $f_2'$.}
  \label{fig:sub-first}
\end{subfigure}
\begin{subfigure}{.5\textwidth}
  \centering
  \includegraphics[width=0.7\linewidth]{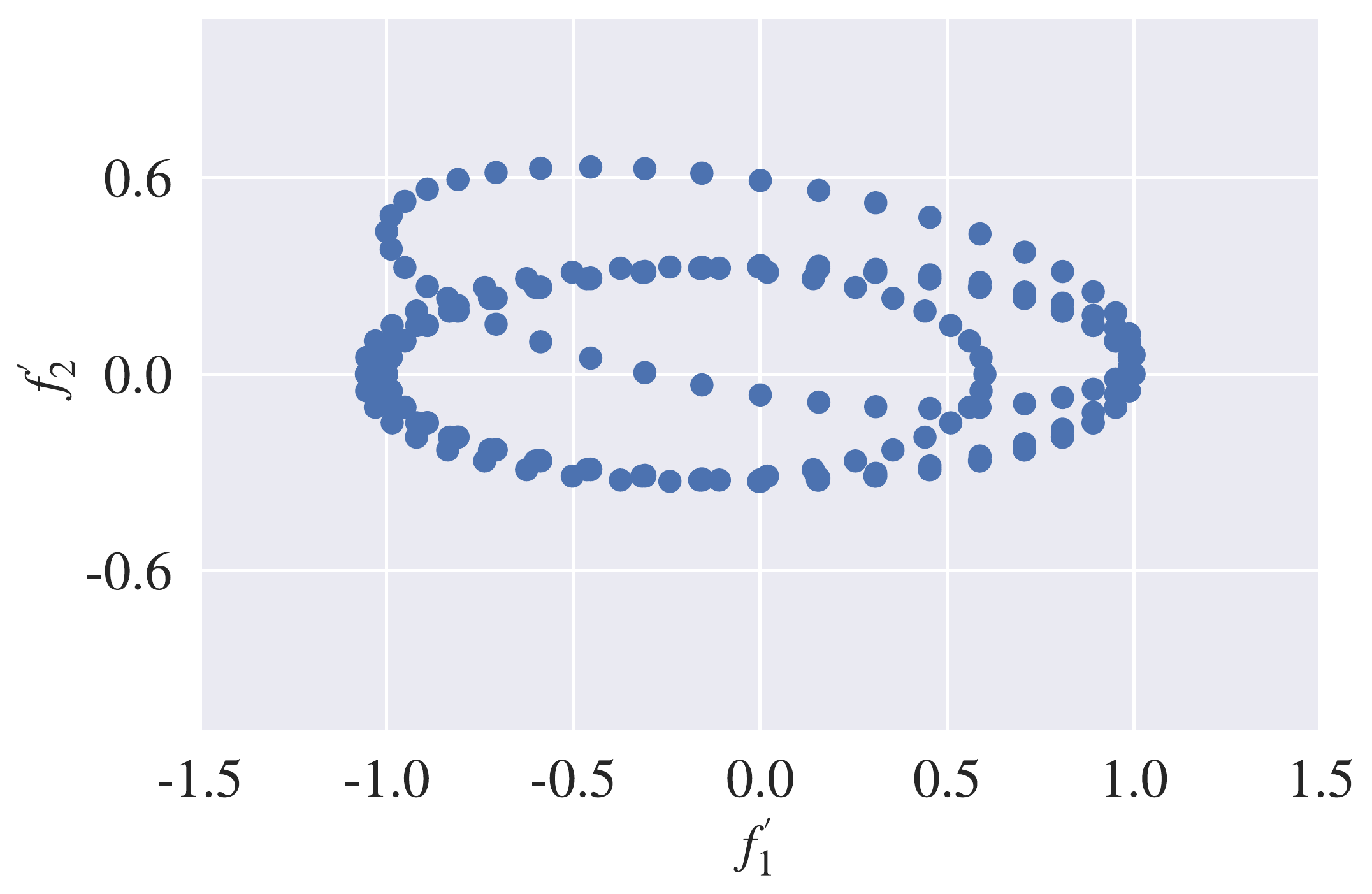}  
  \caption{Phase plane for functions $f_1'$ and $f_2'$}
  \label{fig:sub-second}
\end{subfigure}
\caption{Phase plane for perturbed periodic orbit with two states. The geometry of the plane can no longer be represented as an ellipse, but still represents a sampling from a more complex geometric object.}
\label{fig:SS2}
\end{figure}

Persistence homology has seen applications in the area of time series analysis \cite{seversky2016time,umeda2017time,wang2018topological,perea2019topological,khasawneh2018chatter}. A simple example of the application of persistence homology is the analysis of the topology of phase-planes generated by a dynamical system. An example for two state variables $f_1,f_2$ is shown in Figure \ref{fig:SS1}. The phase plane for functions $f_1$ and $f_2$ is created by plotting the two functions against each other (Figure \ref{fig:SS1b}). The phase plane for this periodic system defines an ellipse, which is easy to characterize (e.g., in terms of its axes). We can add complexity to the topology of the phase plane by perturbing the dynamical system. For example, by adding a perturbation, we change the phase-plane to that shown in Figure \ref{fig:SS2}. The topology of the new plane cannot be fully characterized using simple ellipsoids.  

\begin{figure}[h!]
\begin{subfigure}{.5\textwidth}
  \centering
  \includegraphics[width=0.7\linewidth]{fig2.pdf}  
  \caption{Phase plane of functions $f_1$ and $f_2$ }
  \label{fig:sub-first}
\end{subfigure}
\begin{subfigure}{.5\textwidth}
  \centering
  \includegraphics[width=0.7\linewidth]{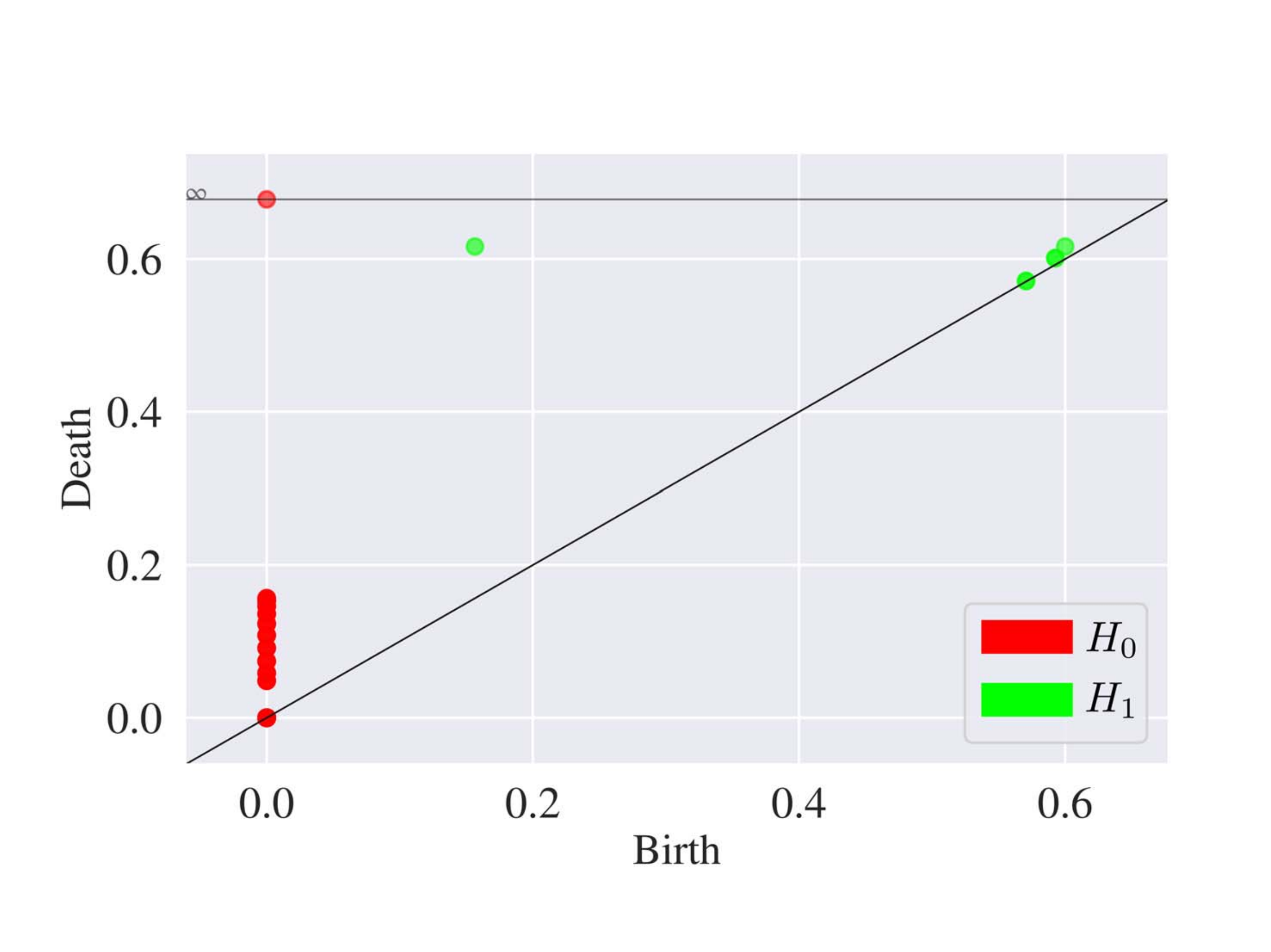}  
  \caption{Persistence Diagram of the phase plane.}
  \label{fig:sub-second}
\end{subfigure}
\caption{PD for the phase plane reflects the presence of a single, persistent loop in the diagram. The persistence diagram captures the important geometric aspects of the data.}
\label{fig:SS1_pd}
\end{figure}

\begin{figure}[h!]
\begin{subfigure}{.5\textwidth}
  \centering
  \includegraphics[width=.7\linewidth]{fig4.pdf}  
  \caption{Phase plane of functions $f_1'$ and $f_2'$ }
  \label{fig:sub-first}
\end{subfigure}
\begin{subfigure}{.5\textwidth}
  \centering
  \includegraphics[width=0.7\linewidth]{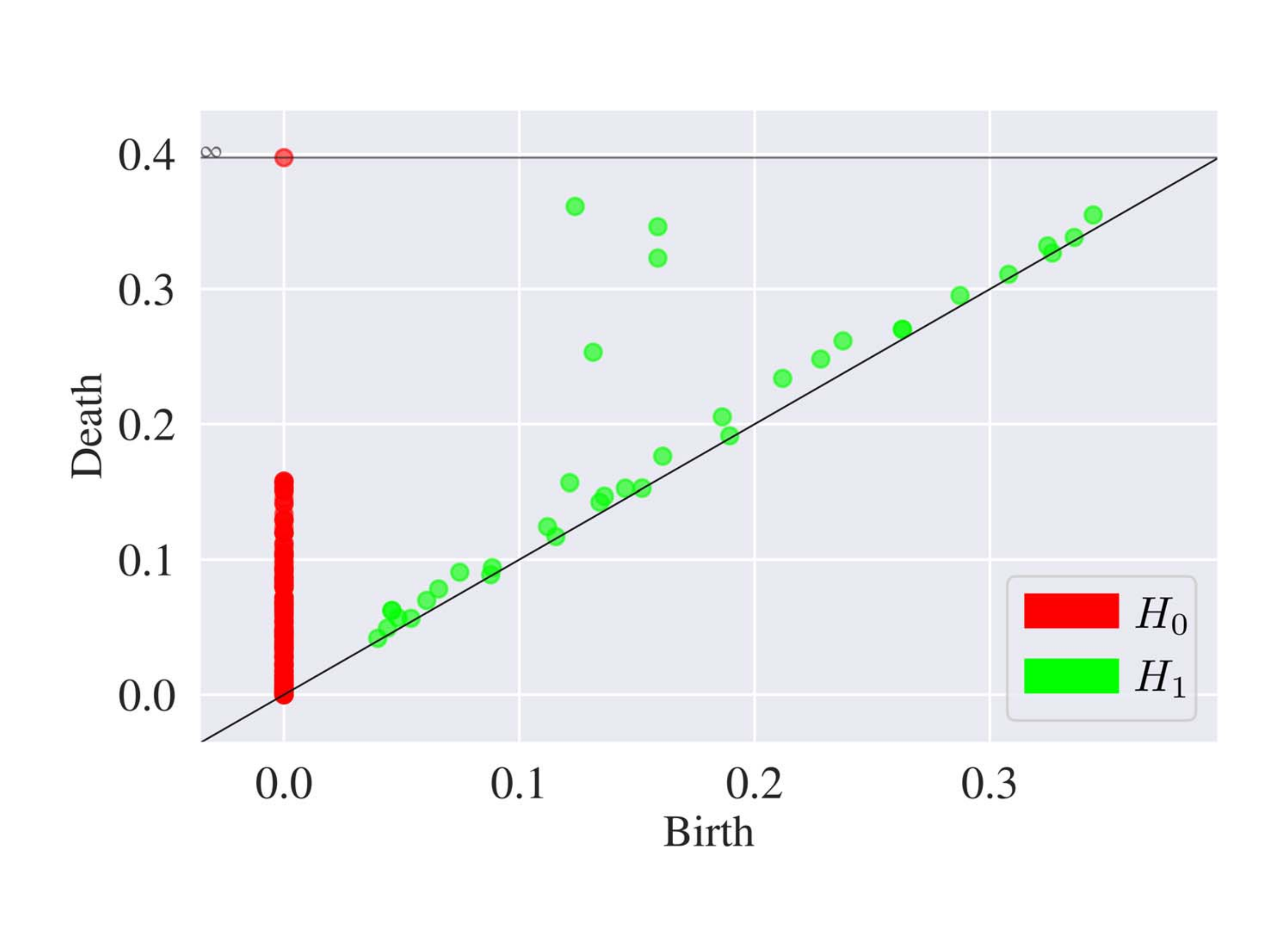}  
  \caption{Persistence diagram of the phase plane.}
  \label{fig:sub-second}
\end{subfigure}
\caption{PD for the phase plane reflects the presence of four loops. The persistence diagram captures the important geometric aspects of the phase plane without having to fit a complex geometric model to the data. }
\label{fig:SS2_pd}
\end{figure}

The analysis of the phase plane topology through an epsilon ball filtration allows us to differentiate the dynamics of the perturbed and unperturbed systems. We compare their PDs in Figure \ref{fig:SS1_pd} and \ref{fig:SS2_pd};  the unperturbed system contains a single highly persistent cycle while the perturbed system contains four cycles that are less persistent. 

\begin{figure}[!htp]
    \centering
    \includegraphics[width = .8\textwidth]{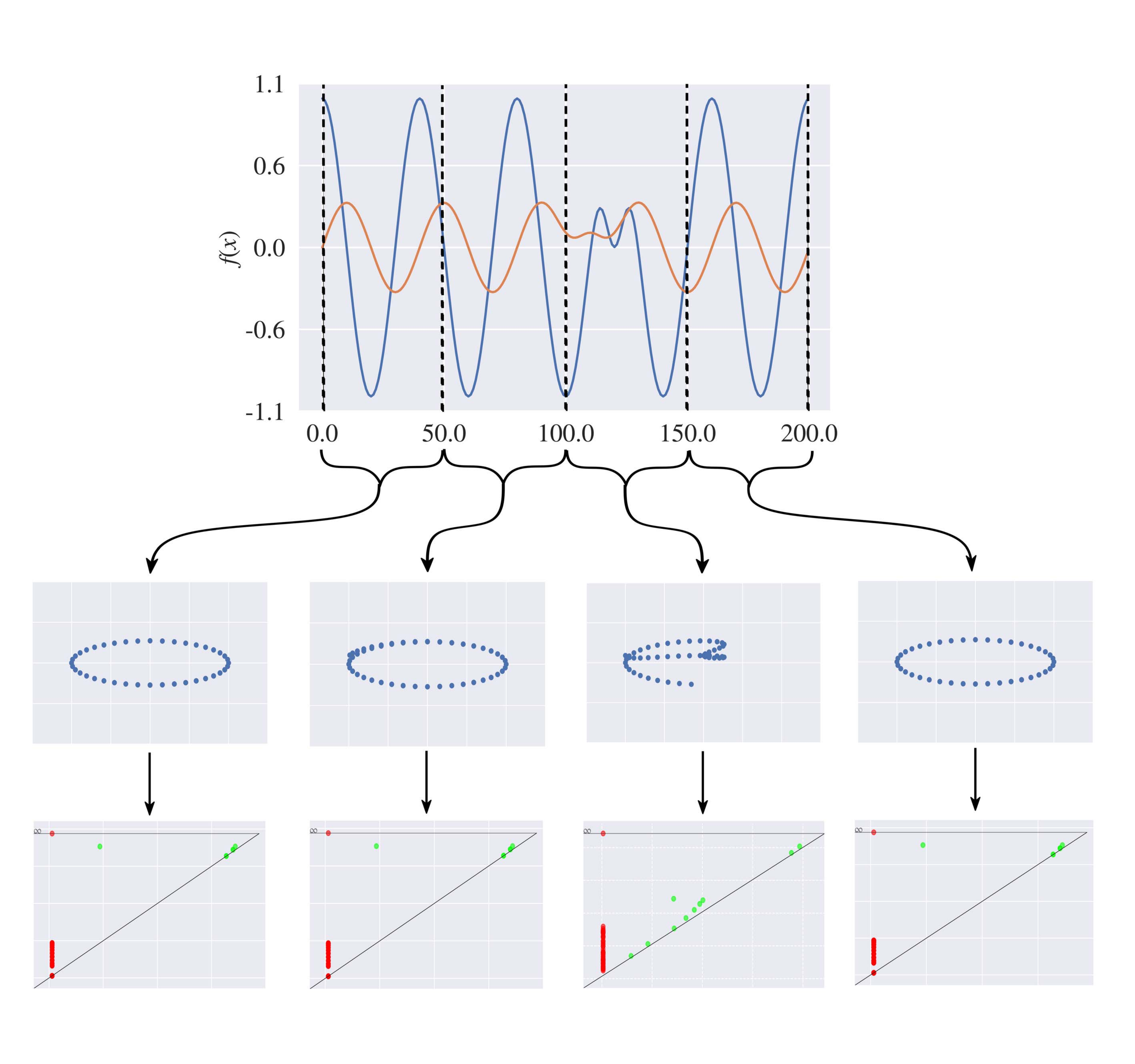}
    \caption{Sliding window method applied to time series. The first, second, and fourth windows have phase planes with similar topologies while the third window has an obvious shift (which introduces a change in the topology of its plane). This change in topology is captured by the persistence diagrams.}
    \label{fig:window1}
\end{figure}

We can also use persistent homology within a {\em sliding windowing method} to detect when a change has occurred in the dynamics of a system (e.g., for fault detection). This concept is illustrated in Figure \ref{fig:window1} where we only show $4$ windows for illustrative purposes; here we can see that the first, second, and fourth windows have similar topology while the third window has a different topology.  We compare the PDs of each window in Figure \ref{fig:window1}; the PDs clearly reveal \hl{that the} phase plane of the third window is different. \hl{ This example demonstrates that the geometry of time series are highly informative and can capture the behavior of the system with minimal information. The small windows used in this example immediately demonstrate the cyclic/periodic nature of the time series and show a large difference in shape when a perturbation occurs. Many current methods require a large amount of information to model the given system via statistics or machine learning methods \cite{malhotra2015long,laptev2015generic,chan2005modeling}.}

\begin{figure}[!htp]
\begin{subfigure}{.5\textwidth}
  \centering
  \includegraphics[width=0.7\linewidth]{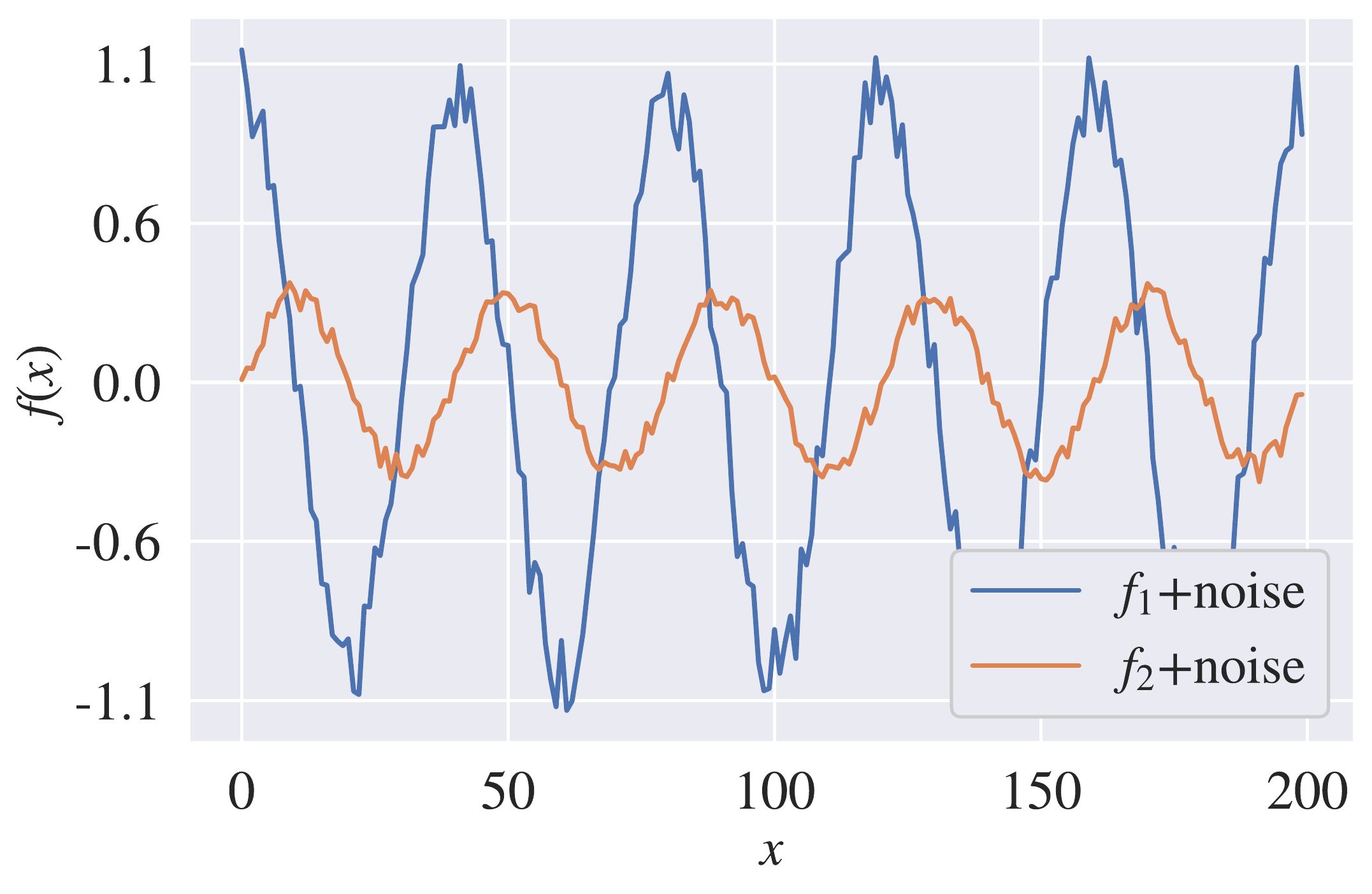}  
  \caption{Time series with added Gaussian noise.}
  \label{fig:sub-first}
\end{subfigure}
\begin{subfigure}{.5\textwidth}
  \centering
  \includegraphics[width=0.7\linewidth]{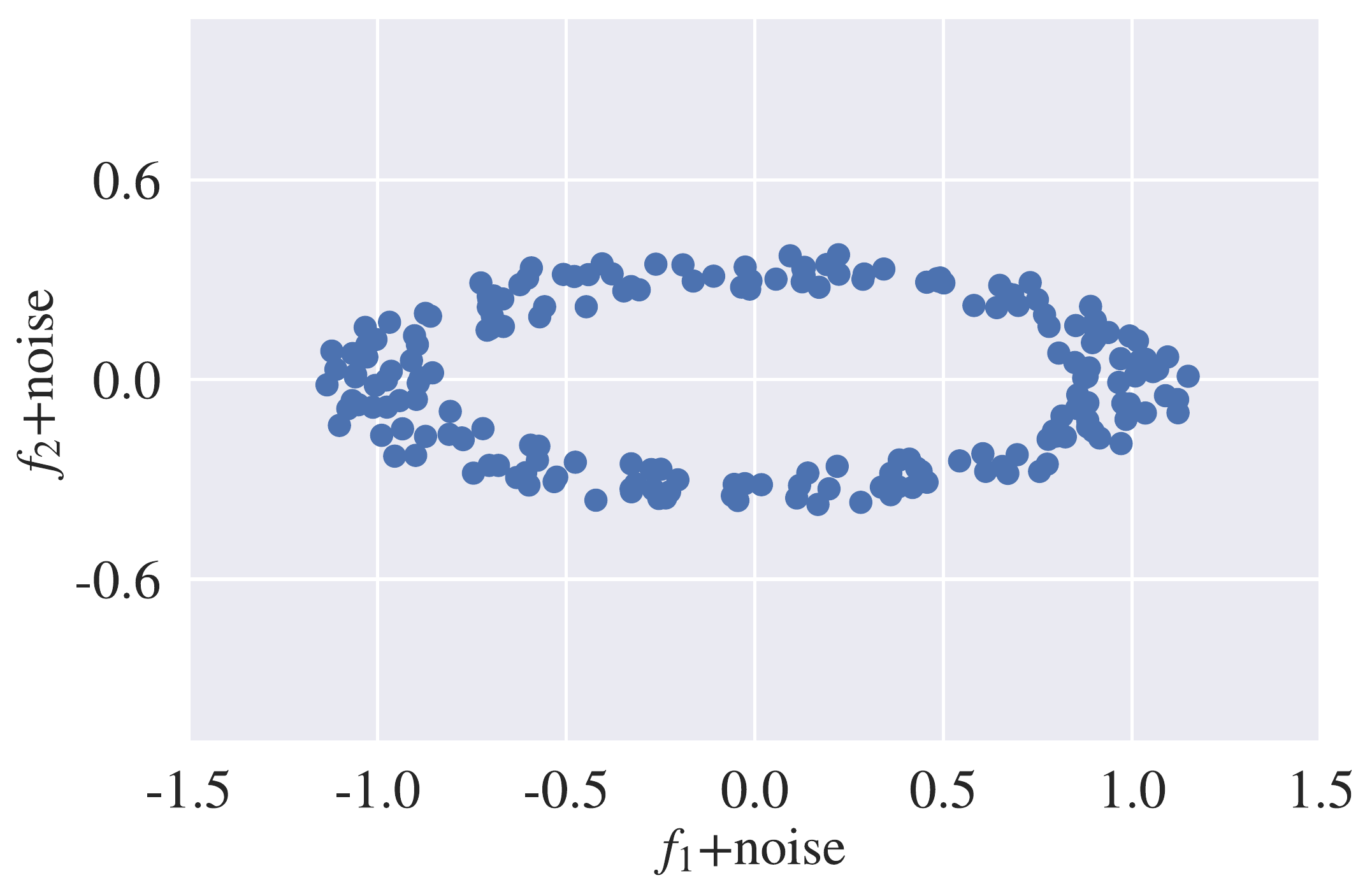}  
  \caption{Phase plane representation.}
  \label{fig:sub-second}
\end{subfigure}
\caption{Phase plane for noisy $f_1$ and $f_2$  shows a similar topology to the noiseless counterparts.}
\label{fig:SSnoise}
\end{figure}

\begin{figure}[!htp]
\begin{subfigure}{.5\textwidth}
  \centering
  \includegraphics[width=0.7\linewidth]{fig7.pdf}  
  \caption{PD for functions $f_1$ and $f_2$}
  \label{fig:sub-first}
\end{subfigure}
\begin{subfigure}{.5\textwidth}
  \centering
  \includegraphics[width=0.7\linewidth]{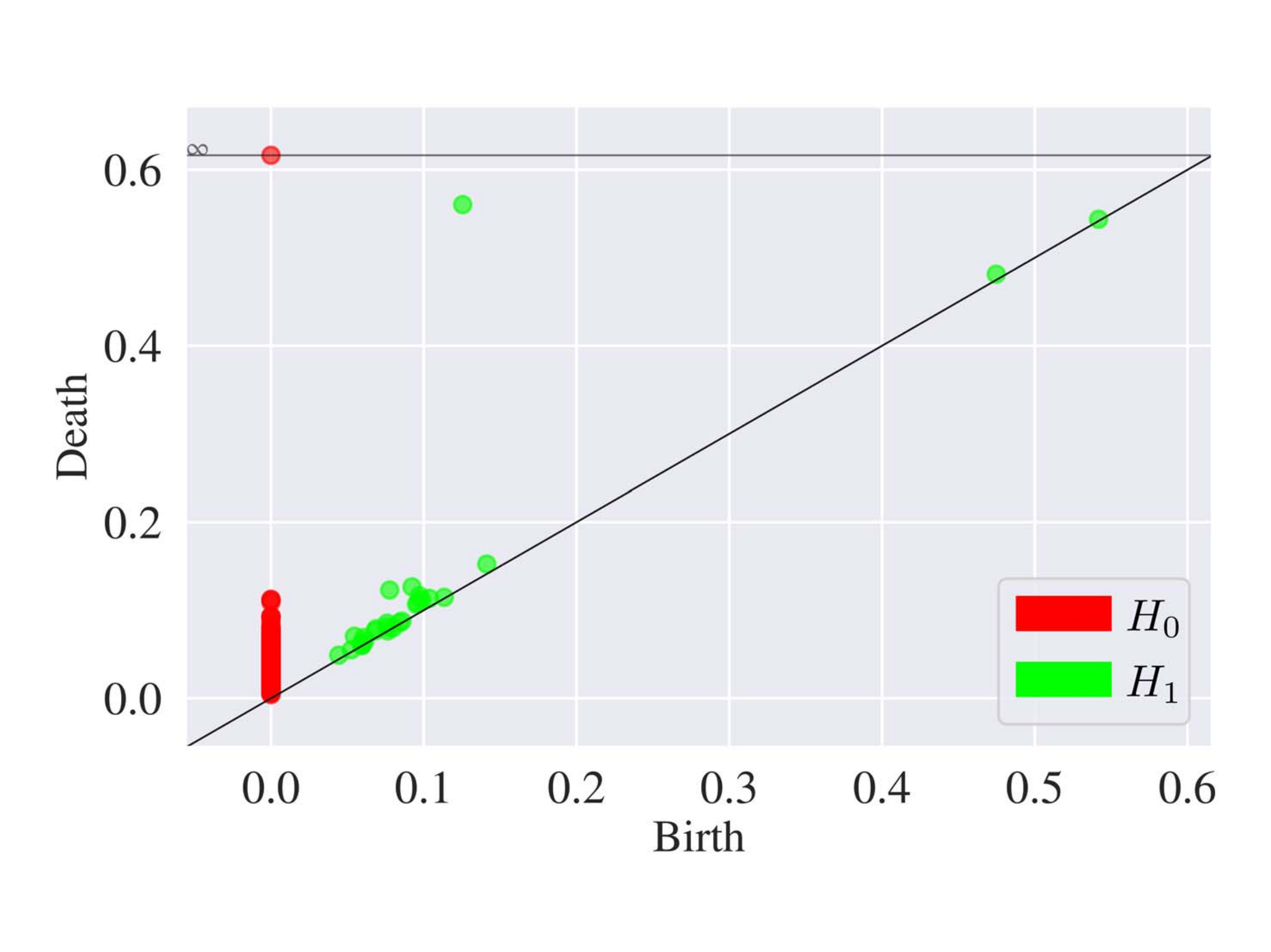}  
  \caption{PD for functions $f_1$ and $f_2$ with added noise}
  \label{fig:sub-second}
\end{subfigure}
\caption{Persistence diagrams for noisy and noiseless functions. Note that the dominant feature (cycle) persists.}
\label{fig:SSnoise_pd}
\end{figure}

Another common perturbation of a system is {\em noise}; an important observation here is that noise usually has the effect of introducing local effects to a trajectory but does not distort the overall topology of the trajectory. This is illustrated in Figure \ref{fig:SSnoise}, where we can see that the ellipse shape of the phase plane is retained. In other words, the topology of the phase plane is stable in the presence of noisy perturbations. In a real system, one may wish to characterize whether noise-free and noisy systems exhibit similar dynamics. One way to do this is to perform persistent homology calculations on both phase planes and to compare the resulting PDs. We show the results of this analysis in Figure \ref{fig:SSnoise_pd};  the PDs reveal that both phase planes contain a persistent cycle (which indicates that they have phase planes with similar topologies). Features created by noise are shown at the bottom of  the diagonal (these features have short persistence). These results demonstrate the stability of PDs  \cite{chazal2009proximity}, which is an important concept in TDA. These results also highlight that TDA is a powerful tool for the classification of time series \cite{umeda2017time}, identification of periodic orbits and shifts \cite{perea2015sliding}, and for change point detection \cite{gidea2018topological}. The paper of Perea \cite{perea2019topological} provides an excellent overview of TDA in signal processing.

\subsection{Topology of 2D Scalar Fields}
\label{sec:ex3}

\begin{figure}[h!]
     \centering
     \includegraphics[width = .9\textwidth]{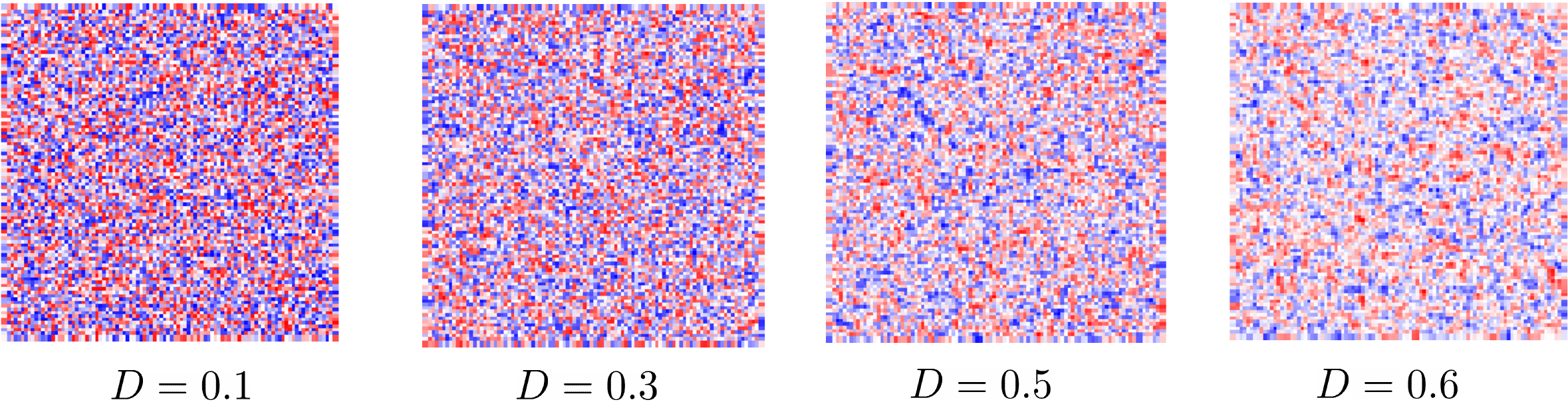}
     \caption{Samples from 2D scalar fields with their corresponding diffusion coefficient ($D$) value, where red represents high intensity values and blue represents small intensity values.}
     \label{fig:2ddiff}
 \end{figure}

\begin{figure}[h!]
     \centering
     \includegraphics[width = .5\textwidth]{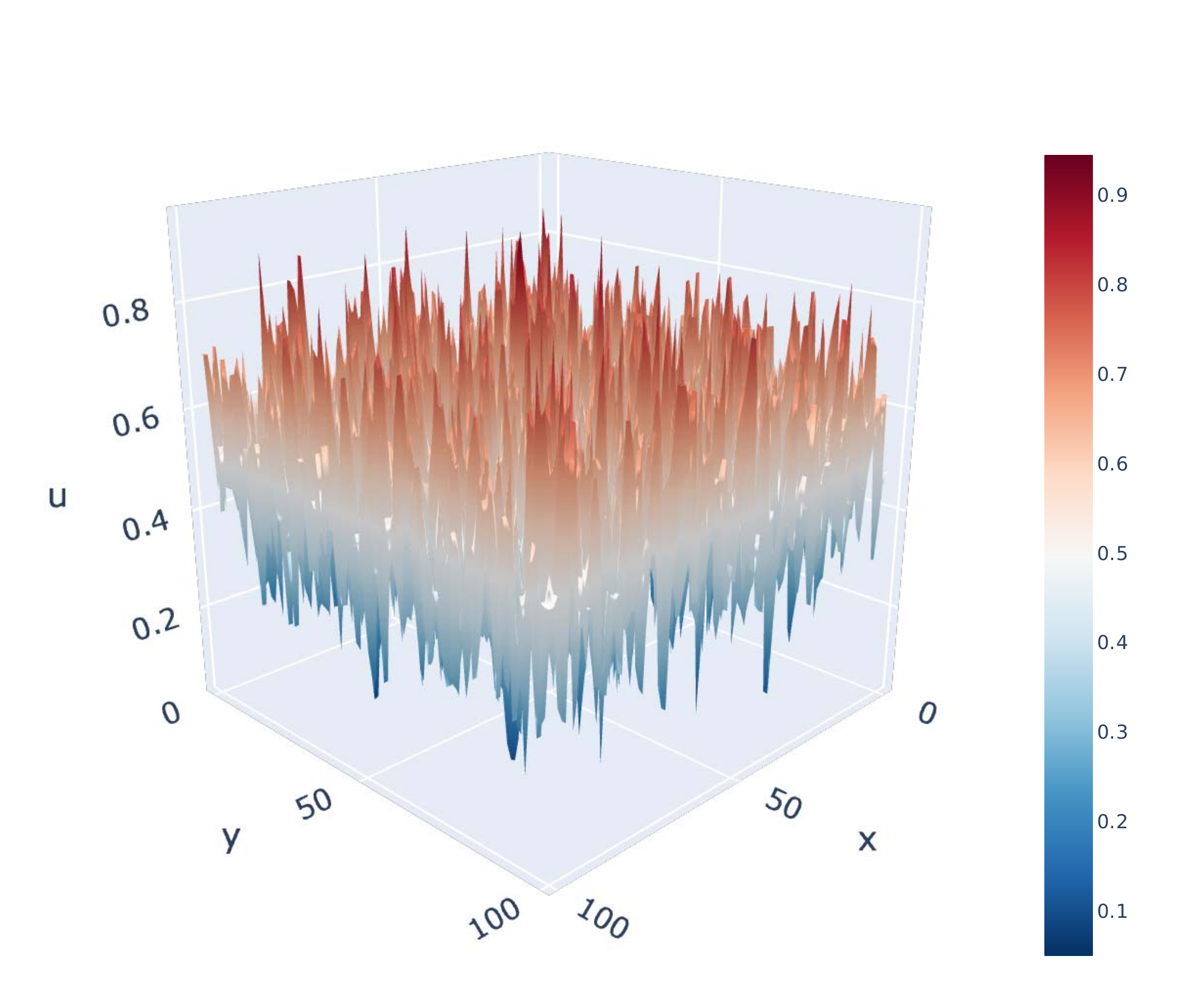}
     \caption{3D-functional representation of a scatter field with diffusion coefficient $D = 0.6$. The function is treated as a cubical complex and the filtration is performed over the scalar value.}
     \label{fig:diffsurf}
 \end{figure}

In this example we \hl{use TDA} to analyze the topology of 2D scalar fields over a discrete $n$ by $n$ domain $U$. We generate data by applying propagating a random field $\{u_{i,j} : u_{i,j} \in \mathcal{U}(0,1)\}$ using the dynamic 2D diffusion equation \eqref{eq:2ddiff} and we obtain the final steady-state field. We generate fields with different textures by using different diffusion coefficients $D$ and \hl{independent, random intializations}. The resulting fields are used as the datasets; illustrative examples for different diffusion coefficients are shown in Figure \ref{fig:2ddiff}. Here, the blue color are points of small intensity (small values of $u_{i,j}$) while the red color are points of high intensity. We see that small coefficients generate textures that are more granular while large coefficients generate smoother textures. Our goal is to characterize the geometry of the fields to investigate if their underlying structure can be correlated to the diffusion coefficient. In our analysis, we represent the scalar field as a function (in 3D), as shown in Figure \ref{fig:diffsurf}. This functional representation of the field reveals that low and high intensity points define critical points of the function. This also reveals that the field has complex topological features (many critical points with no obvious patterns are present). 

\begin{equation}
    \frac{u_{i,j}^{(n+1)} - u_{i,j}^{(n)}}{\Delta t} = D\left[ \frac{u_{i+1,j}^{(n)} - 2u_{i,j}^{(n)} + u_{i-1,j}^{(n)}}{(\Delta x)^2} + \frac{u_{i,j+1}^{(n)} - 2u_{i,j}^{(n)} + u_{i,j-1}^{(n)}}{(\Delta y)^2} \right]
\label{eq:2ddiff}
\end{equation}

\begin{figure}[!htp]
\begin{subfigure}{.5\textwidth}
  \centering
  \includegraphics[width=.9\linewidth]{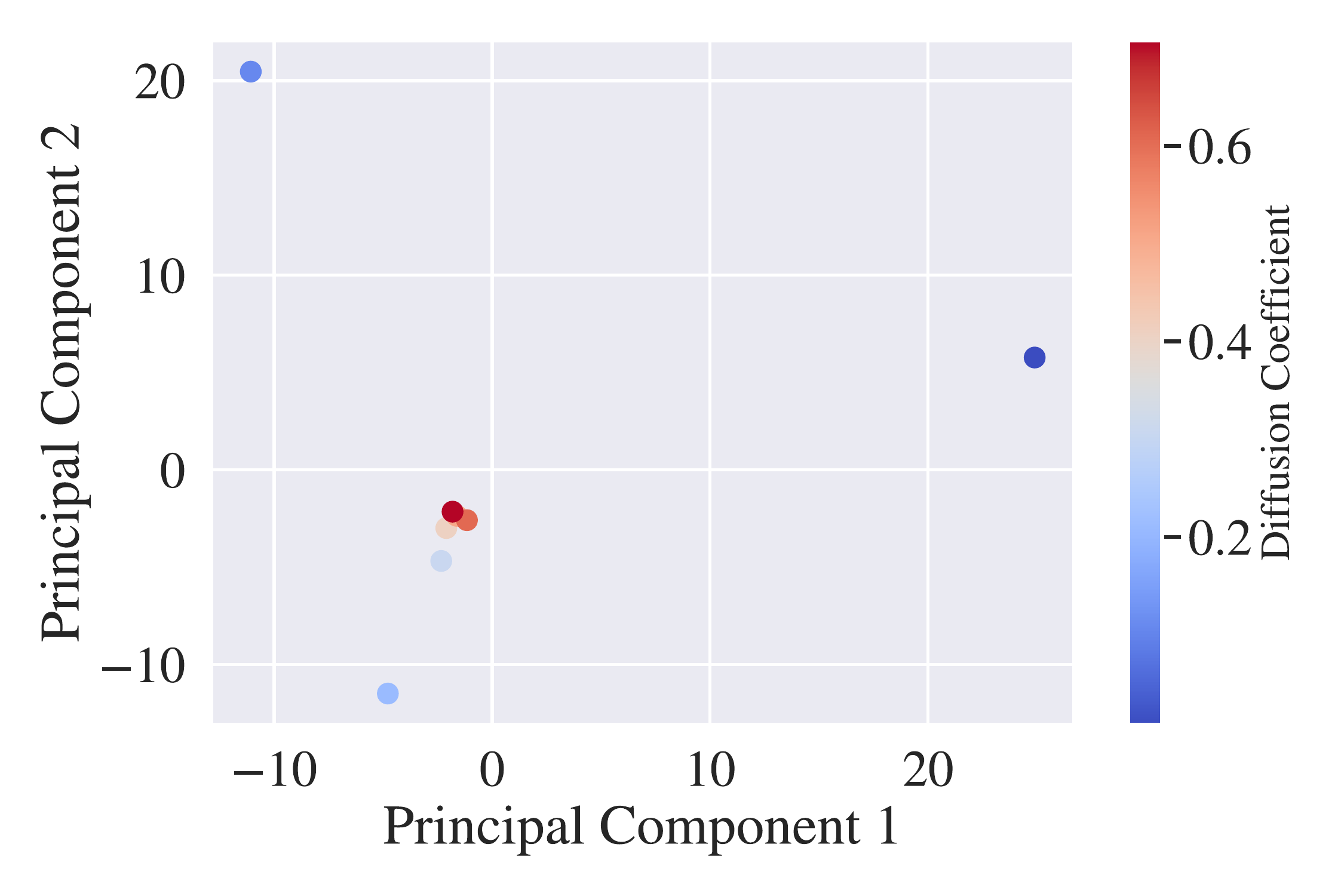}  
  \caption{Dominant principal components from PCA.}
  \label{fig:pcadiff}
\end{subfigure}
\begin{subfigure}{.5\textwidth}
  \centering
  \includegraphics[width=.9\linewidth]{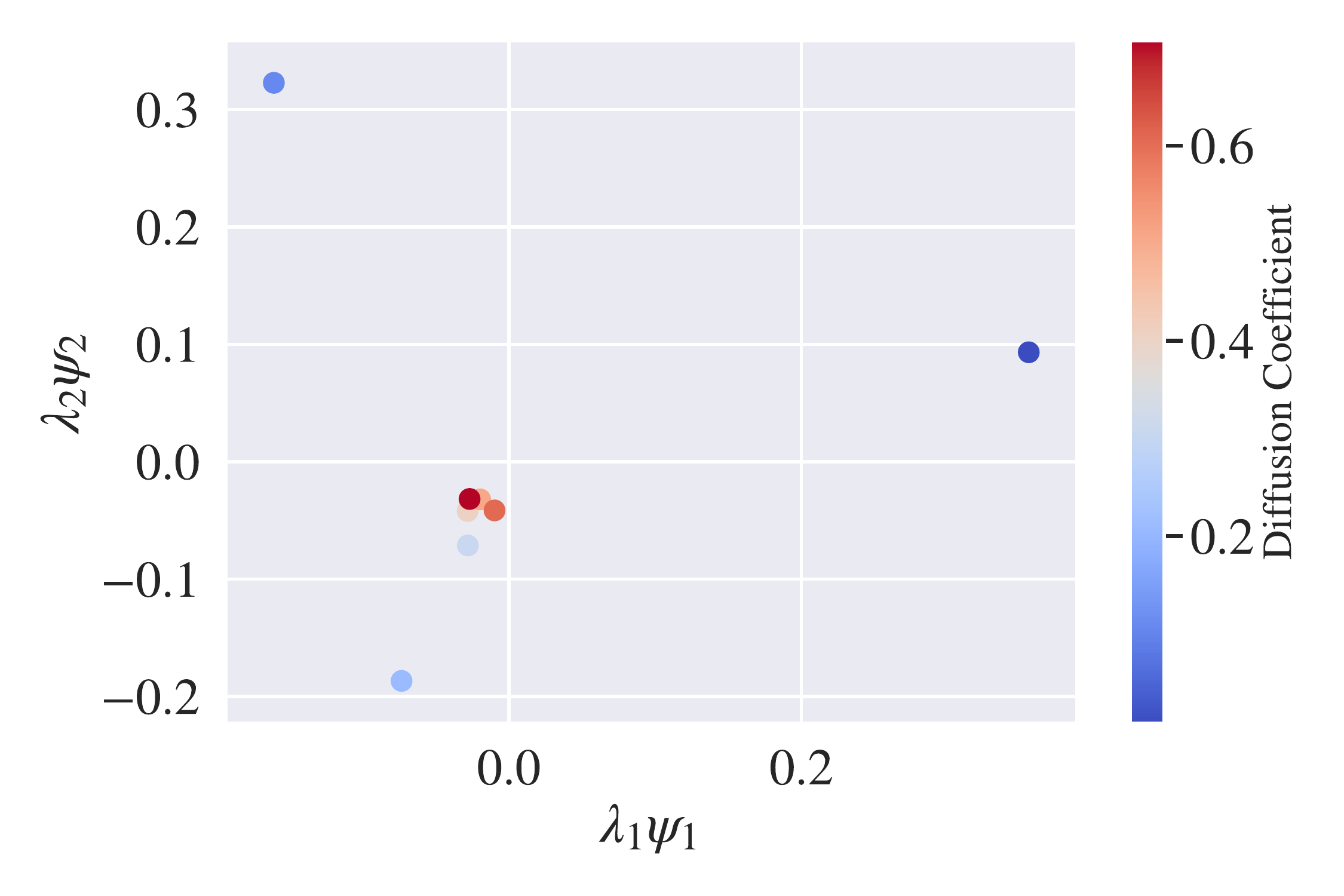}  
  \caption{Dominant features from diffusion map.}
  \label{fig:dmapdiff}
\end{subfigure}
\caption{Dimensionality reduction for the 2D dimensional scalar fields using PCA and diffusion maps.}
\label{fig:dimred}
\end{figure}

To illustrate the benefits of using TDA against other techniques, we investigate whether the structure of the dataset can be revealed by direct application of PCA  \cite{jolliffe2016principal} and diffusion maps  \cite{coifman2006diffusion}. This is a simple, naive comparison but will demonstrate that the datasets have  nontrivial structure. The projection of the data onto the first two principal components is shown in Figure \ref{fig:pcadiff}. Here, we highlight the points based on the associated diffusion coefficient. We also apply the diffusion maps, which is a nonlinear dimensionality reduction technique, and we obtain similar results (see Figure \ref{fig:dmapdiff}).  From these results we see that the features extracted by PCA and diffusion maps do not correlate to the diffusion coefficient.
 
\begin{figure}[!htp]
    \centering
    \includegraphics[width = 1\textwidth]{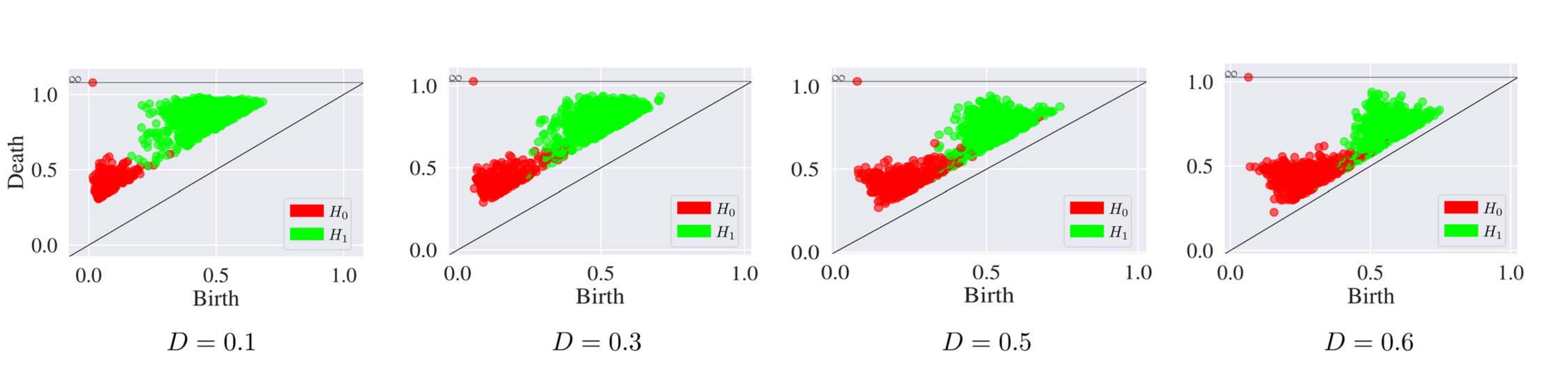}
    \caption{Evolution of PDs with the diffusion coefficient. A dependence of the topology with the diffusion coefficient emerges. }
    \label{fig:2ddiff_pd}
\end{figure}

\begin{figure}[!htp]
    \centering
    \includegraphics[width = .5\textwidth]{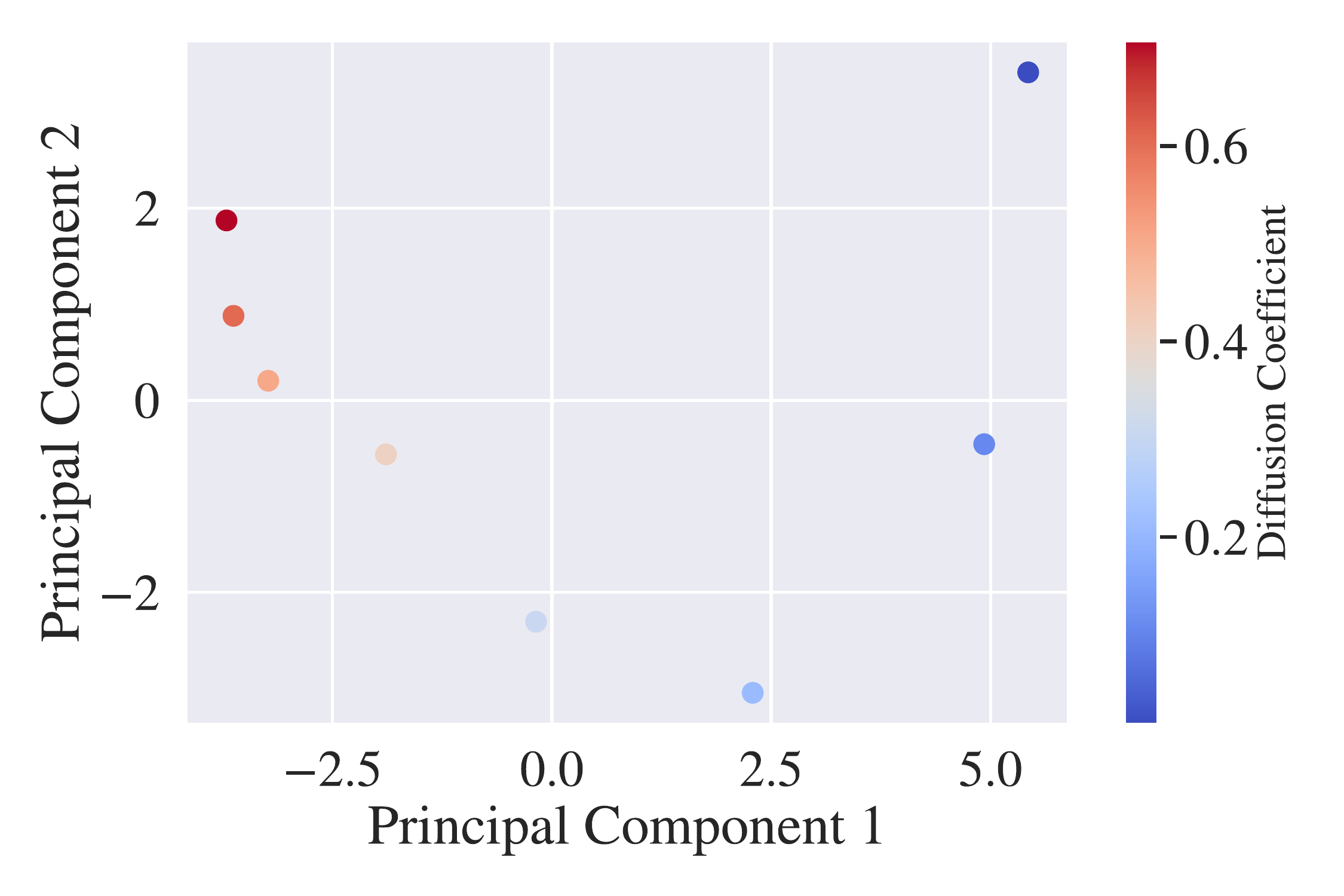}
    \caption{PCA performed on the PD for two-dimensional scalar fields. This reveals that the geometry of the dataset is directly related to the diffusion coefficient.}
    \label{fig:pca_tda}
\end{figure}

We now apply TDA to the 3D field functions and extract persistence diagrams. Example persistence diagrams for $H_0$ and $H_1$ are shown in  Figure \ref{fig:2ddiff_pd}. We see that there is a visual shift in the PDs as we increase the diffusion coefficient. This seems to indicate that the PDs vary continuously with the diffusion coefficient. To verify this, we vectorize the PDs and apply PCA to the vectorized diagrams. The projection of the vectorized PDs onto the dominant principal components is shown in Figure \ref{fig:pca_tda}. It is clear that there is a continuous dependence of the PD on the diffusion coefficient (it forms a continuous manifold). This result provides another demonstration of the stability of persistence diagrams and on how topology varies continuously under perturbations. Specifically, stability indicates that small changes in a given function ($f,g$) results in bounded  changes in the associated persistence diagrams ($PD_f,PD_g$). Thus, because our perturbations to each function are based on changes in $D$, we can guarantee that the distance between PDs is bounded by the size in the perturbation in the diffusion coefficient (i.e., the distance is not arbitrary).

\subsection{Topology of Images}
\label{sec:ex4}

\begin{figure}[h!]
    \centering
    \includegraphics[width = .4\textwidth]{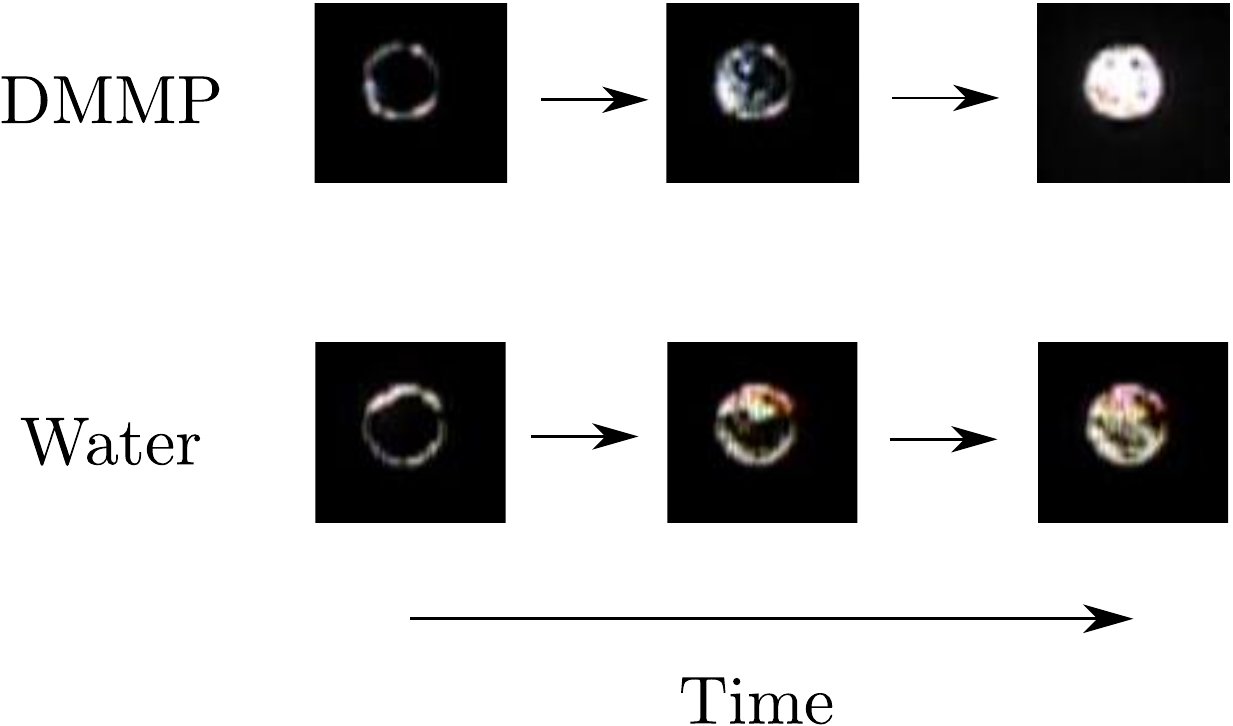}
    \caption{Optical patterns for a liquid crystal sensor when exposed to DMMP or water.}
    \label{fig:LC_resp}
\end{figure}

We now illustrate how to use TDA to analyze images. Specifically, we analyze the optical response of liquid crystal sensors to air contaminants, in particular dimethyl methylphosphonate (DMMP) \cite{smith2020convolutional,shah2001principles}. We analyze the \hlp{spatial} response of a sensor in the presence of DMMP and in the presence of humid nitrogen (nitrogen + water), which we will refer to as water. The sensor responds to both DMMP and water, but there are subtle spatial differences in the optical response of the sensor to these different environments. The working principle of these sensors relies on a change in orientation of liquid crystal molecules in a film when exposed to an air contaminant (analyte). The change in orientation results in optical fields with different spatial and color features (see Figure \ref{fig:LC_resp}). We use TDA to investigate if the topological features of these patterns present a dependence on the air environment. Such information can be used to design sensors (i.e., we can calibrate the sensor by correlating the optical response to the presence of DMMP). 

%

 \begin{figure}[!htp]
  \begin{subfigure}{.33\textwidth}
   \centering
   \includegraphics[width=1\linewidth]{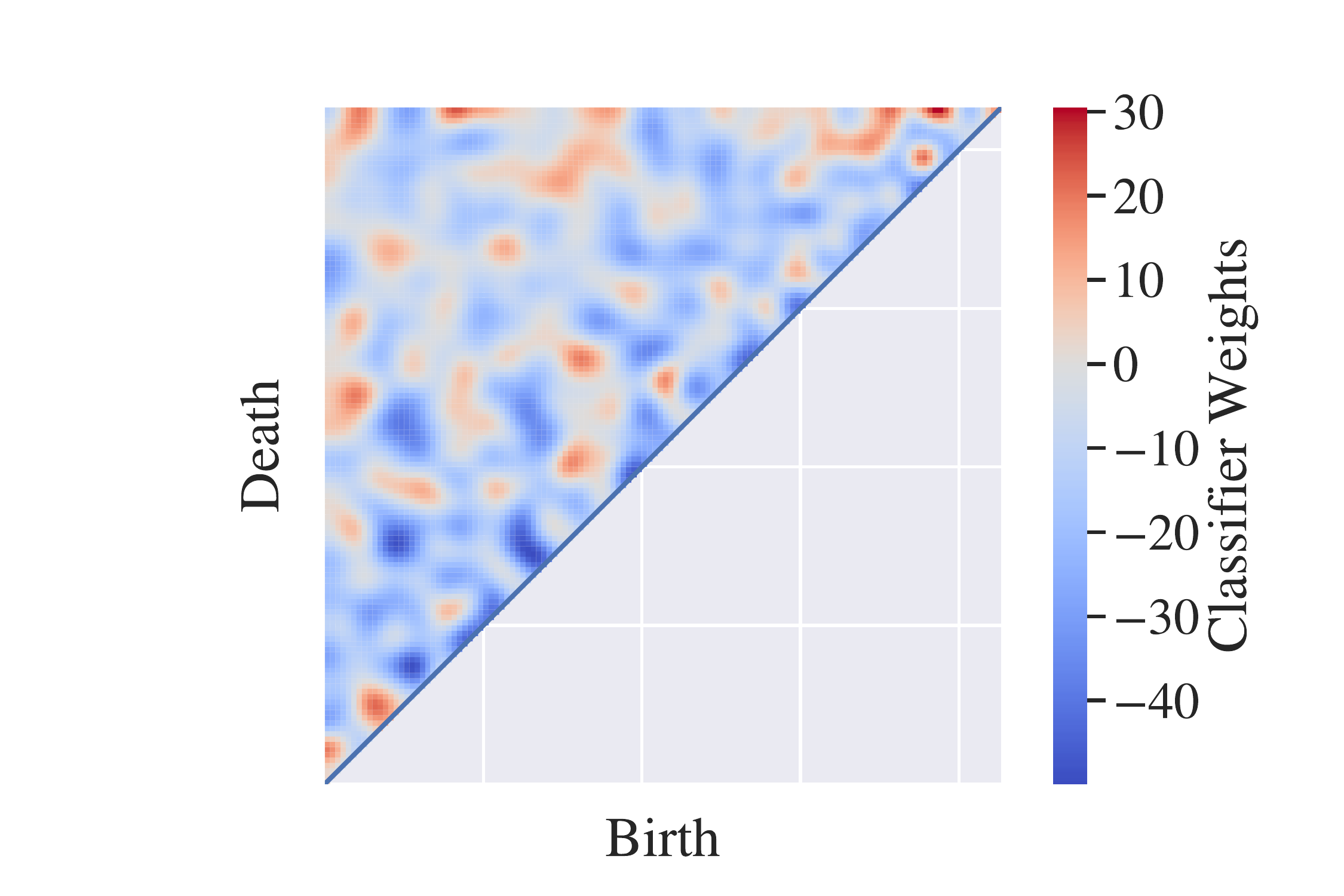}  
   \caption{}
 \end{subfigure}
 \begin{subfigure}{.33\textwidth}
   \centering
   \includegraphics[width=1\linewidth]{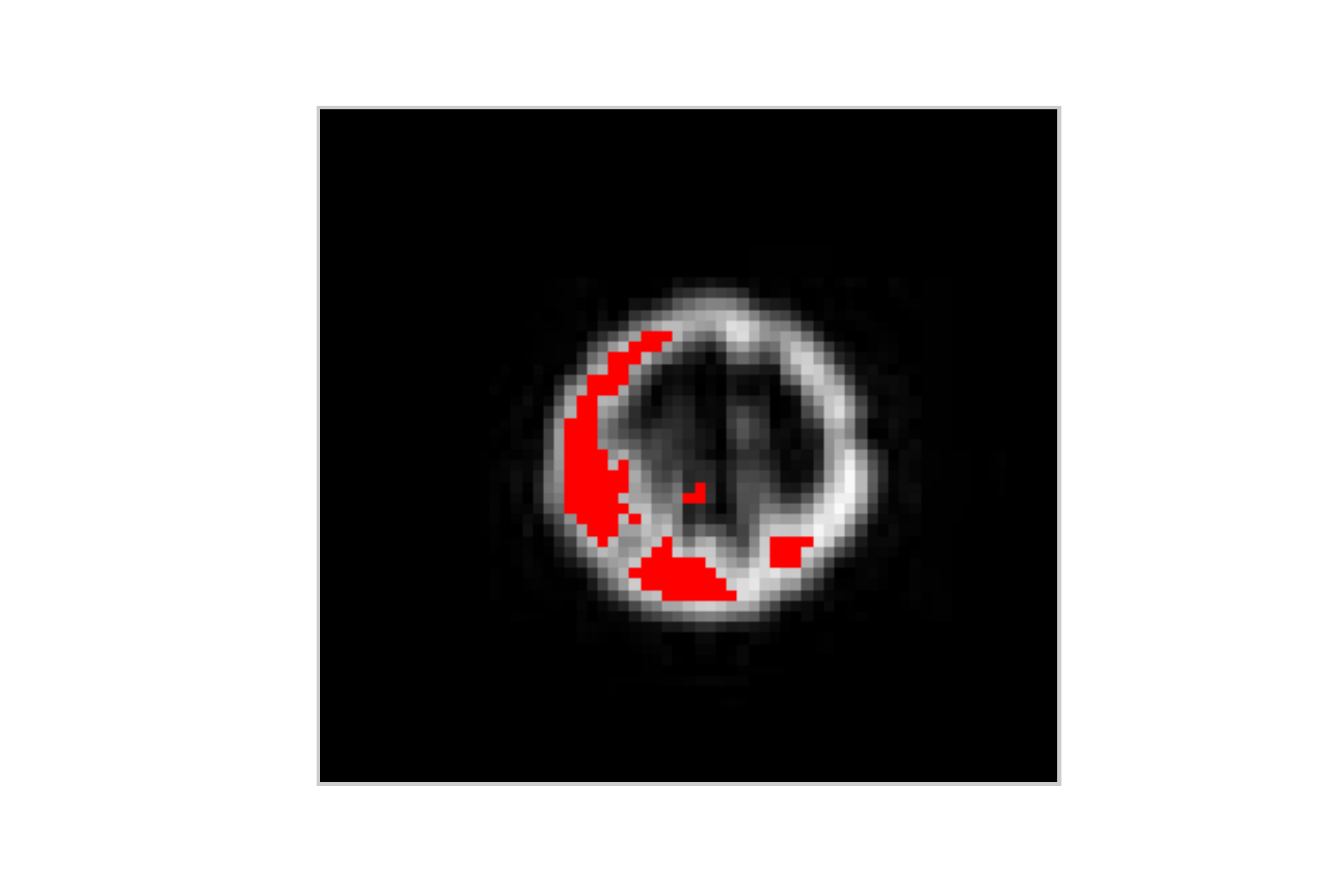}  
   \caption{}
 \end{subfigure}
 \begin{subfigure}{.33\textwidth}
   \centering
   \includegraphics[width=1\linewidth]{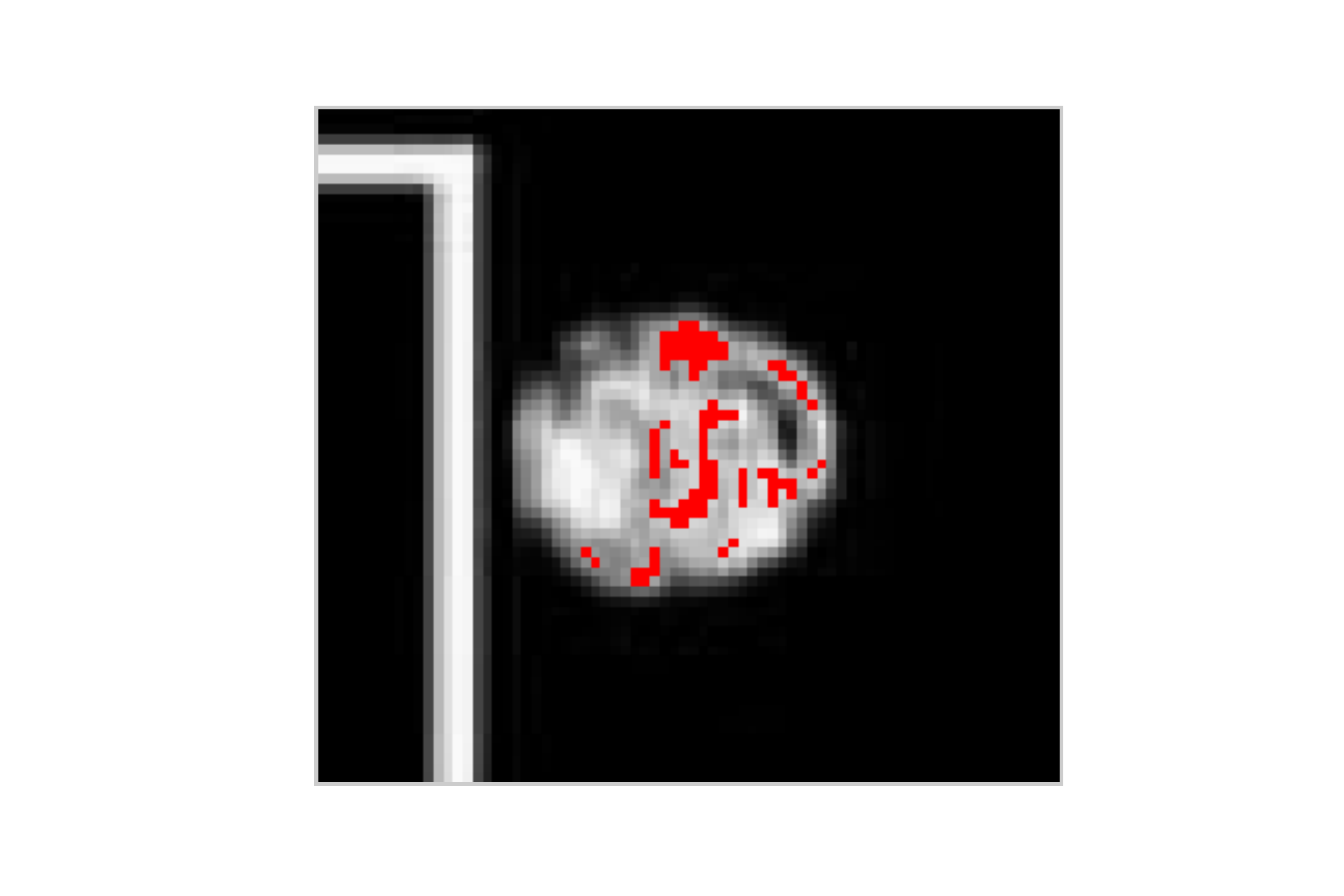}  
   \caption{}
 \end{subfigure}
\caption{(a) Areas corresponding to optimal weights from linear SVM classification. The areas of the PD that distinguish DMMP are shown in red and the ones that distinguish water are shown in blue. Inverse analysis based on SVM weights for (b) DMMP and (c) water responses. Note that the camera artifact in (c) has no highlighted areas, demonstrating that the extracted features are physically relevant.} \label{fig:invLC}
\end{figure}

Each optical micrograph is an image, \hl{taken from the endpoint of the sensor response}, with three channels (\textit{Red}, \textit{Green}, and \textit{Blue}).  We project these three channels onto a single grayscale channel by computing the total intensity of each pixel in the image. The conversion of the image to a single channel allows us to treat the image as a cubical complex over which we can perform a simple Morse filtration (level sets defined in terms of intensity, as done in the previous diffusion field example). A grayscale image was used to simplify the computations; however, more complicated approaches could be taken to deal with the three color channels. 
In order to understand the important of the information contained in the Morse filtration analysis, we apply a linear SVM to the vectorized persistence diagrams associated to our images. We find that the topological features of the images gives us a classification accuracy of $85 \pm 2\%$ (for a dataset of more than 1,000 images).  In order to identify the characteristic features for the responses at high and low concentration, we utilize the classification weights of the linear SVM model. The classification weights are visualized in Figure \ref{fig:invLC}. We apply a masking method to identify the portions of the persistence diagram that are critical for defining whether a response pattern is a result of high or low concentration. We can use these masked areas to identify the features of the images that separate the high concentration patterns from the low concentration patterns. To visualize the geometric differences between the patterns associated with DMMP and water, we again utilize the \texttt{Homcloud} software to perform the inverse analysis via volume optimal cycles. Here, we focus on inverse analysis for the $H_1$ homology group. The results of this analysis for a couple of sample images is found in Figure \ref{fig:invLC}. Inverse analysis reveals that, when the sensors are exposed to water, the pattern exhibits many small distinct clusters; in contrast, when the sensor is exposed to DMMP, there are few large clusters. This shows how inverse analysis allows us to pinpoint topological features of the image that drive classification. \hl{The ability to classify as well as extract meaningful information from the topology of the images provides an important advantage over machine learning methods. The work of Smith, Cao, and co-workers demonstrate that machine learning models can be used to classify the responses of these sensors with higher accuracy, but provide minimal interpretability. Interpretability is needed to understand the physics of these sensors which can improve the sensitivity, and broaden the applicability, of the sensors \cite{smith2020convolutional,cao2018machine}}.

\subsection{Topology of Probability Density Functions}
\label{sec:ex5}

\begin{figure}[!htp]
 \begin{subfigure}{.24\textwidth}
   \centering
   \includegraphics[width=1\linewidth]{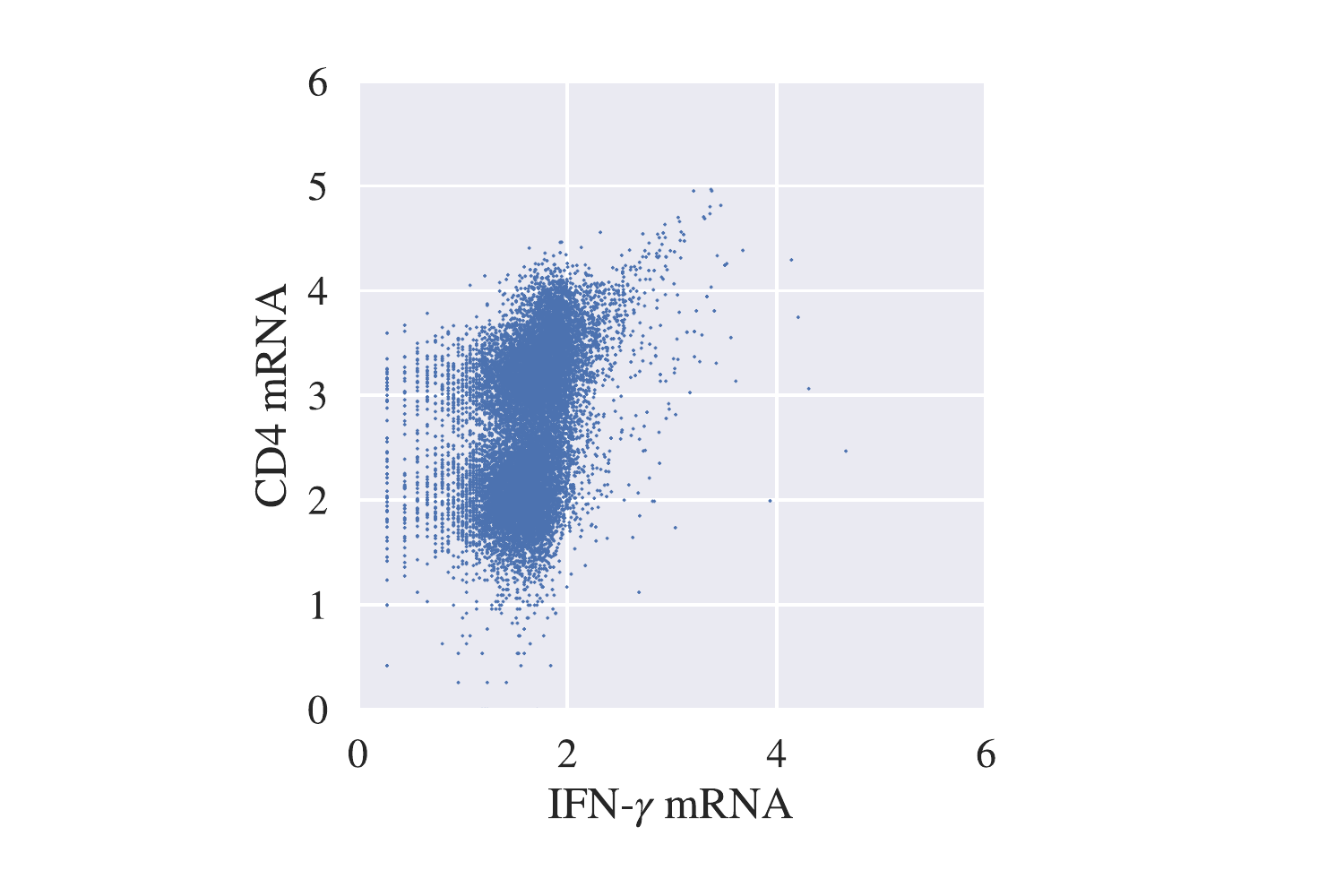}  
   \caption{Time = 0 Minutes}
 \end{subfigure}
 \begin{subfigure}{.24\textwidth}
   \centering
   \includegraphics[width=1\linewidth]{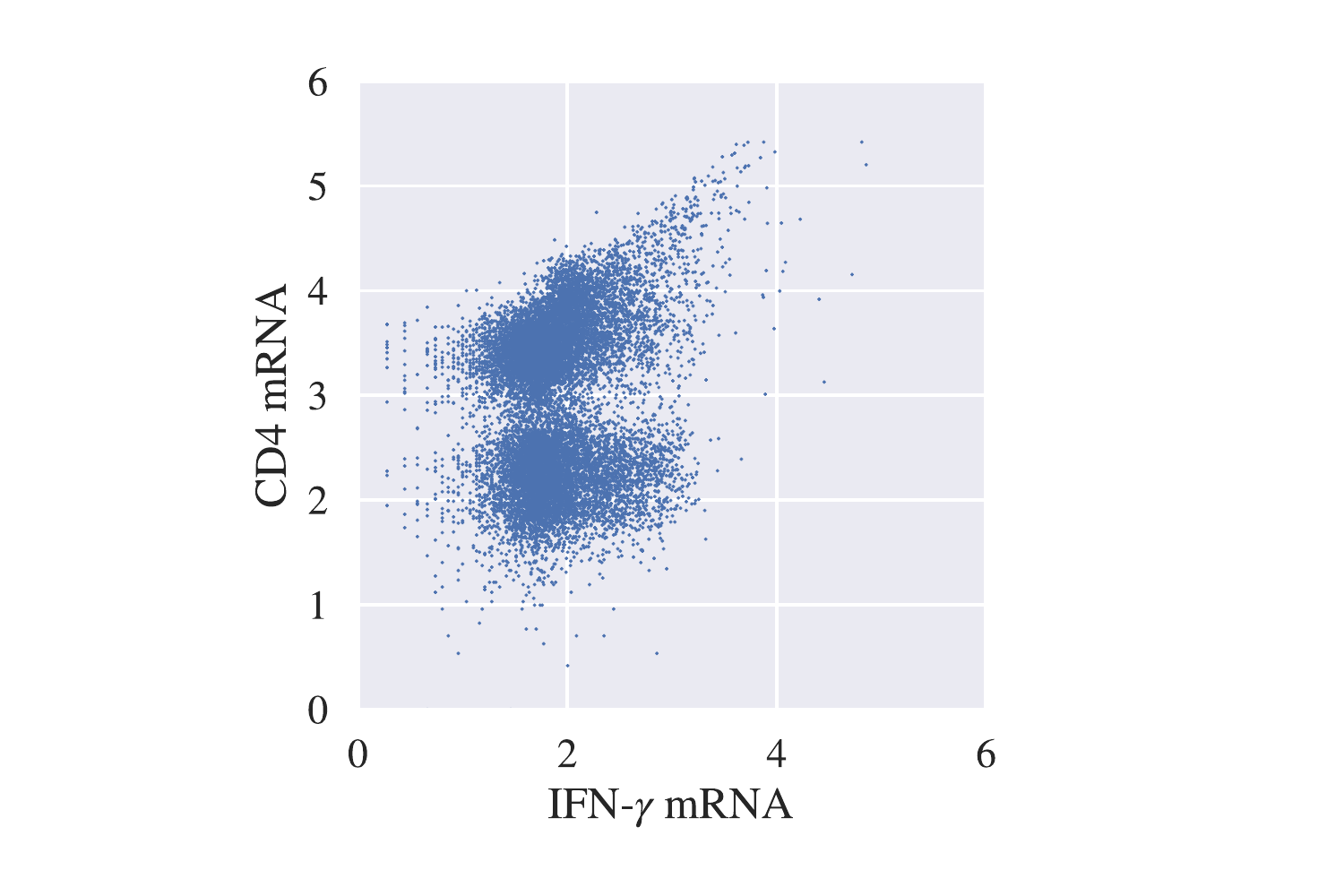}  
   \caption{Time = 30 Minutes}
 \end{subfigure}
 \begin{subfigure}{.24\textwidth}
   \centering
   \includegraphics[width=1\linewidth]{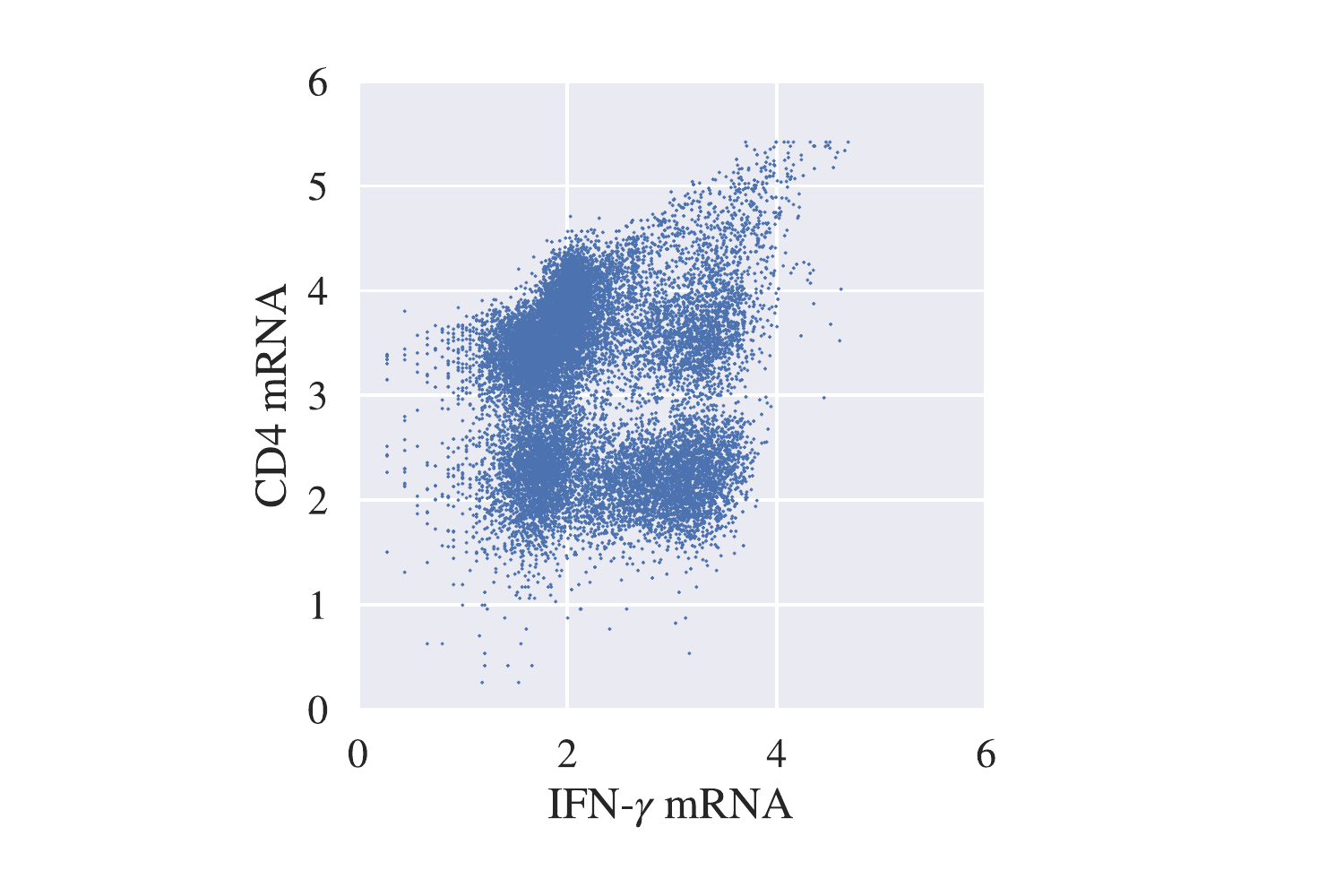}  
   \caption{Time = 60 Minutes}
 \end{subfigure}
 \begin{subfigure}{.24\textwidth}
   \centering
   \includegraphics[width=1\linewidth]{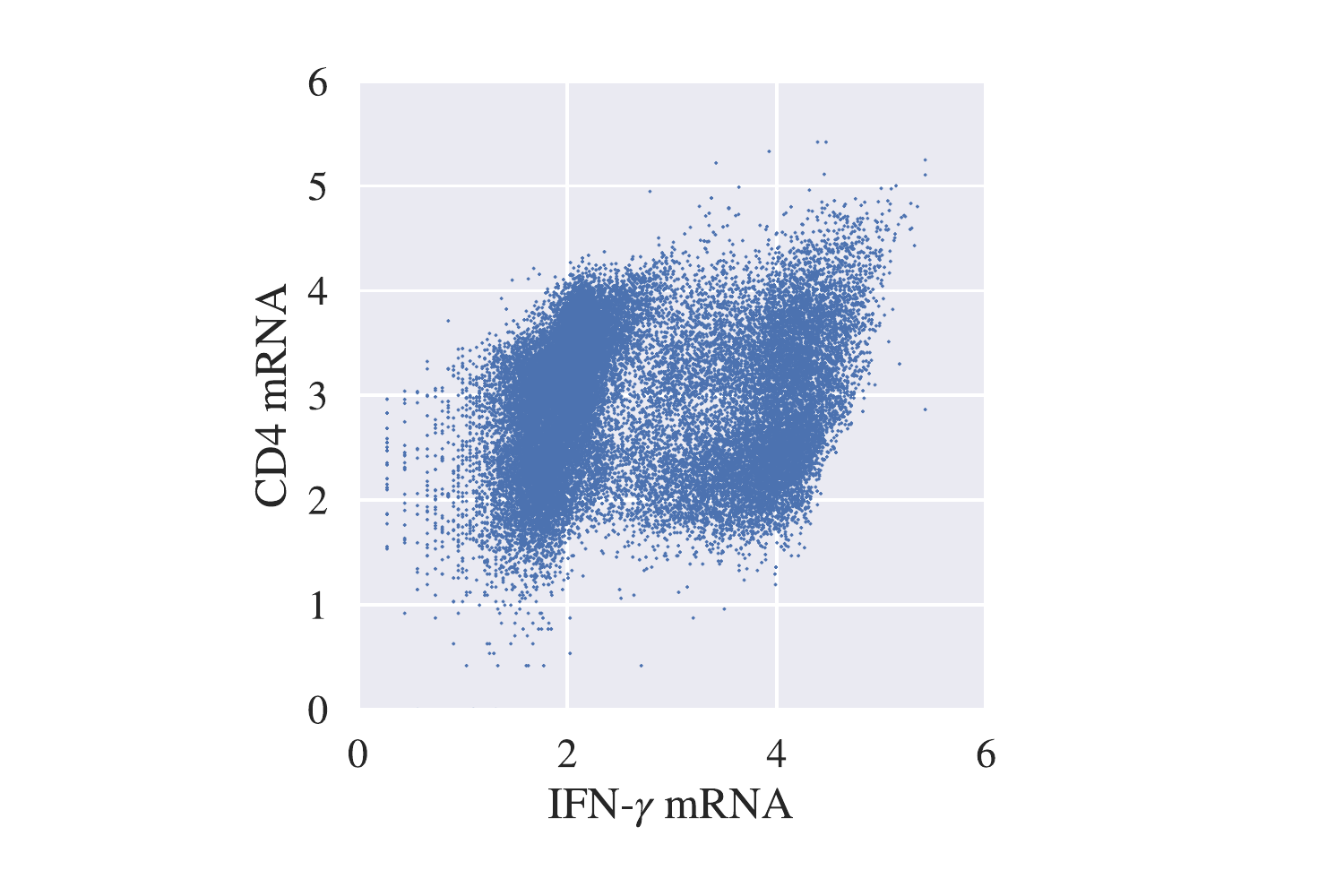}  
   \caption{Time = 90 Minutes}
 \end{subfigure}
\caption{Deformation of a $2$D scatter field over time.}
\label{fig:cytodata}
\end{figure}

The analysis of the shape of probability density functions is typically done using summarizing statistics (e.g., moments such as skewness and kurtosis) or through parametric techniques (fitting a parametric model such as a Gaussian mixture to the data) \cite{walter1997identification}. These models are powerful in their simplicity but might not be flexible enough to capture complex features of density functions (particularly in high dimensions). In this example we explore the shape of complex density functions by using topological techniques. We use an experimental flow cytometry dataset to illustrate how this can be done. The flow cytometry dataset was obtained through the FlowRepository \cite{flowrepository} (Repository ID: FR-FCM-ZZC9). This dataset represents a temporal study of the kinetics of gene transcription and \hl{protein} translation within stimulated human blood mononuclear cells through the quantification of \hl{proteins} (CD4 and IFN-$\gamma$) and mRNA (CD4 and IFN-$\gamma$) \cite{van2014simultaneous}. In our study, we focus on the evolution of the concentration of CD4 mRNA and IFN-$\gamma$ mRNA in a given cell which is measured via a flow cytometer. At each time point in the study, a number of cells ($\sim 15,000$)  are passed through the flow cytometer, each one of these cells provides a vector of scalar values $x \in \mathbb{R}^n$ corresponding to each measurable variable. In this case we set $n=2$ as we are only utilizing the scalar values that represent the measure of both CD4 mRNA and IFN-$\gamma$ mRNA, samples of these distributions are found in Figure \ref{fig:cytodata}. From this we obtain a $2$D scatter field. We have restricted ourselves to two dimensions for illustrative purposes, but this same analysis can be conducted for a point cloud of higher dimensions, accounting for all variables measured by a flow cytometer. The goal in this analysis is to use TDA to quantify the temporal evolution of the shape of the scatter field during the kinetic response of human blood mononuclear cells during stimulation. \hl{This approach provides an alternative to traditional parametric or heuristic methods such as gating, which are difficult to tune as they are highly sensitive to potential noise and outliers in the data. The gate selection may also require complicated multivariate mixture models to identify the correct gating values \cite{lo2008automated}.} 

 \begin{figure}[!htp]
 \begin{subfigure}{.33\textwidth}
   \centering
   \includegraphics[width=.9\linewidth]{FC_Scattert120.pdf}  
   \caption{Sample Distribution.}
 \end{subfigure}
 \begin{subfigure}{.33\textwidth}
   \centering
   \includegraphics[width=.9\linewidth]{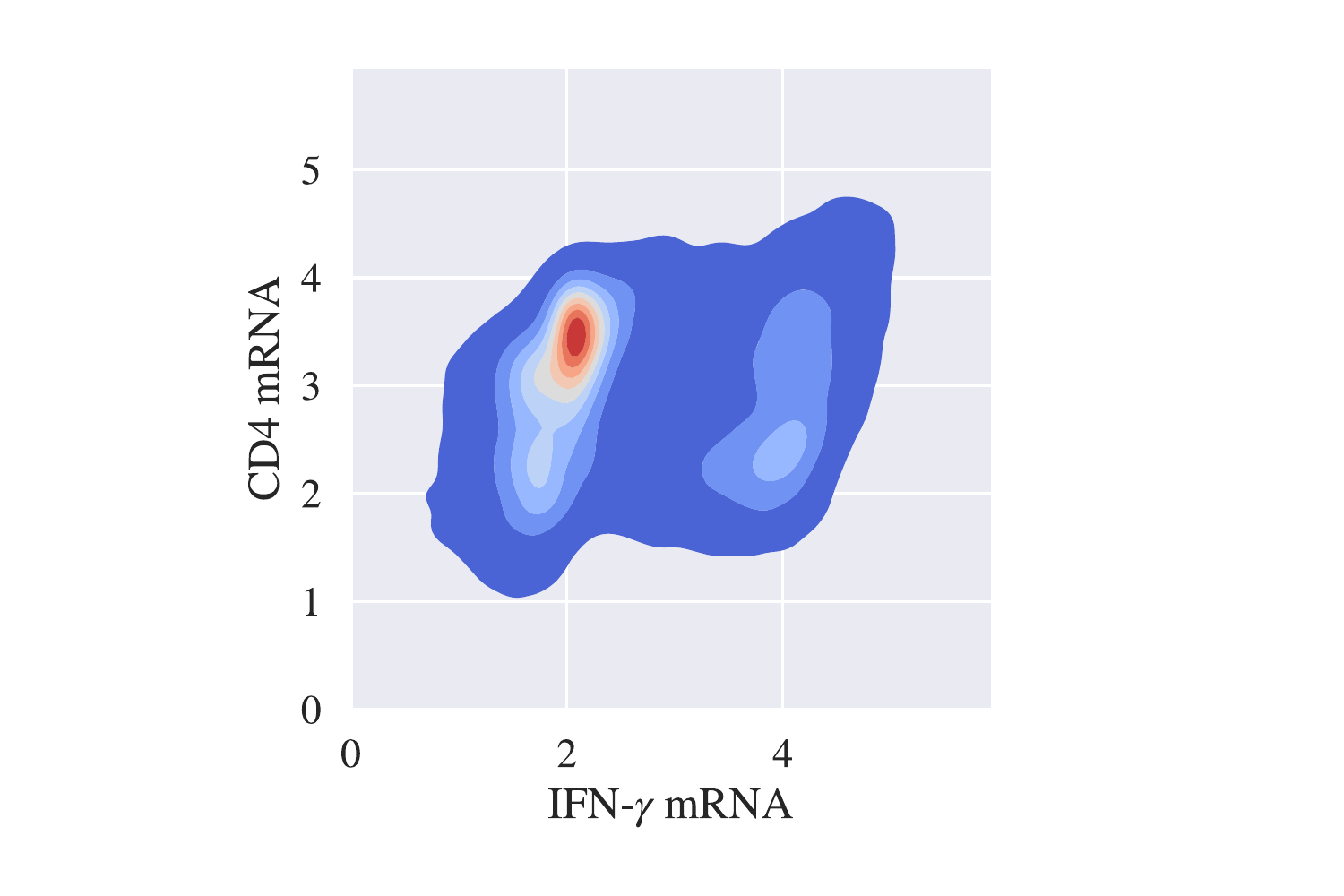}  
   \caption{Gaussian Density Estimate}
 \end{subfigure}
 \begin{subfigure}{.33\textwidth}
   \centering
   \includegraphics[width=.9\linewidth]{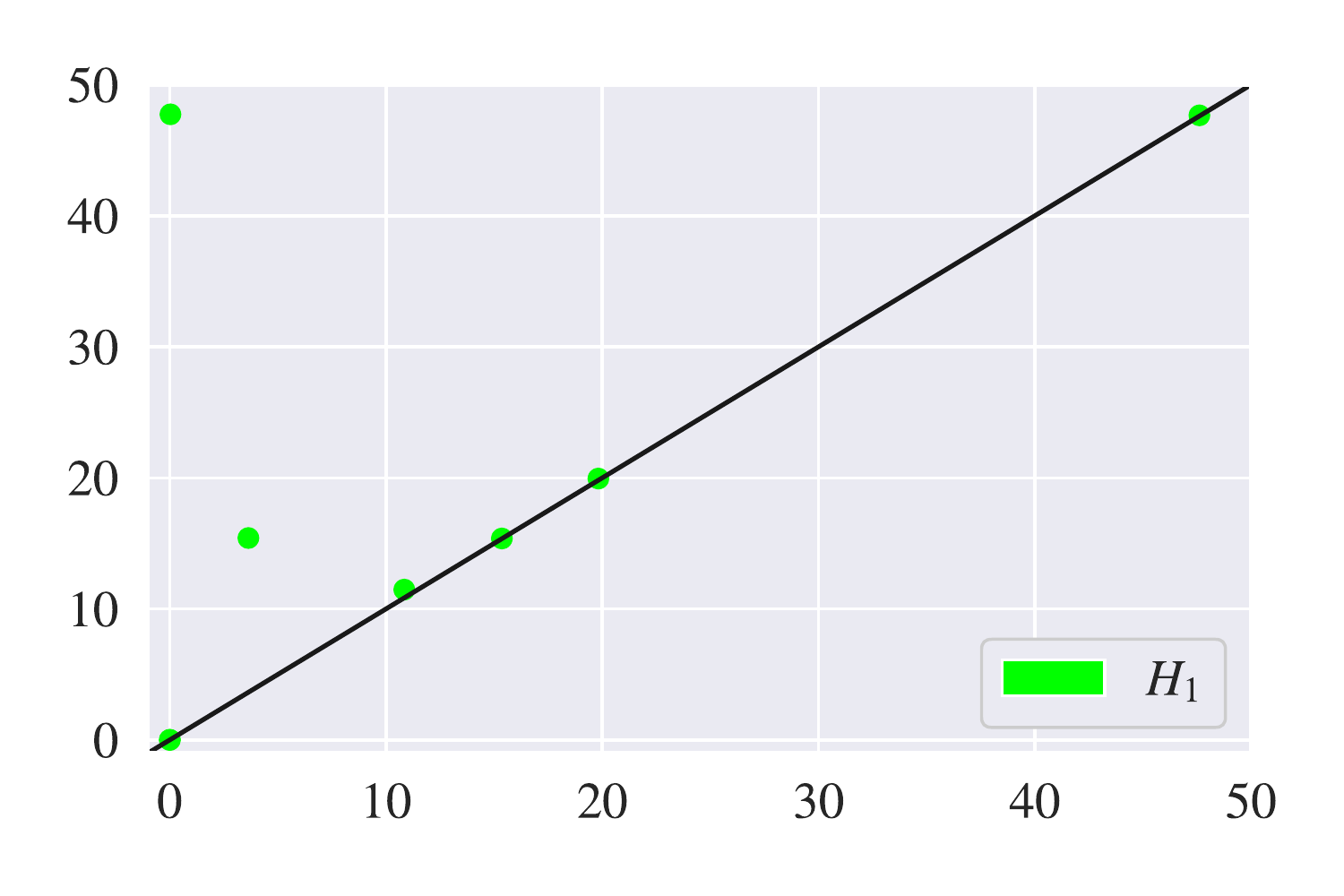}  
   \caption{Persistence Diagram}
 \end{subfigure}
\caption{Processing of the flow cytometry scatter field. The raw data is first smoothed via a Gaussian kernel and then the smoothed diagram is processed via a Morse filtration.}
\label{fig:cytodata2}
\end{figure}

 \begin{figure}[!htp]
 \begin{subfigure}{.5\textwidth}
   \centering
   \includegraphics[width=.9\linewidth]{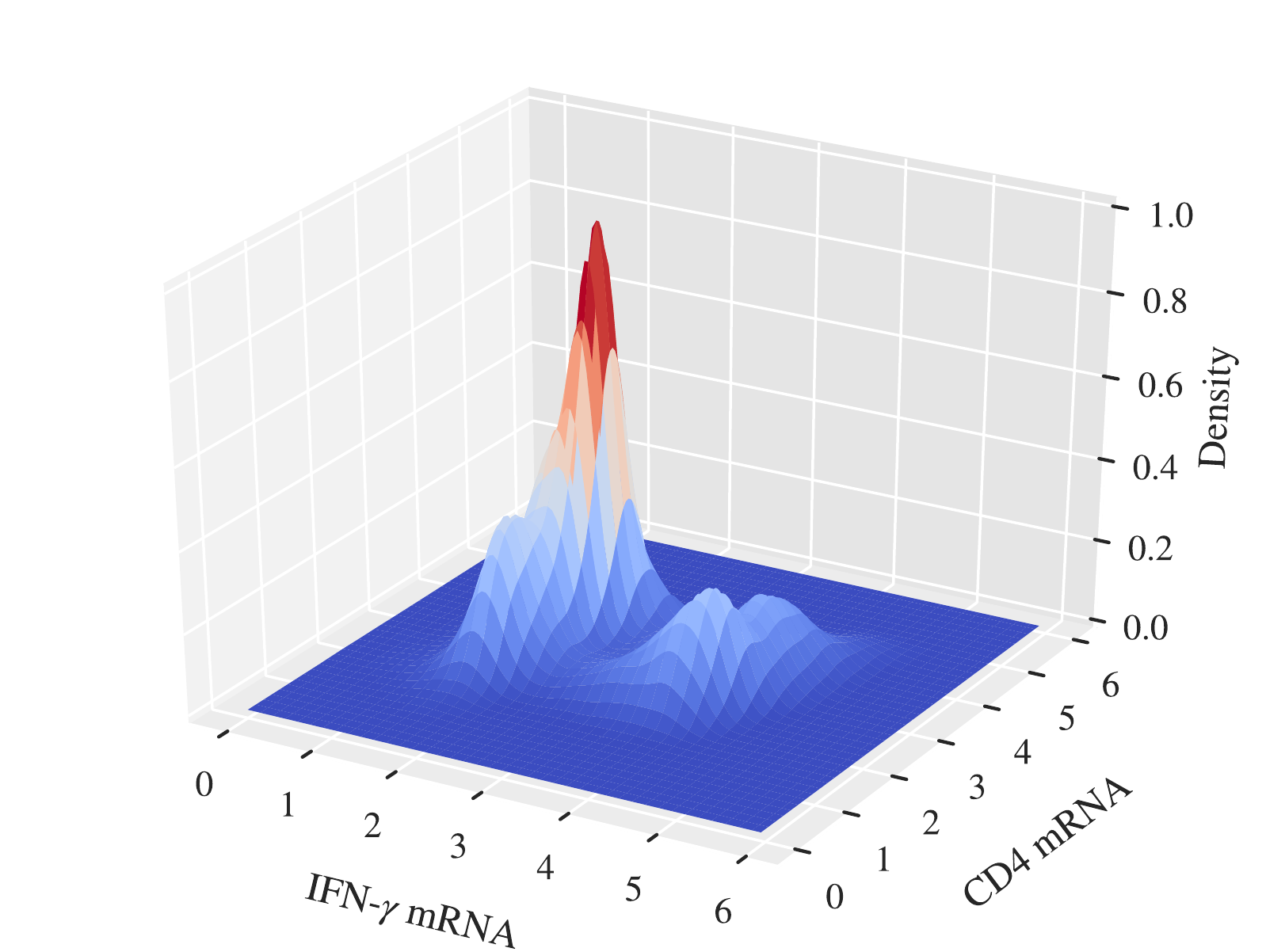}  
   \caption{}
   \label{fig:sub-first2}
 \end{subfigure}
 \begin{subfigure}{.5\textwidth}
   \centering
   \includegraphics[width=.8\linewidth]{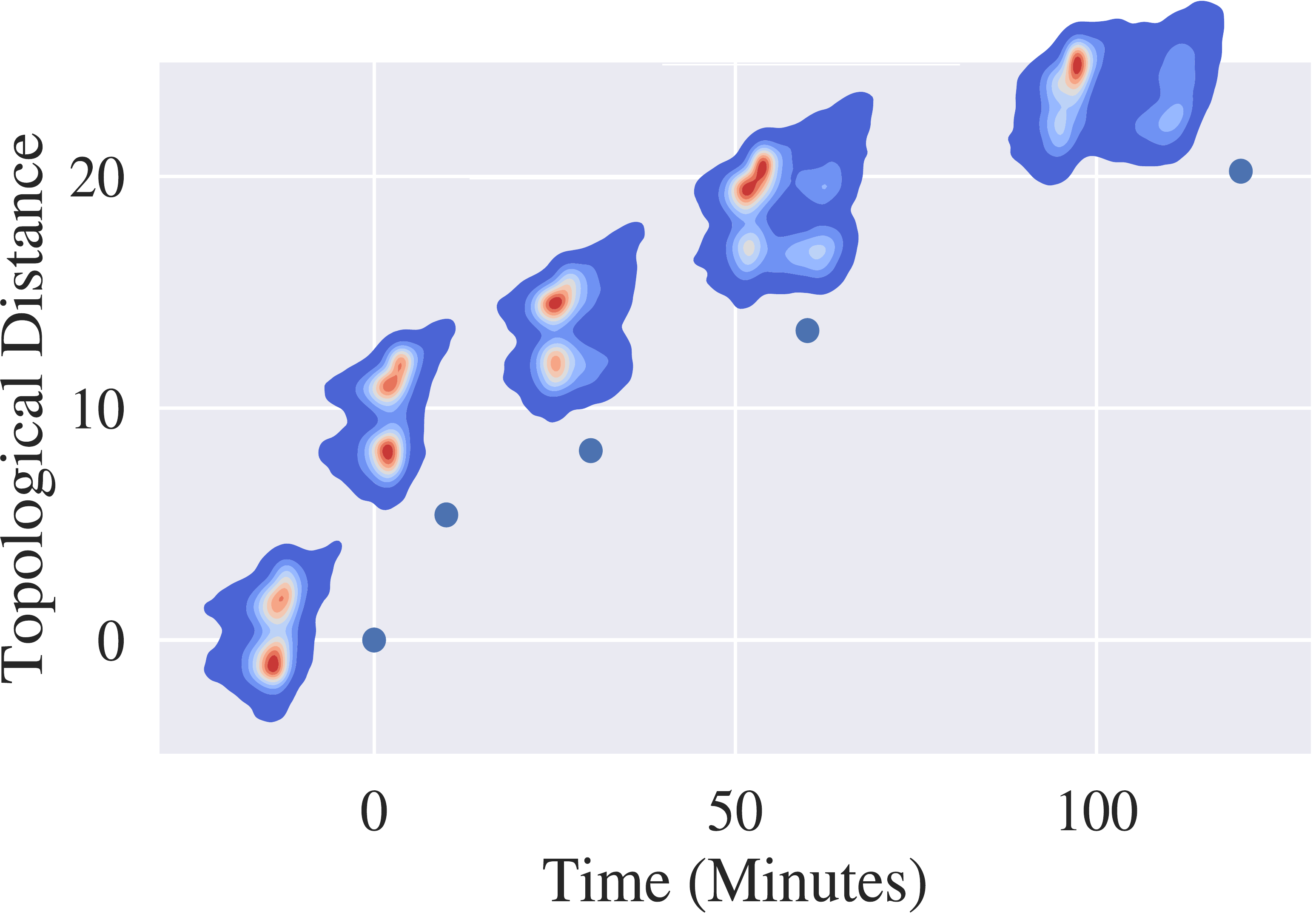}  
   \caption{}
   \label{fig:Endo_Plot2}
 \end{subfigure}
\caption{(a) Scatter field as a probability density function in 3D; Morse filtration is performed using level sets of the probability density. (b) Wasserstein distance between the persistence diagrams as a function of time. There is a clear continuous evolution of the Wasserstein distance that characterizes the change in topology.}
\label{fig:invLC2}
\end{figure}

In order to analyze the topology of the scatter fields, we utilize Gaussian kernel smoothing \cite{sheather2004density}. The work of \cite{bobrowski2017topological} demonstrates that provided a large enough sample, the homology of the Gaussian kernel density estimate derived from a sample is equivalent to the homology of the true density. Figure \ref{fig:cytodata2} shows the Gaussian kernel smoothing of a flow cytometry scatter plot. An example persistence diagram for the function is shown in Figure \ref{fig:cytodata2}. 
As in the case of the diffusion example, in 3D we can represent the scatter field as a continuous function (in this case a probability density function). This probability density function (Figure \ref{fig:sub-first2}) can be analyzed using Morse filtration.  Our goal in this analysis is to quantify the time evolution of the probability density functions during stimulation. Our strategy consists of computing the Wasserstein distance ($w_d$) between the persistence diagram of a given time point to the persistence diagram of the sample at time zero \cite{carriere2017sliced}.  Specifically,  given flow cytometry samples $X_{1}$, $X_{2}$ and the time zero sample $X_{0}$ as well as their corresponding persistence diagrams $PD_1,PD_2,PD_{0}$ and time points $T_0 \leq T_1 \leq T_2$, we observe that: 

\begin{equation}
    w_d(PD_1,PD_{0}) \leq w_d(PD_2,PD_{0})
\end{equation}
From Figure \ref{fig:Endo_Plot2}, we can see that the distance exhibits a  strong dependence on time. This suggests that there exists a continuous mapping between the persistent diagram and time (the topological deformation is continuous with respect to time). This again reveals the continuity of the persistence diagrams (of topology) to perturbations.  These results highlight how topological data analysis provides a quantifiable approach to characterize complex probability density functions and their evolution over time.

 \subsection{Topology of 3D Fields}
 \label{sec:ex6}

\begin{figure}[!htp]
    \centering
    \includegraphics[width = .5\textwidth]{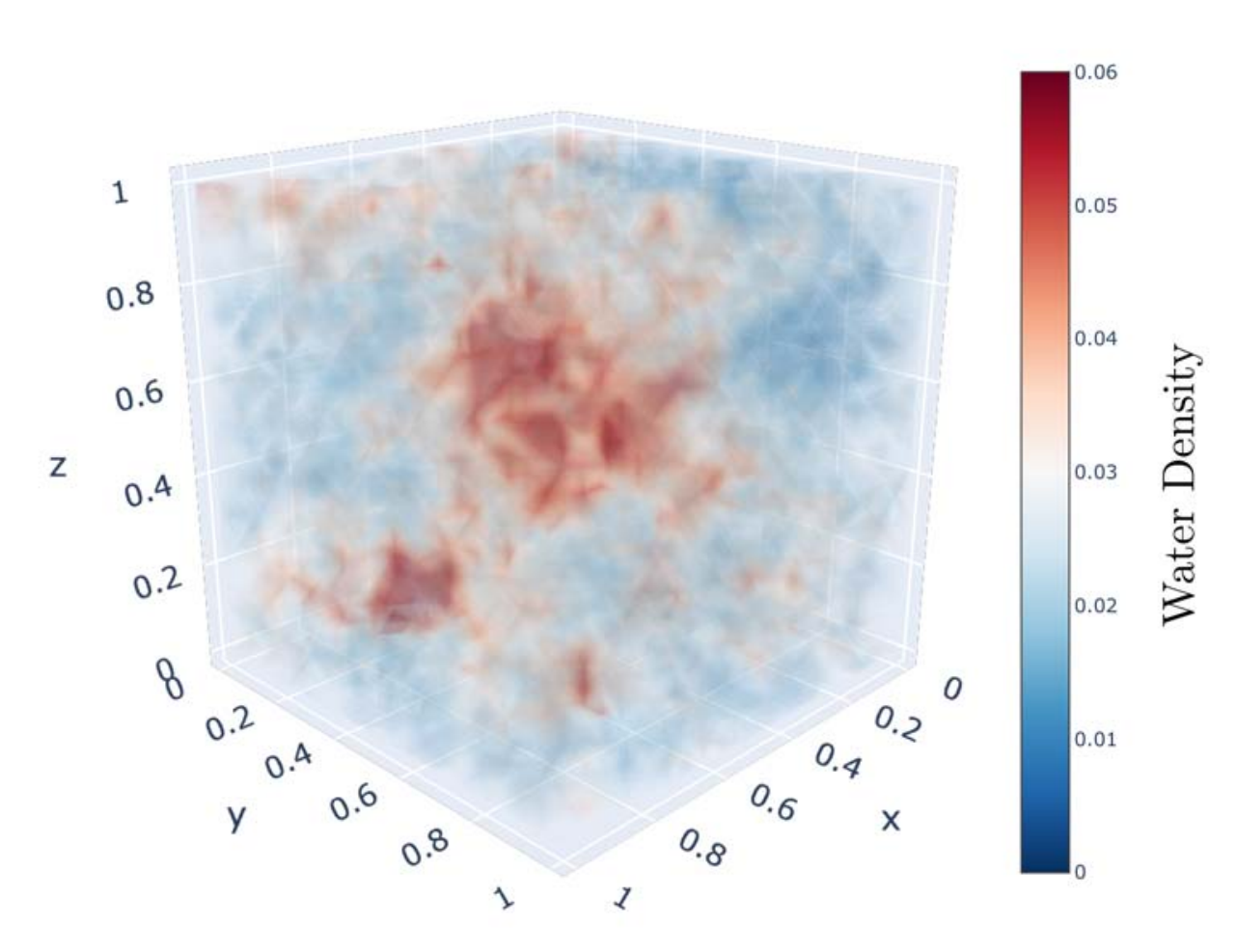}
    \caption{Visualization of 3D water density field generated by MD simulation.}
    \label{fig:dens_plot}
\end{figure}

We now illustrate how to use TDA to analyze the topology induced by 3D point clouds. Specifically, we study datsets generated by molecular dynamics (MD) simulations  \cite{chew2019fast,walker2018universal}. The dataset under study analyzes the influence of the 3D liquid-phase environment formed by molecules of a solvent, co-solvent, and a reactant on  reactivity \cite{chew2019fast}. The reactivity is quantified via a \textit{kinetic solvent parameter} $\sigma$ that is obtained from experiments. Experiments suggest that reactivity is influenced by hydrophilicity of the solvent. The main hypothesis is that, as the solvent concentration is increased, the water in the system is concentrated around the solvent molecule, and that molecules with high hydrophilicity are able to take advantage of this effect. In order to study this hypothesis, \hl{molecular} dynamics computations were performed in \cite{chew2019fast}.  The data output of an MD simulation has both spatial and temporal dimensions. Each simulation gives atomic positions $X_t \in \R^{M\times3}$ ($M$ is the number of species) at multiple times $t$ (measured in nanoseconds). In our analysis, we utilize a 3D point cloud of water molecule positions that result from a time average of 100 nanoseconds. An example of this point cloud (visualized as a field) is provided in Figure \ref{fig:dens_plot}. Each density field is labeled with a reactivity $\sigma$ obtained from experiments.  Recent work by Chew and co-workers has analyzed the 3D density by using 3D convolutional neural networks (CNNs) and has shown that the features extracted from the CNN are strongly correlated to reactivity \cite{chew2019fast}. \hl{CNNs are highly effective tools but require a large number of parameters ($\sim$160,000 in this case) and are difficult to interpret. 
Our goal is to study if the topology of the 3D density can be characterized in a more straightforward manner and to explore whether topology changes in correlation with reactivity.}

 \begin{figure}[!htp]
 \begin{subfigure}{.33\textwidth}
   \centering
   \includegraphics[width=.9\linewidth]{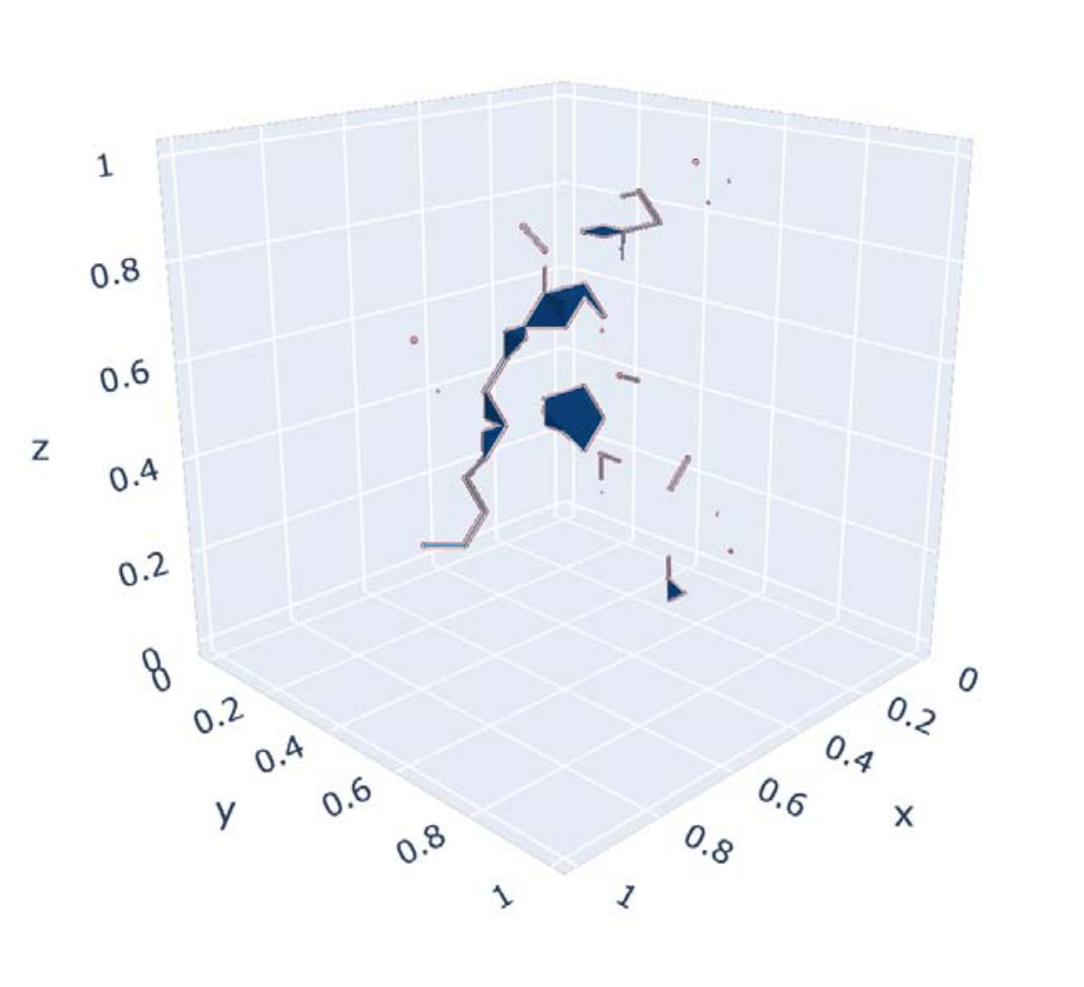}  
   \caption{Density Level Set of 0.01}
   \label{fig:sub-first}
 \end{subfigure}
 \begin{subfigure}{.33\textwidth}
   \centering
   \includegraphics[width=.9\linewidth]{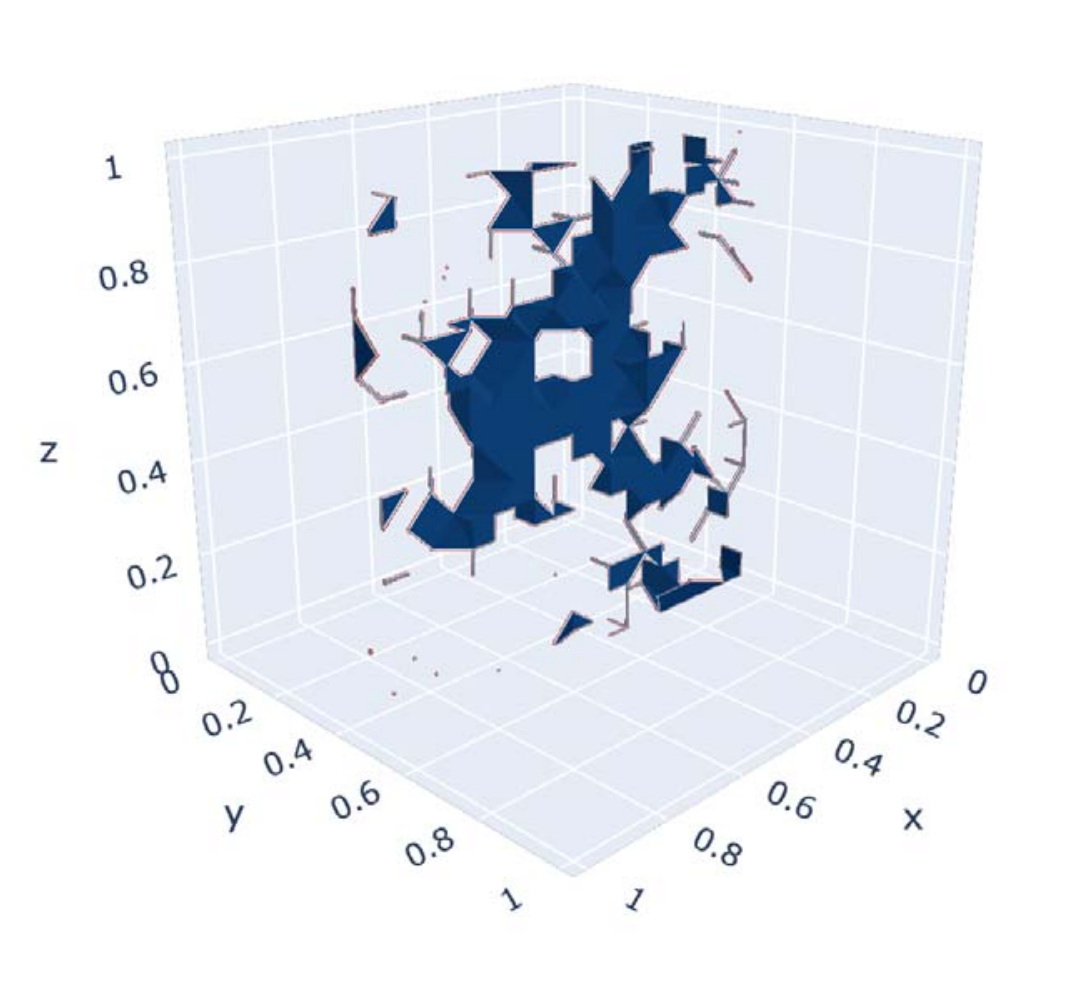}  
   \caption{Density Level Set of 0.02}
   \label{fig:sub-second}
 \end{subfigure}
 \begin{subfigure}{.33\textwidth}
   \centering
   \includegraphics[width=.9\linewidth]{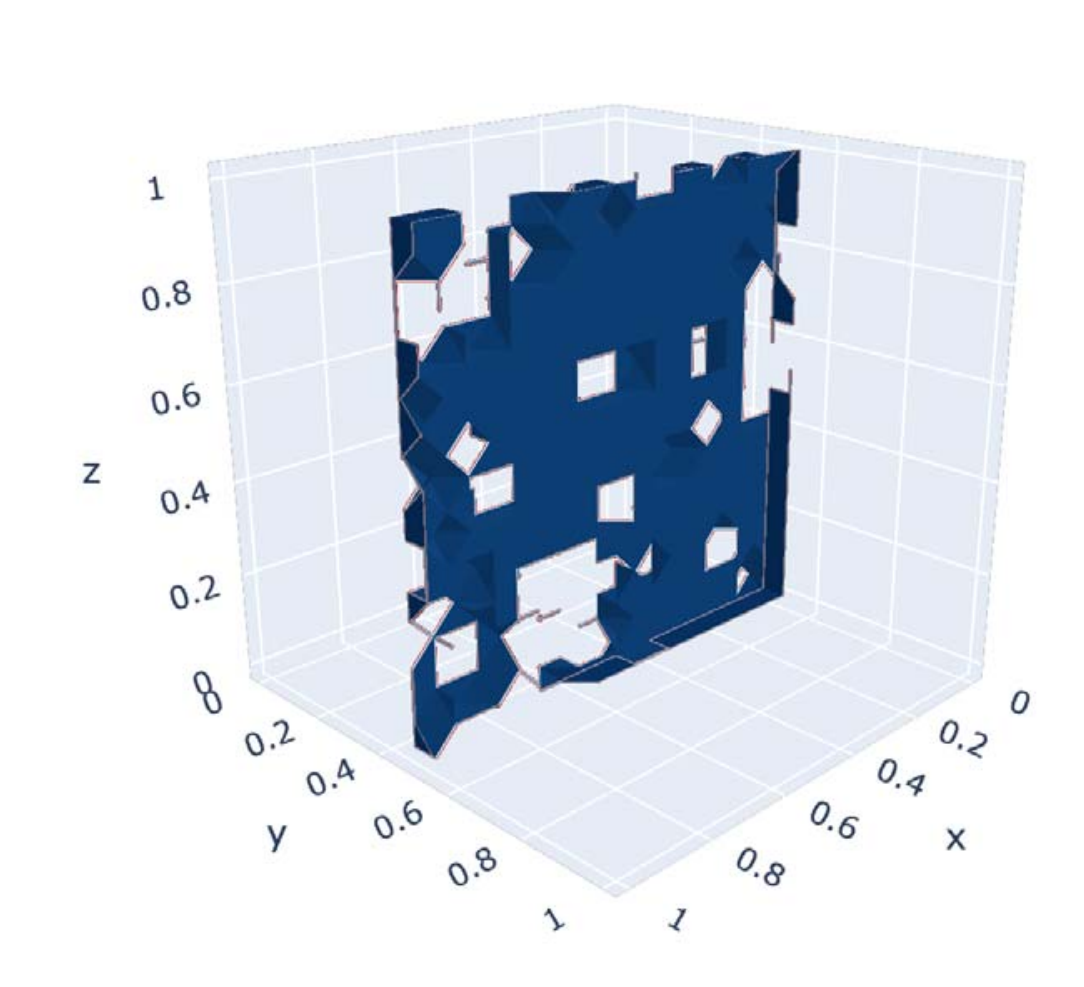}  
   \caption{Density Level Set of 0.03}
   \label{fig:sub-second}
 \end{subfigure}
\caption{Slices of 3D water density field as filtration proceeds for different density values. The filtration reveals the presence of voids in the data associated with high concentrations of water molecules.}
\label{fig:slices}
\end{figure}

 \begin{figure}[!htp]
 \begin{subfigure}{.5\textwidth}
   \centering
   \includegraphics[width=.7\linewidth]{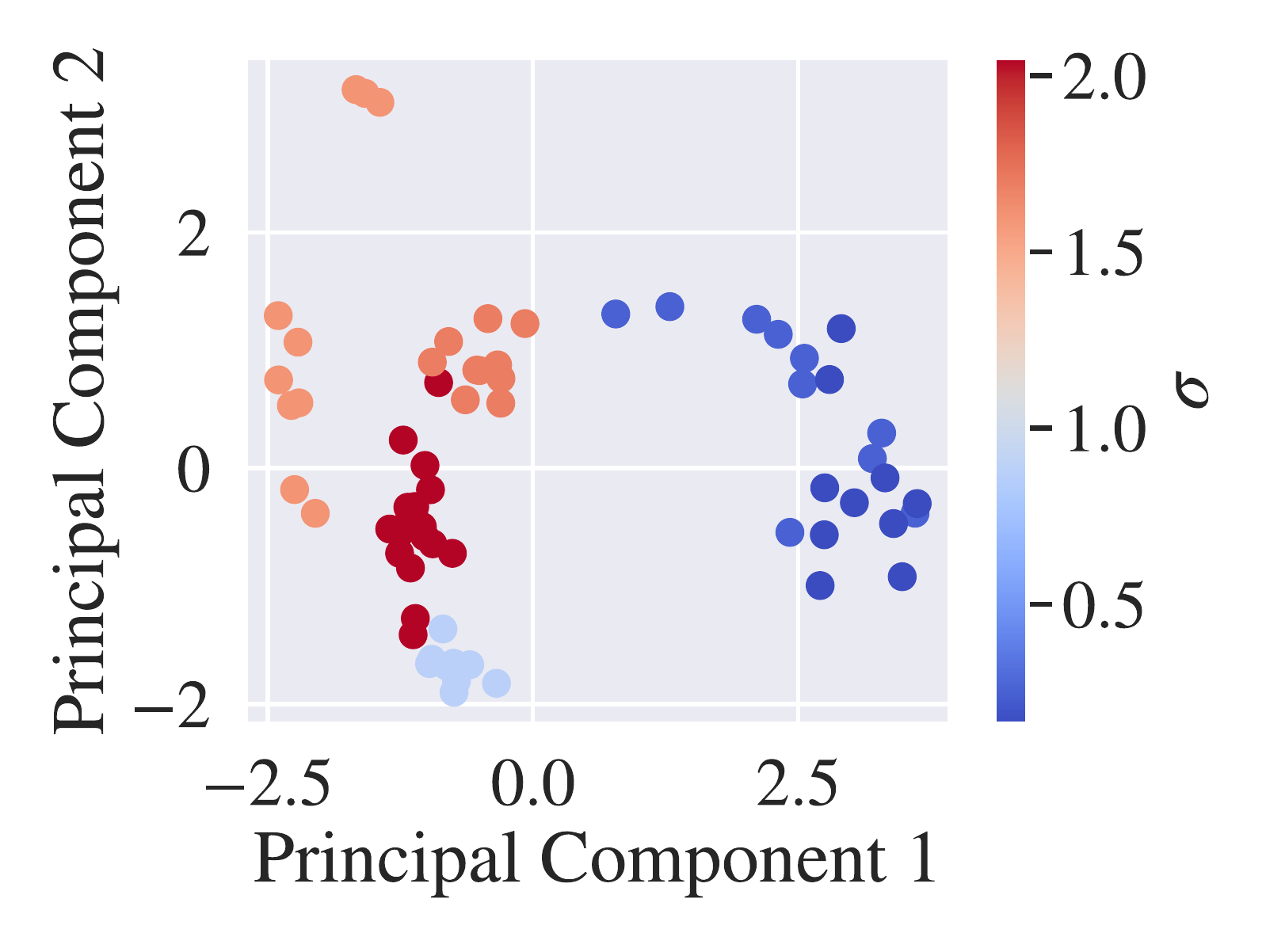}  
   \caption{PCA analysis}
   \label{fig:pca_md}
 \end{subfigure}
 \begin{subfigure}{.5\textwidth}
   \centering
   \includegraphics[width=.7\linewidth]{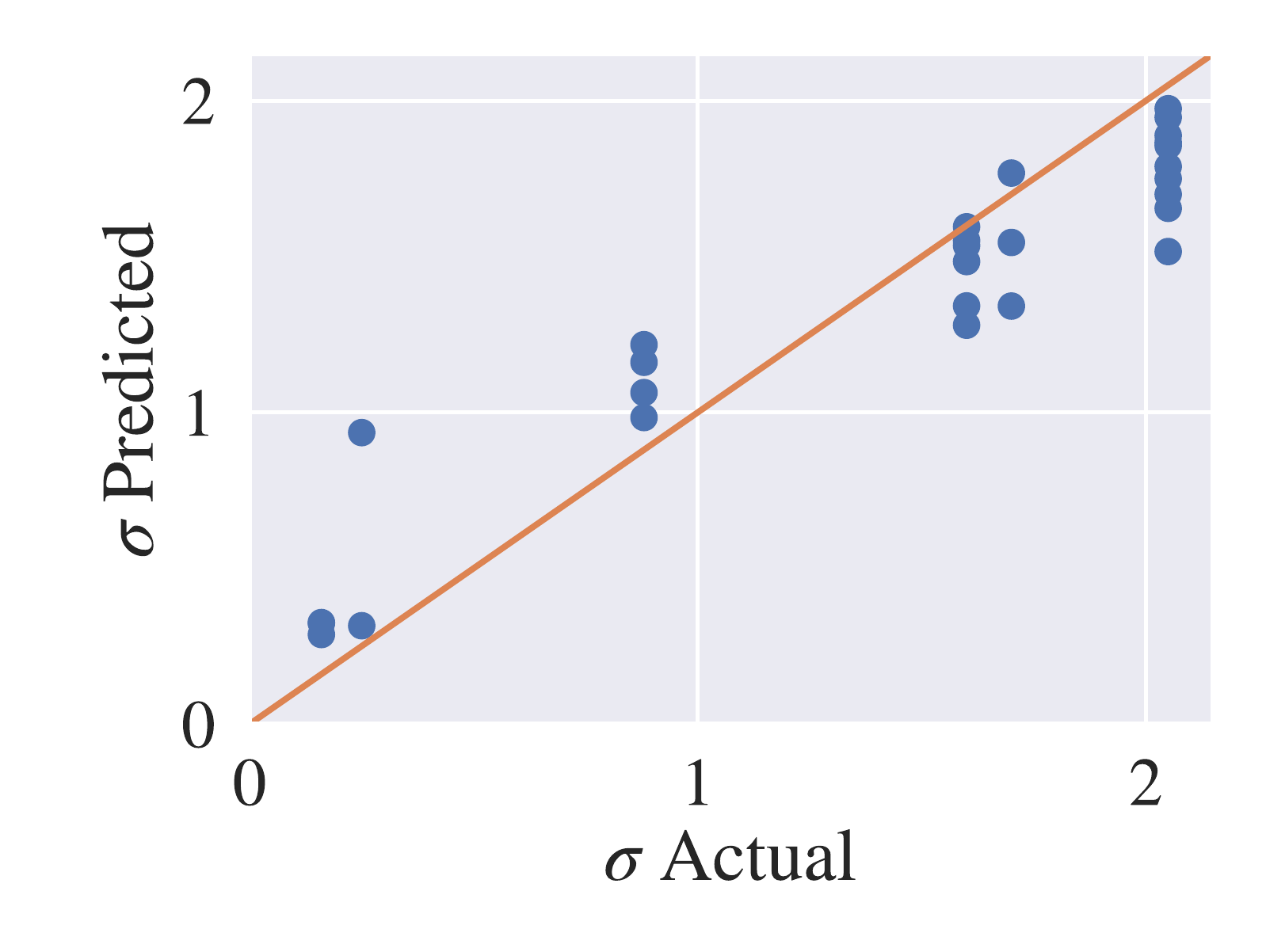}  
   \caption{SVM Regression}
   \label{fig:reg_md}
 \end{subfigure}
\caption{(a) PCA analysis on persistence diagrams for MD simulations. The analysis reveals strong dependence of reactivity $\sigma$. (b) Regression plot for SVM with a radial basis function. The predictions over 5-fold cross validation yield an MSE of $0.07 \pm 0.003$.}
\label{fig:md_analysis}
\end{figure}

To perform our analysis, we treat the 3D point cloud as a continuous field (function). Here, we perform a Morse filtration and treat the data as a 3D cubical complex (a voxel). The filtration will be done by exploring levels sets for the water density. Note that this is a filtration in a higher dimension that our previous examples for the diffusion field and for liquid crystal sensors. The main focus of this approach is to capture the clustering of water near hydrophilic molecules, and the lack of clustering near non-hydrophilic molecules. We visualize the water density filtration in Figure \ref{fig:slices} via a $2$-dimensional slice. We can see that \hl{voids} in the data are generated as we increase the filtration value. The voids represent areas of high water density which is precisely what we wish to quantify. From these filtrations, we produce persistence diagrams focused on $H_2$ (since this homology group quantifies these voids). The persistent  diagrams are then vectorized via the persistence image and we use PCA reduction to visualize them (see Figure \ref{fig:pca_md}). We use SVM regression with a radial basis kernel to predict reactivity as a function of the persistence diagrams (see Figure \ref{fig:reg_md}). 

The PCA projection reveals that there is \hl{strong dependence of the reactivity on the topology}. This suggests that the information gained via persistent homology extracts informative features of the 3D field that explain reactivity. This also suggests that a simple regression method (as opposed to a complex neural network) would be effective at predicting the reactivity. In order to test this hypothesis, the experimental dataset with $70$ points is split into train/test sets with 49/21, points to build the SVM regression model. A 5-fold cross validation is performed to estimate the performance of the SVM regression. We found that this model captures the trend of reactivity well; moreover, this simple model yields a mean square error of $0.07 \pm .003$ (Figure \ref{fig:md_analysis}). These results are relevant because they indicate that it is possible to predict experimental reactivity directly from MD simulations. 

\section{Conclusions}

Topological data analysis (TDA) provides a set of powerful methods and tools for understanding the underlying geometry of data. These techniques represent data such, as point clouds and functions, as geometric objects and explores these objects in terms of basic geometric features. We have shown that TDA offers a number of important theoretical \hl{properties} (such as stability), offers flexibility to extract features from different types of data, and how these extracted features can be exploited using statistical and machine learning techniques. \hl{In particular, we demonstrate the capability of TDA in characterizing the shape of 2-dimensional point clouds and distributions without parameter estimation (Sections \ref{sec:ex1},\ref{sec:ex5}). We also utilize TDA in understanding the phase space of coupled time series, and how it can be used to detect anomalies based upon the geometry of the time series alone (Section \ref{sec:ex2}). The TDA methodology also provides a new perspective on the analysis of continuous data, such as an image, as a manifold embedded in higher dimensional space. This perspective allows for the analysis of the geometric and topological structure of data such as 2-dimensional images (Section \ref{sec:ex3},\ref{sec:ex4}) or 3-dimensional density fields from MD simulations (Section \ref{sec:ex6}}). The TDA field is in rapid development from both a theoretical and applied perspective. Moreover, there are easily accessible software packages to conduct scalable computations (such as such as \texttt{Gudhi} \cite{maria2014gudhi}, \texttt{Homcloud} \cite{homcloud}, \texttt{Javaplex-Matlab} \cite{adams2014javaplex} and \texttt{R-TDA} \cite{fasy2014introduction}). We believe that TDA can bring a new geometric and topological perspective to challenging problems in chemical engineering.

\section*{Acknowledgements}

We acknowledge funding from the U.S. National Science Foundation  (NSF) under BIGDATA grant IIS- 1837812 and also acknowledge partial support from the NSF through the University of Wisconsin Materials Research Science and Engineering Center (DMR-1720415). \hlp{We acknowledge the Dioscouri Center for Topological Data Analysis for providing expert advice on Topological Data Analysis.}  We thank Alex Chew and Prof. Reid Van Lehn for providing datasets on MD simulations and thank Prof. Nicholas Abbott for providing the liquid crystal dataset.

\bibliography{References}

\end{document}